\documentclass[review,onefignum,onetabnum]{siamart190516}

\usepackage[toc,page]{appendix}
\usepackage[utf8]{inputenc}
\usepackage{amsfonts}
\usepackage{graphicx}
\usepackage{amsmath}
\usepackage[normalem]{ulem}
\usepackage{comment}
\usepackage{dsfont}
\usepackage{subfig}
\usepackage{mathrsfs}
\usepackage{float}
\usepackage{amssymb}
\usepackage{amsopn}
\usepackage{bbm}
\usepackage{xcolor}
\usepackage{cleveref}
\usepackage{cancel}

\reversemarginpar
\usepackage{todonotes}

\graphicspath{{Fig/}}
\newcommand{\T}{\mathcal{T}}
\newcommand{\Pp}{\mathscr{P}}
\newcommand{\R}{\mathbb{R}}

\usepackage{soul} 

\newsiamremark{remark}{Remark}
\newsiamremark{hypothesis}{Hypothesis}
\crefname{hypothesis}{Hypothesis}{Hypotheses}
\newsiamthm{claim}{Claim}
\newsiamthm{example}{Example} %
\newsiamthm{assumptions}{Assumptions} %

\ifpdf
  \DeclareGraphicsExtensions{.eps,.pdf,.png,.jpg}
\else
  \DeclareGraphicsExtensions{.eps}
\fi

\newcommand{\TheShortTitle}{Optimal Transport for Parameter Identification}
\newcommand{\TheAuthors}{Y. Yang, L. Nurbekyan, E. Negrini, R. Martin, M. Pasha}

\headers{\TheShortTitle}{\TheAuthors}

\title{Optimal Transport for Parameter Identification of Chaotic Dynamics via Invariant Measures}

\author{Yunan Yang\thanks{Institute for Theoretical Studies, ETH Z\"urich, Z\"urich, Switzerland 8092. (\email{yunan.yang@eth-its.ethz.ch})}
\and Levon Nurbekyan\thanks{Department of Mathematics, UCLA, Los Angeles, CA 90095. (\email{lnurbek@math.ucla.edu})}
\and Elisa Negrini\thanks{Department of Mathematical Sciences, Worcester Polytechnic Institute, Worcester, MA 01609. (\email{enegrini@wpi.edu})}
\and Robert Martin\thanks{Air Force Research Laboratory, Edwards AFB, CA 93524. (\email{robert.martin.101@us.af.mil})}
\and Mirjeta Pasha\thanks{School of Mathematical and Statistical Sciences, Arizona State University, Tempe, AZ 85281. (\email{mpasha3@asu.edu})
} 
}

\ifpdf
\hypersetup{
  pdftitle={\TheShortTitle},
  pdfauthor={\TheAuthors}
}
\fi

\begin{document}

\maketitle

\begin{abstract}
We study an optimal transportation approach for recovering parameters in dynamical systems with a single smoothly varying attractor. We assume that the data is not sufficient for estimating time derivatives of state variables but enough to approximate the long-time behavior of the system through an approximation of its physical measure. Thus, we fit physical measures by taking the Wasserstein distance from optimal transportation as a misfit function between two probability distributions. In particular, we analyze the regularity of the resulting loss function for general transportation costs and derive gradient formulas. Physical measures are approximated as fixed points of suitable PDE-based Perron--Frobenius operators. Test cases discussed in the paper include common low-dimensional dynamical systems.
\end{abstract}



\section{Introduction}\label{sec:intro}
The problem of parameter identification in dynamical systems is common in many areas of science and engineering, such as signal processing \cite{feng2008reconstruction}, optimal control \cite{fradkov2005control, nakagawa1999chaos}, secure communications \cite{ruan2003chaotic, feng2008reconstruction}, as well as and biology \cite{rodriguez2006novel, gabor2015robust}, to mention a few. The main idea of parameter identification for a dynamical system is to identify a mathematical model of the real-world system and adapt its parameters until the simulations obtained with the mathematical model are close to experimental data. The models usually represent time-dependent processes with numerous state variables and many interactions between variables. In many applications, one can derive the form of the mathematical model from some knowledge about the process under investigation, but in general, the parameters of such a model must be inferred from empirical observations of time series data. The initial parameter values are usually based on, for instance, some preliminary knowledge of the real-world system. The type of mathematical model and the parameter identification algorithm chosen strongly influence the accuracy of the estimates. 

More formally, suppose that we have noisy observations
\[
{\bf X}^*=\left({\bf x}^*(t_0)+\eta_0,{\bf x}^*(t_1)+\eta_1,\cdots,{\bf x}^*(t_n)+\eta_n\right),    
\]
where $\{t_0,t_1,\cdots,t_n\}$ are sampling times, ${\bf x}^*$ is the solution of the autonomous dynamical system $\dot{\bf x}=v({\bf x},\theta^*)$, and $\{\eta_0,\eta_1,\cdots,\eta_n\}$ are measurement errors or uncertainties. The goal is to find $\theta^*$ from ${\bf X}^*$.

Most common parameter estimation techniques estimate $\theta$ by integrating $\dot{\bf x}=v({\bf x},\theta)$ and fitting the resulting trajectory ${\bf X}(\theta)=\left({\bf x}(t_0),{\bf x}(t_1),\cdots,{\bf x}(t_n)\right)$ to data ${\bf X}^*$ via optimization
\[
\inf\limits_{\theta \in \Theta} \|{\bf X}(\theta)-{\bf X}^*\|^2    
\]
for a suitably chosen norm $\|\cdot\|$. For a linear map $\theta \mapsto v({\bf x},\theta)$ and a quadratic norm, the problem above reduces to the least-squares problem that tends to overfit measurement errors~\cite{kostelich1992problems,jaeger1996unbiased}. For a nonlinear map $\theta \mapsto v({\bf x},\theta)$, this approach leads to a so-called single shooting method~\cite{michalik2009incremental} that uses a single initial condition to produce a trajectory. However, relying only on one trajectory may not result in meaningful approximations of the desired solution for chaotic systems due to their sensitivity to initial data. The multiple shooting algorithm deals with this issue by using multiple trajectories to estimate parameters~\cite{baake1992fitting}. For a more complete review we refer to~\cite{aguirre2009modeling} and~\cite{mcgoff2015statistical}.
Because of their universal approximation properties, neural networks and combinations of the above methods with neural networks have also been used recently for parameter identification of dynamical systems~\cite{bakker2000learning,lu2021deepxde, negrini2021system, negrini2021neural}. 

An alternative approach is to fit the time derivatives of the state. More precisely, assume that $\dot{{\bf x}}^*$ is either measured directly or estimated from ${\bf X}^*$ yielding
\[
{\bf V}^*=\left(\dot{\bf x}^*(t_0)+\xi_0,\dot{\bf x}^*(t_1)+\xi_1,\cdots,\dot{\bf x}^*(t_n)+\xi_n\right),    
\]
where $\{\xi_0,\xi_1,\cdots,\xi_n\}$ are measurement or estimation errors. The parameter estimation is then performed via an optimization problem
\[
\inf_{\theta\in \Theta} \|{\bf V}^*-v({\bf X}^*,\theta)\|^2+R(\theta)    
\]
for a suitably chosen norm $\|\cdot\|$ and a regularization $R(\theta)$, where we denote $v({\bf X}^*,\theta)=(v({\bf x}^*(t_0),\theta),v({\bf x}^*(t_1),\theta),\cdots,v({\bf x}^*(t_n),\theta))$ by slightly abusing the notation. Sparse identification of nonlinear dynamics (SINDy)~\cite{brunton2016discovering} is one such notable method, where one has a linear model $v({\bf x},\theta)=\sum_i \theta_i \psi_i({\bf x})$ with a suitably chosen dictionary of basis functions $\{\psi_i\}$ and a sparsity enforcing regularization term $R(\theta)=\|\theta\|_1$. 

We are interested in parameter estimation problems where trajectories are sensitive to initial conditions and estimation parameters. In particular, we consider the case where the time derivatives ${\bf V}^*$ cannot be estimated due to the lack of observational data, slow sampling, discontinuous or inconsistent time trajectories, and noisy measurements~\cite{bezruchko2010extracting}. The methods described above incur many challenges or are inapplicable in such settings. Hence, following~\cite{greve2019data}, we ``suppress'' the time variable and consider the state-space distribution of data
\[
    \rho^*=\frac{1}{n+1}\sum_{i=0}^n \delta_{{\bf x}^*(t_i)}.
\]

We say that a dynamical system $\dot{{\bf x}}=v({\bf x},\theta)$ admits a \textit{physical measure} $\rho(\theta)$~\cite[Definition 2.3]{young2002srb},\cite[Section 9.3]{medio2001}, if for a Lebesgue positive set of initial conditions ${\bf x}(0)=x$, one has that
\[
    \rho(\theta)=\lim_{T\to \infty} \frac{1}{T}\int_0^T \delta_{{\bf x}(t)}dt.
\]
Therefore, as an alternative, we can fit physical measures instead of trajectories for systems admitting such measures. In this work, we focus on dynamical systems with a unique physical measure.
More precisely, the parameter estimation problem reduces to the optimization problem
\begin{equation}\label{eq:param_estim_gen}
    \inf_{\theta \in \Theta} f(\theta):= d(\rho_\epsilon(\theta),\rho^*),
\end{equation}
where $\rho_\epsilon(\theta)$ is an approximation of $\rho(\theta)$ with an approximation (regularization) parameter $\epsilon>0$, and $d$ is a suitable metric in the space of probability measures.

Note that the definition of physical measures reflects their stability with respect to perturbations of initial conditions. Additionally, $\rho^*$ can provide an accurate estimate of $\rho(\theta^*)$ even if we perform slow sampling; that is, when the time derivatives ${\bf V}^*$ cannot be estimated (\Cref{sec:invariant_measure_sample}).

The difficulty and efficiency of the parameter estimation problem \eqref{eq:param_estim_gen} depend significantly on the choice of the approximation method $\rho_\epsilon$ and the metric $d$. The Wasserstein metric from optimal transportation (OT)~\cite{villani03} has recently gained popularity as a metric of choice in numerous fields such as image processing~\cite{haker2004optimal}, machine learning~\cite{arjovsky2017wasserstein}, large-scale inverse problems~\cite{engquist2020optimal}, and statistical inference~\cite{bernton2019parameter}, only to mention a few. Interested readers may further refer to~\cite{peyre2019computational}. The Wasserstein metric is beneficial for several reasons. First, it is well-defined for singular measures and, unlike the Kullback--Leibler divergence, reflects both the local intensity differences and the global geometry mismatches~\cite{engquist2020optimal}. Additionally, the $L^p$ and total variation (TV) norms lead to weakly pronounced minima with small basins of attraction when the supports are disjoint or only partially intersect.
Second, recent works in both deterministic and Bayesian inverse problems have demonstrated that the Wasserstein metric is robust to noise~\cite{engquist2020quadratic,dunlop2020stability}, a preferred property to help avoid overfitting when we have noisy time trajectories, modeling error, and numerical errors.

Our first main goal of this work is to study OT distances as the objective function for parameter identification problems in dynamical systems building on insights from~\cite{greve2019data}. An important element of the method~\eqref{eq:param_estim_gen} is the surrogate model $\rho_\epsilon(\theta)$. In~\cite{greve2019data}, the authors build a histogram from a single long-time trajectory, where $\epsilon$ is the bin width. Although effective, one drawback of this approximation method is the inability of differentiating $\rho_\epsilon(\theta)$ with respect to $\theta$. Consequently, it relies on a potentially slow derivative-free optimization method to solve~\eqref{eq:param_estim_gen}. Our second main goal is to explore an alternative scheme for the approximation $\rho_\epsilon(\theta)$ that is differentiable in $\theta$, and rigorously study the regularity of $f(\theta)$ in~\eqref{eq:param_estim_gen}. One can then devise more efficient gradient-based optimization algorithms to solve~\eqref{eq:param_estim_gen}.

In this work, we propose a partial differential equation (PDE)-based approximation method for $\rho(\theta)$. Note that $\rho(\theta)$ is a distributional solution of the stationary continuity PDE
\begin{equation}\label{eq:continuity}
    -\nabla \cdot (v({\bf x},\theta) \rho({\bf x}))=0.
\end{equation}
Hence, we consider a regularized solution $\rho_\epsilon(\theta)$ of~\eqref{eq:continuity} and turn \eqref{eq:param_estim_gen} into a PDE-constrained optimization problem. We choose the teleportation regularization from Google's PageRank algorithm~\cite{gleich2015pagerank} because of its simplicity in implementation and other favorable properties such as the uniqueness, absolute continuity, and differentiability (with respect to $\theta$) of $\rho_\epsilon(\theta)$. The numerical method for computing $\rho_\epsilon(\theta)$ is based on its representation as a fixed point of a suitable Perron--Frobenius operator.

Approximating physical measures by PDE and fixed points of Perron--Frobenius operators instead of directly simulating single long-time trajectories is not new~\cite{dellnitz2002set,allawala2016statistics}. Some of these methods come with rigorous convergence guarantees, especially for uniformly hyperbolic systems~\cite[Theorem 4.14]{dellnitz2002set}, and are more computationally efficient because of considering $\rho_\epsilon(\theta)$ that are supported on tight covers of $\operatorname{supp}(\rho)$~\cite[Section 4]{dellnitz2002set}. However, the differentiability of the resulting approximations with respect to the parameters is unclear and warrants separate careful analyses. Here, we do not analyze the convergence of $\rho_\epsilon(\theta)$ to $\rho(\theta)$, but the numerical evidence in \Cref{sec:invariant_measure_tele} and the discussion in \Cref{sec:FPE_FWD} suggest that this convergence occurs for a suitable class of dynamical systems. Instead, we focus on studying the properties of OT-based distances and the viability of the overall approach at the expense of employing a less accurate yet more straightforward approximation method for the differentiability analysis. Thus, our work serves as a foundation for possibly other OT-based techniques with different but differentiable approximation methods for the physical measures. Formally, we assume that (1) the dynamical system of interest, $\dot{x} = v(x,\theta)$ where $\theta\in \Theta$,  has one unique physical invariant measure, and (2) the distributional solution to~\eqref{eq:continuity} with the same $v(x,\theta)$ is unique and recovers the physical invariant measure to the dynamical system. 
We refer to \Cref{sec:FPE_FWD} for more details.

The discussion above leads to our next essential contribution: the regularity analysis of the optimal transport cost with respect to the inference parameter for generic cost functions; see~\Cref{sec:f_diff}. Although the gradient formula is well known in the literature, its validity analysis seems to be missing except in special cases where the optimal transport cost can be calculated explicitly \cite[Lemma 2.4]{rippl16}. In the non-parametric setting, such analysis can be found in \cite[Theorem 2.4]{sommerfeld18} for probability measures on finite spaces and~\cite[Proposition 7.17]{santambrogio15} for probability measures on $\R^d$.  For probability measures modeled by push-forward maps, see~\cite{arjovsky2017wasserstein}.

Similar to related results in the literature, we rely on Kantorovich's formulation of the OT problem and the regularity theory of optimal value functions \cite{bonnans00}. Under rather mild conditions, we prove that the transportation cost is directionally differentiable everywhere. In general, the directional derivative is nonlinear and depends on the structure of Kantorovich potentials. 
To this end, we find a sufficient condition in terms of the geometry of the optimal transport plans that guarantees the linearity of the directional derivative providing a descent direction for the optimal transport cost. To the best of our knowledge, this condition is new in the literature.




The paper is arranged as follows. In~\Cref{background}, we review challenges of the chaotic dynamics, the advantages provided by the PDE perspective~\eqref{eq:continuity}, and a short introduction to optimal transport. In~\Cref{sec:FPE_FWD}, we describe a regularized forward problem based on the PDE perspective and discuss the numerical scheme that enforces positivity and strict mass conservation. The solution to the forward problem is computed as finding the dominant eigenvector of a Markov matrix.
In~\Cref{sec:f_diff}, we present theoretical regularity analysis for evaluating gradients of optimal transport costs with respect to the model parameters. In~\Cref{gradient_descent}, we introduce two different ways to compute gradients for our PDE-constrained optimization problem using the implicit function theorem and the adjoint-state method. Numerical results for the Lorenz, R\"ossler, and Chen systems are presented in~\Cref{results}. In~\Cref{conclusion}, we summarize our results and describe several future research directions.

\section{Background}\label{background}
In this section, we present the essential background of dynamical systems and optimal transportation theory.

\subsection{Dynamical Systems}
This section reviews some basic terminologies in the field of dynamical systems that will appear throughout the paper.
\subsubsection{Chaotic Dynamical Systems} \label{SubSec:Dynamical}
A continuous-time dynamical system represents the behavior of a system in which the time-dependent flow of a point in a geometrical state space, $\textbf{x}$, is governed by a function of that state, $v(\textbf{x})$, such that
\begin{equation}\label{eq: dynamicalSystem}
\frac{d\textbf{x}}{dt}=\dot{\textbf{x}}=v(\textbf{x}).
\end{equation}

This first-order ordinary differential equation (ODE) can be viewed as the trajectory of a point in Lagrangian coordinates. While linear first-order dynamical systems, $\dot{\textbf{x}}=A\textbf{x}$, admit only stable, unstable, and periodic solutions, the more general class of nonlinear dynamical systems can exhibit a range of more complex long time behaviors due to locally bounded regions of instability.  It is this local region of instability that enables the emergence of chaotic behavior.

While a formal definition of chaos remains elusive, it is generally characterized by bifurcation and sequences of period doubling, transitivity and dense orbit, sensitive dependence to initial conditions, and expansivity; see \cite{effah2018study} for more details.  In particular, it is this sensitive dependence on initial conditions that results in the apparent randomness characteristic of chaotic systems. This randomness results from a combination of local instability causing exponential divergence of nearby trajectories and state-space mixing that occurs when this exponential divergence is re-stabilized such that a nontrivial attractor forms. This combination makes long-time predictions impossible despite the purely causal nature of the governing system.  It is also this sensitivity that makes the classical trajectory-based parameter inference problem challenging when the observed dynamics are obscured by noise, slow sampling, and other corruption, as described in~\Cref{sec:intro}.

\subsubsection{From Trajectory Samples to the Physical Measure}
We shift from the trajectory-based to distribution-based perspective to remedy the aforementioned stability and data availability issues. Mathematically, statistical properties of~\eqref{eq: dynamicalSystem} can be characterized by the occupation measure $\rho_{x,T}$ defined as
\begin{equation} \label{eq:our_data}
\rho_{x,T}(B) =  \frac{1}{T}\int_{0}^{T} \mathds{1}_B ({\bf x}(s)) ds = \frac{\int_{0}^{T} \mathds{1}_B ({\bf x}(s)) ds}{\int_{0}^{T} \mathds{1}_{\R^d} ({\bf x}(s)) ds},
\end{equation}
where $T>0$, $\mathds{1}$ is the indicator function, $B$ is any Borel measurable set, and ${\bf x}(\cdot)$ is the time-dependent trajectory starting at $x$. System \eqref{eq: dynamicalSystem} has robust statistical properties if there exists a set of positive Lebesgue measure $U$ and an invariant probability measure $\rho$ such that $\rho_{x,T}$ converges weakly to $\rho$ for all initial conditions $x\in U$. Such $\rho$ are called \textit{physical}~\cite[Definition 2.3]{young2002srb}, \cite[Section 9.3]{medio2001}. For suitable classes of dynamical systems, such as Axiom A, physical measures reduce to so called \textit{Sinai-Ruelle-Bowen (SRB) measures}~\cite{dellnitz2002set,young2002srb,medio2001}.

In general, the existence and properties of such measures are rather intricate and require careful analysis. For a more detailed account on these topics, we refer to \cite{young2002srb} for general systems, and \cite{tucker02,tucker99} for the Lorenz system. Furthermore, in some cases, one can recover $\rho$ as the zero-noise limit of stationary measures of the corresponding stochastic dynamical systems \cite{cowieson_young_2005,huang18,kifer1986small,dellnitz2002set}.

As we will show in \Cref{sec:FPE_FWD}, direct simulation of $\rho$ for parameter identification faces the difficulty of not having access to the gradients of the loss function. Consequently, one has to rely on gradient-free space-search methods. Motivated by these challenges, we take a PDE perspective on $\rho$ and formulate the parameter inference problem as a PDE-constrained optimization.

\subsection{Optimal Transportation}
In this subsection, we give a brief overview of the topic of optimal transportation (OT), first brought up by Monge in 1781. 

We first introduce the original Monge's problem. Let $\Omega\subset \mathbb{R}^d$ be an arbitrary domain, and $\mu,\nu \in \Pp(\Omega)$ arbitrary probability measures supported in $\Omega$. 
A transport map $T: \Omega \rightarrow \Omega$ is mass-preserving if for any measurable set $B \subseteq  \Omega$
\[
\mu (T^{-1}(B)) = \nu(B).
\]
If this condition is satisfied, $\nu$ is said to be the push-forward of $\mu$ by $T$, and we write $\nu = T_\sharp \mu$. In case $\mu,\nu$ are absolutely continuous; that is, $d\mu(x)=f(x)dx $ and $d\nu(y)=g(y) dy $, we have that $T$ is a mass-preserving map if
\[
    f(x)=g(T(x)) \cdot |\det \left(\nabla T(x)\right)|,\quad x\in \Omega.
\]
The transport cost function $c(x,y)$ maps pairs $(x,y) \in \Omega\times \Omega$ to $\mathbb{R}\cup \{+\infty\}$, which denotes the cost of transporting one unit mass from location $x$ to $y$. The most common choice of $c(x,y)$ is $|x-y|^p$, $p\in \mathbb{N}$, where $|x-y|$ denotes the Euclidean distance between vectors $x$ and $y$. Given a mass-preserving map $T$, the total transport cost is
\[
    \int_\Omega c(x,T(x))f(x)\,dx.
\]
While there are many maps $T$ that can perform the relocation, we are interested in finding the optimal map that minimizes the total cost. So far, we have informally defined the optimal transport problem, which induces the so-called Wasserstein distance defined below, associated to cost function $c(x,y) = |x-y|^p$.
\begin{definition}[The Wasserstein distance]
We denote by $\mathscr{P}_p(\Omega)$ the set of probability measures with finite moments of order $p$. For all $p \in [1, \infty)$,   
\begin{equation}\label{eq:static}
W_p(\mu,\nu)=\left( \inf _{T_{\mu,\nu}\in \mathcal{M}}\int_{\Omega}\left|x-T_{\mu,\nu}(x)\right|^p d\mu(x)\right) ^{\frac{1}{p}},\quad \mu, \nu \in \mathscr{P}_p(\Omega),
\end{equation}
where $\mathcal{M}$ is the set of all maps that push forward $\mu$ into $\nu$.
\end{definition}

The definition~\eqref{eq:static} is the original static formulation of the optimal transport problem with a specific cost function. In mid 20th century, Kantorovich relaxed the constraints, turning it into a linear programming problem, and also formulated the dual problem~\cite{santambrogio15}. Instead of searching for a map $T$, a transport plan $\pi$ is considered, which is a measure supported in the product space $\Omega\times \Omega$. The Kantorovich problem is to find an optimal transport plan as follows:
\begin{equation}\label{eq:T_c}
\T_c(\mu,\nu)=\inf_{\pi} \bigg\{ \int_{\Omega \times \Omega} c(x,y) d\pi\ |\ \pi \geq 0\ \text{and}\ \pi \in \Pi(\mu, \nu) \bigg\} ,
\end{equation}
where $\Pi (\mu, \nu) =\{ \pi \in \mathscr{P}( \Omega \times \Omega)\ |\ (P_1)_\sharp \pi = \mu, (P_2)_\sharp \pi = \nu  \}$. Here, $\mathscr{P}(\Omega \times \Omega)$ stands for the set of all the probability measures on $\Omega \times \Omega$, functions $P_1(x,y) = x$ and $P_2(x,y) = y$ denote projections over the two coordinates, and $(P_1)_\sharp \pi$ and $(P_2)_\sharp \pi$ are two measures obtained by pushing forward $\pi$ with these two projections.

Since every transport map determines a transport plan of the same cost, Kantorovich's problem is weaker than the original Monge's problem. If the cost function $c(x,y)$ is of the form $|x-y|^p$ and $\mu$ and $\nu$ are absolutely continuous with respect to the Lebesgue measure, solutions to the Kantorovich and Monge problems coincide under certain conditions. When $p>1$, the strict convexity of $|x-y|^p$ guarantees that there is a unique solution to Kantorovich's problem~\eqref{eq:T_c} which is also the unique solution to Monge's problem~\eqref{eq:static}.

\section{The Forward Model}\label{sec:FPE_FWD}

While matching shadow state-space density in~\cite{greve2019data} provided a potential route to resolve issues related to the chaotic divergence of state-space trajectories and data availability, the direct estimation of state-space density from trajectory data still retained two major challenges. 
One significant issue was the inability to efficiently calculate a gradient of the Wasserstein metric with respect to the parameters, forcing the reliance on evolutionary or other gradient-free optimization methods.
Another major issue was related to the time required to converge to the density estimate asymptotically, as particularly highlighted in~\cite[Fig.\ 7]{greve2019data}, where the self-self Wasserstein metric 
is observed to oscillate as it converges with more ODE time steps. This slow convergence is related to the long and intermittent switching times between lobes of the butterfly attractor. While the invariant measure of the Lorenz system is known to exist~\cite{tucker99}, the long measurement times with respect to the switching times complicate the parameter inference problem. The problem is exacerbated in more expensive and complicated dynamics such as the thruster model~\cite{greve2019data}.

To address these challenges, we instead directly solve for the solution of the stationary continuity equation~\eqref{eq:continuity}.  This choice not only removes the issue of slow convergence with respect to the slowest system processes but also provides a forward model that can be differentiated for building the required gradients needed to tackle the parameter inference problem directly.  This alternative forward model follows the approach described in~\cite{bewley2012efficient} in converting from the trajectory samples to the probability measure for the Bayesian estimation problem, as detailed in~\Cref{subsec:LinearAdvect}, but then recasts this forward Perron--Frobenius operator as a Markov process for determining the steady-state solution as described in~\Cref{sec:finding_stationary}.

Our approach is close in spirit to other cell-based or grid-based frameworks that introduce a suitable Perron--Frobenius operator and compute its fixed points~\cite[Section 4]{dellnitz2002set}. Some of these methods, such as the software package GAIO proposed in \cite{dellnitz1998adaptive,dellnitz1999approximation} by Dellnitz and Junge, represent the attractors via a hierarchy of covers by cells: cells that do not intersect the support of the invariant measure are ignored so that the data structures and computational requirements for this method are smaller than the ones required for our grid-based approach. In some cases, such as uniformly hyperbolic systems, these methods come with convergence guarantees~\cite[Theorem 4.14]{dellnitz2002set}. Many other subdivision methods have been successfully applied to the numerical analysis of complex dynamical behavior, see for instance  \cite{dellnitz1997almost,eidenschink1997exploring,siegmund2006approximation}. A more comprehensive list of examples can be found in \cite{dellnitz2001algorithms,givon2004extracting}.

We regularize our Perron--Frobenius operator via teleportation regularization from Google's PageRank method~\cite{gleich2015pagerank}, which ensures the uniqueness and regularity of the fixed point. This step is similar to stochastic perturbation techniques for approximating physical measures~\cite{cowieson_young_2005,huang18,kifer1986small,dellnitz2002set}. Intuitively, teleportation amounts to stopping the dynamics at a random time and restarting it from a randomly chosen initial point. The regularization parameter $\epsilon$ controls the restarting frequency: the smaller $\epsilon$, the rarer we restart. This regularization is somewhat similar to ``snapshot attractors'' described in~\cite{robin17} where attractors are estimated by following the dynamics from randomly chosen initial conditions for a fixed time. Here, we do not analyze the convergence of $\rho_\epsilon(\theta)$ to the physical invariant measure, but the numerical evidence in Section \ref{sec:invariant_measure_tele} suggest that this convergence does take place for the tested examples. Intuitively, if we restart the dynamics from the basin of attraction and do so very rarely, we should approximate the physical measure. Additionally, general results in~\cite{kifer1986small} hint at a convergence result similar to~\cite[Theorem 4.14]{dellnitz2002set} for uniformly hyperbolic attractors. Analyzing the convergence of our model and the differentiability of other forward models described here is an exciting future research direction that we plan to pursue. For additional methods based on Markov partitions and chains we refer to~\cite{bowen1975equilibrium,froyland2001extracting,fiedler2002handbook}. Formally, we assume that (1) the dynamical system of interest, $\dot{x} = v(x,\theta)$ where $\theta\in \Theta$,  has one unique physical invariant measure, and (2) the distributional solution to~\eqref{eq:continuity} with the same $v(x,\theta)$ is unique and recovers the physical invariant measure to the dynamical system.



\subsection{From Linear Advection to Stationary Eigenvectors} \label{subsec:LinearAdvect}

In converting the dynamical system from the trajectory samples to the probability measure, the governing equation is converted from a nonlinear ODE for the system state ``point'', ${\bf x}$, to a linear PDE~\eqref{eq:continuity} for the state space density $\rho({\bf x})$. 


Note that a causal dynamical system includes no diffusion. It then corresponds to~\eqref{eq:continuity}, a linear advection of probability density in state space. \Cref{subsec:finite_volume} describes a particular simple low-order discretization of this linear advection problem.  While adding physical diffusion is a relatively simple modification of the numerical method, the more significant issue with this approach relates to excess diffusion. Although the zero diffusion case can be relaxed for stochastic dynamical systems where $D_{ij}\neq0$, the upwinding scheme required to stabilize the advection introduces an artificial diffusion, which is the predominant numerical error as described in~\cite{bewley2012efficient}.  This numerical diffusion is expected to dominate physical diffusion for the moderate spatial resolution that is tractable for the forward model unless the dynamics of the system are highly stochastic.  As this numerical diffusion is irreducible at finite computational cost, the addition of finite diffusion to the ODE model is explored in~\Cref{subsec:Invariant} when attempting to understand the class of problems for which inference with respect to the binned direct ODE solution is viable.

\subsection{Finite Volume Discretization}\label{subsec:finite_volume}
A finite volume discretization of the resulting continuity equation defined on the domain $\Omega$, as described in \cite{bewley2012efficient}, is then obtained. The finite volume discretization combined with a zero-flux boundary condition, $v=0$ on the boundaries $\partial \Omega$, enforces strict mass conservation whenever the discrete integration by parts formulation is used~\cite{fernandez2014generalized}. Only the first-order operator split upwind discretization is used in this work to enforce positivity of the probability density, as will be shown to be a consequence of the form of the discrete operator. 


We first discretize~\eqref{eq:continuity} on a $d-$dimensional uniform mesh in space and time with no added diffusion, which gives us the following equation for the explicit time evolution of the probability density,
\[
    \frac{\rho^{(l+1)}(x_i)-\rho^{(l)}(x_i)}{\Delta t} = -\sum_{i_d=1}^{d} \frac{F_{(i_d)}^{(l)}(x_i+\Delta x_{(i_d)}/2) - F_{(i_d)}^{(l)}(x_i-\Delta x_{(i_d)}/2)} {\Delta x_{(i_d)}}.
\]
Here, the point $x_i$ refers to the $i^{th}$ cell center vector and $\Delta x_{(i_d)}$ refers to the mesh spacing in the $i_d$-th direction, $i_d = 1,\ldots,d$. The upwind $i_d-$direction flux at the $l$-th time step, $F_{(i_d)}^{(l)}$, is then approximated using face center velocity assuming uniform density within the cell centered at $x_i$ as follows:
\[
    F_{(i_d)}^{(l)}\left(x_i-\frac{\Delta x_{(i_d)}}{2}\right) = v_{(i_d-\frac{1}{2})}^+ \rho^{(l)}( x_i-\Delta x_{(i_d)}) + v_{(i_d-\frac{1}{2})}^-  \rho^{(l)}( x_i ),
\]
where the upwind velocities $v_{(i_d)}^+=\max(v_{(i_d)},0)$ and $v_{(i_d)}^-=\min(v_{(i_d)},0)$ refer to the $i_d$-th component of the velocity vector split between positive and negative values, and 
\[
v_{(i_d-\frac{1}{2})}^+ := v_{(i_d)}^+\left( x_i-\frac{\Delta x_{(i_d)}}{2} \right),\quad v_{(i_d-\frac{1}{2})}^- := v_{(i_d)}^-\left( x_i-\frac{\Delta x_{(i_d)}}{2} \right).
\]
Inserting these fluxes into the discrete equation yields the following expression for the future time density, $\rho^{(l+1)}$.
\[
    \rho^{(l+1)}_0= \rho^{(l)}_0 + \Delta t \sum_{i_d=1}^{d} \frac{ \left( v_{(i_d-\frac{1}{2})}^+ \rho^{(l)}_- +v_{(i_d-\frac{1}{2})}^- \rho^{(l)}_0 \right) - \left( v_{(i_d+\frac{1}{2})}^+ \rho^{(l)}_0 +v_{(i_d+\frac{1}{2})}^- \rho^{(l)}_+ \right)  }{\Delta x_{(i_d)}},
\]
where $\rho^{(l)}_0  = \rho^{(l)}(x_i) $, $\rho^{(l)}_- = \rho^{(l)} (x_i - \Delta x_{(i_d)})$ and $\rho^{(l)}_+ = \rho^{(l)} (x_i + \Delta x_{(i_d)})$. The equation above can be rewritten as a matrix-vector format:
\[
    \rho^{(l+1)} = \rho^{(l)} + K_{mat} \rho^{(l)} = (I+K_{mat})  \rho^{(l)}.
\]
For steady state distributions, $\rho^{(l+1)}=\rho^{(l)}=\rho^{eq}$. This corresponds to finding a nonzero solution  $\rho^{(eq)}$ to the following linear system 
\[
    K_{mat} \rho^{(eq)} =   \left[ \sum_{i_d=1}^{d} \frac{\Delta t}{\Delta x_{(i_d)}}K_{(i_d)} \right] \rho^{(eq)} = 0
\]
where for $i_d= 1,\ldots,d$ we have
\begin{equation}\label{eq:Kd}
K_{(i_d)} = 
\begin{bmatrix}
  \ddots & & &  & \\
  & -v_{(i_d-\frac{3}{2})}^-  & &  & \\
  \ddots &\vdots& & &\\
 &  v_{(i_d-\frac{3}{2})}^- - v_{(i_d-\frac{1}{2})}^+ &  -v_{(i_d-\frac{1}{2})}^- & &  \\
   \ddots & \vdots &\vdots  &  &\\
  & +v_{(i_d-\frac{1}{2})}^+ &   v_{(i_d-\frac{1}{2})}^- - v_{(i_d+\frac{1}{2})}^+& -v_{(i_d+\frac{1}{2})}^-  &\\
 &    & \vdots & \vdots  & \ddots\\
 & &  +v_{(i_d+\frac{1}{2})}^+ &  v_{(i_d+\frac{1}{2})}^- - v_{(i_d+\frac{3}{2})}^+&\\
 &  & & \vdots &\ddots   \\
 &  & & +v_{(i_d+\frac{3}{2})}^+ & \\
 &  & &  & \ddots  
\end{bmatrix}.
\end{equation}

We remark that each $K_{(i_d)}$, $i_d = 1,\dots, d$, is a tridiagonal matrix, while the offsets for the three diagonals vary for different $i_d$. For example, consider the case that $\Omega \subseteq \mathbb{R}^3$ is a cuboid, discretized with grid size $n_x,n_y,n_z$ in the $x,y,z$ dimension, respectively. Then, $K_{(1)}$ is nonzero at the first lower diagonal, the main diagonal, and the first upper diagonal; $K_{(2)}$ is nonzero at the $n_x$-th lower diagonal, the main diagonal, and the $n_x$-th upper diagonal; $K_{(3)}$ is nonzero at the $(n_x\times n_y)$-th lower diagonal, the main diagonal, and the $(n_x\times n_y)$-th upper diagonal.

We highlight that the solution $\rho^{(l)}$ at any $l$-th time step satisfies the mass conservation property. That is,
\[
    \rho^{(l)} \cdot \mathbf{1}= \rho^{(l+1)} \cdot \mathbf{1} = \rho^{(eq)} \cdot \mathbf{1},\quad \text{where } \mathbf{1} = [1,1,\ldots,1]^\top.
\]
It is a direct consequence of the fact that columns of $K_{mat}$ sum to zero. Note also that the off-diagonal terms are all positive or zero while the diagonal terms are all negative or zero by construction. One can construct a \textit{column-stochastic matrix} $M$
\[    
M = I+cK_{mat}.
\]
$M$ can be positive definite if we ensure that $c$ is small enough.

Since the main focus of this paper is parameter identification, the velocity field $v$ is parameter-dependent. Thus, we will highlight the dependency on the parameter $\theta$ by using notation $v(\theta)$, $K_{mat}(\theta)$, $K_{(i_d)}(\theta)$, and $\rho^{(eq)}(\theta)$ hereafter.

The upper bound on $c$ unsurprisingly also depends on $\theta$. Nevertheless, if we assume that $v$ depends continuously on $\theta$ and we operate in a bounded domain $\Omega$, we can choose $c$ small enough to serve all $\theta$-s of interest. For instance, we can choose
\begin{equation}\label{eq:CFL}
   0< c<\min_{i_d}\frac{\Delta x_{(i_d)}}{2 \Delta t \max\limits_{x\in \Omega, \theta \in \Theta} |v_{(i_d)}(x,\theta)|}.
\end{equation}

\subsection{Finding the Stationary Distribution of a Markov Chain} \label{sec:finding_stationary}
From the previous section, we learned that $\rho(\theta)$ is the solution of
\begin{equation}\label{eq:main_direct-1}
    M(\theta) \rho = \rho,\quad \rho \cdot \mathbf{1}=1,
\end{equation}
where 
\[
   M(\theta) = I+cK_{mat}(\theta),\quad K_{mat}(\theta)=\sum_{i_d=1}^d \frac{\Delta t}{\Delta x_{(i_d)}} K_{(i_d)}(\theta),
\]
with $K_{(i_d)}(\theta)$ given in \eqref{eq:Kd}, and $c$ is chosen to satisfy~\eqref{eq:CFL}. While the matrix, $M$, was built from a finite volume causal flow model, it was noted that this flux also approximates a discrete cell-to-cell transition probability for a point randomly sampled from the volume of one cell to its neighbor cells, which mirrors the propagator of a Markov chain as described in~\cite{KaiserJFM2014}. 

A priori we have that the off-diagonal entries of $M(\theta)=I+cK_{mat}(\theta)$ are non-negative. Additionally, we know that $M(\theta)$ is column stochastic. Thus, by Gershgorin's theorem~\cite{GVL12} we have that the spectral radius of $M$ is not greater than one. On the other hand, we know $\mathbf{1}=[1,1,\cdots,1]^\top$ is an eigenvector for $M^\top$ which is a row-stochastic matrix, and so $\lambda=1$ is an eigenvalue for both $M$ and $M^\top$. The spectral radius of $M$ has to be equal to 1. Furthermore, by a limiting argument, we can show that the eigenspace of $M$ corresponding to the eigenvalue $\lambda=1$ contains vectors with non-negative entries.

However, the dimension of this eigenspace may be bigger than one, which complicates our analysis. Thus, we regularize $M$ via the so-called teleportation trick, which is well-known from Google's PageRank method~\cite{gleich2015pagerank}. That is, given a small positive constant $\epsilon$, we consider  
\begin{equation}\label{eq:M_eps}
    M_{\epsilon}(\theta)=(1-\epsilon) M + \epsilon n^{-1} \mathbf{1}~\mathbf{1}^\top = (1-\epsilon)(I+c K_{mat}(\theta))+\frac{\epsilon}{n} \mathbf{1}~\mathbf{1}^\top.
\end{equation}
Note that the off-diagonal entries of $M_{\epsilon}$ are at least $\frac{\epsilon}{n}>0$. The regularization also connects all cells, achieving similar regularizing effects by having a diffusion term. Moreover, $M_{\epsilon}$ is still column-stochastic. Based on the following Perron--Frobenius Theorem, the spectral radius of $M_\epsilon$ must be $1$. 
\begin{theorem}[Perron--Frobenius Theorem~\cite{meyer2000matrix}] \label{thm:PF}
If all entries of a Markov matrix $A$ are positive, then $A$ has a unique equilibrium:
there is only one eigenvalue equal to $1$. All other eigenvalues are strictly smaller than $1$.
\end{theorem}

Consequently, the eigenspace $\{\rho: M_{\epsilon}(\theta) \rho = \rho \}$ is one-dimensional and has a generator with all positive entries. Hence, the equation
\begin{equation}\label{eq:main_direct}
    M_{\epsilon}(\theta) \rho = \rho,\quad \rho \cdot \mathbf{1}=1,\quad \rho>0,
\end{equation}
has a unique solution that converges to a solution of \eqref{eq:main_direct-1} as $\epsilon \to 0$.
We can analyze the error between $\rho_0$ and $\rho_\epsilon$ where 
\[
 M\rho_0 = \rho_0,\quad M_\epsilon \rho_\epsilon =\rho_\epsilon, \quad \rho_0 \cdot \mathbf{1}=\rho_\epsilon \cdot \mathbf{1}=1.
\]
The error analysis traces back to the classical root-finding problem. We define $\Delta \rho_\epsilon = \rho_\epsilon - \rho_0$. Using the forward error analysis, we obtain that
\[
(M-I)\Delta \rho_\epsilon =  (M-I) \rho_\epsilon = \epsilon \left(M-n^{-1} \mathbf{1}~\mathbf{1}^\top \right) \rho_\epsilon, \quad  \Delta \rho_\epsilon \cdot \mathbf{1} = 0.
\]
Solving for $\Delta \rho_\epsilon$ from the linear system above can improve the current ``root'' $\rho_\epsilon$, which is precisely the principle behind Newton's method. Using backward error analysis, starting from $M_\epsilon \rho_\epsilon =\rho_\epsilon$, we obtain that
\[
\left( (1-\epsilon) M + \epsilon n^{-1} \mathbf{1}~\mathbf{1}^\top - I \right) \left(\rho_0 + \Delta \rho_\epsilon \right) = 0.
\]
Up to the first-order terms, we have
\[
(M-I) \Delta \rho_\epsilon = \epsilon \left( M  - n^{-1} \mathbf{1}~\mathbf{1}^\top \right) \rho_0, \quad  \Delta \rho_\epsilon \cdot \mathbf{1} = 0.
\]
The above equation implies that $\| \Delta \rho_\epsilon \|$ is $\mathcal{O}(\epsilon)$, showing the convergence $\rho_\epsilon \rightarrow \rho_0$ as we decrease $\epsilon$. This is further verified by our numerical examples in~\Cref{sec:invariant_measure_tele}.

Numerically, the problem~\eqref{eq:main_direct} can be solved by mature tools from numerical linear algebra such as the power method and the Richardson iteration~\cite{gleich2015pagerank}. We present one direct solve method in~\Cref{sec:numerical_scheme_1} using the sparsity of $K_{mat}$.

\section{Optimal Transport for Parameter Inference}\label{sec:f_diff}
Here, we discuss gradient evaluation of optimal transport-based costs with respect to the inference parameters. 
Assume that $\Omega \subset \mathbb{R}^d$ is a compact set, and $c:\Omega^2 \to \mathbb{R}$ is a continuous cost function. 
The main goal of this section is to discuss the differentiability of the objective function
\[
    f(\theta)=\T_c \big(\rho(\cdot,\theta),\rho^*\big),\quad \theta \in \Theta,
\]
where  $\{\rho(\cdot,\theta)\}_{\theta \in \Theta}$ is a family of parameter-dependent probability measures on $\Omega$, and $\T_c$ is the optimal transport cost defined in \eqref{eq:T_c}. 
Throughout the paper, we assume that $\Omega \subset \R^d$ is compact, $\rho^* \in \Pp(\Omega)$ is an arbitrary probability measure, and
\begin{itemize}
    \item[\textbf{A1}.] $\Theta \subset \mathbb{R}^m$ is an open set, and $\{\rho(\cdot,\theta)\}_{\theta \in \Theta} \subset \Pp(\Omega)$ is a family of absolutely continuous probability measures.
    \item[\textbf{A2}.] For a.e.\ $x\in \Omega$ the mapping $\theta \mapsto \rho(x,\theta)$ is differentiable, and $|\nabla_\theta \rho(x,\theta)| \leq \eta (x),~\theta \in \Theta$, for some $\eta \in L^1(\Omega)$. Note that by slightly abusing the notation, we use the same notation for probability measures and their densities.
    \item[\textbf{A3}.] $c:\Omega^2 \to \mathbb{R}$ is continuous and nonnegative.
\end{itemize}
Occasionally, we need the following hypothesis.
\begin{itemize}
    \item[\textbf{A4}.] For a.e.\ $x\in \Omega$ the mapping $\theta \mapsto \rho(x,\theta)$ is locally semiconvex, and $\nabla^2_\theta \rho(x,\theta) \geq - h(x),~\theta \in \Theta$, for some $h \in L^1(\Omega)$.
\end{itemize}
Proofs for results of this section can be found in~\Cref{sec:proofs}.

\subsection{Preliminaries}
First, we recall preliminary results from the optimal transportation (OT) theory that can be found in \cite{villani03,ags08,santambrogio15}. A key tool in OT is the Kantorovich duality \cite[Theorem 1.3]{villani03} that states 
\begin{equation}\label{eq:K-duality}
   \T_c(\mu,\nu)=\sup_{(\phi,\psi) \in \Phi_c(\mu,\nu)} \int_\Omega \phi(x) d\mu(x)+\int_\Omega \psi(y) d\nu(y),\quad \mu,\nu \in \Pp(\Omega),
\end{equation}
where $\Phi_c(\mu,\nu) \subset C(\Omega)\times C(\Omega)$ is the set of pairs $(\phi,\psi)$ such that $\phi(x)+\psi(y) \leq c(x,y)$ for all $(x,y)\in \Omega^2$. The maximizing pairs $(\phi,\psi)$ in \eqref{eq:K-duality} are called Kantorovich potentials. The $c$-transform of a function $x \mapsto \phi(x)$ is defined as
\[
    \phi^c(y)=\inf_{x\in \Omega} \left\{ c(x,y)-\phi(x) \right\}.
\]
Similarly, the $c$-transform of a function $y\mapsto \psi(y)$ is defined as
\[
    \psi^c(x)=\inf_{y\in \Omega} \left\{ c(x,y)-\psi(y) \right\}.
\]
A function $x\mapsto \phi(x)$ (resp. $y\mapsto \psi(y)$) is called $c$-concave if there exists a function $\psi$ (resp. $\phi$) such that $\phi=\psi^c$ (resp. $\psi=\phi^c$). 

Since $\Omega$ is compact and $c$ is continuous, we obtain that $c$ is bounded. Thus, the set $\Phi_c(\mu,\nu)$ in \eqref{eq:K-duality} can be further restricted to uniformly bounded pairs of conjugate $c$-concave functions; that is, pairs of $(\phi,\phi^c)\in \Phi_c(\mu,\nu)$, where $\phi=\phi^{cc}$, and $0\leq \phi \leq \|c\|_{\infty}$, $-\|c\|_{\infty}\leq \phi^c \leq 0$ \cite[Remarks 1.12-13]{villani03}. We denote this set by $K_c$.

Since the modulus of continuity of $y \mapsto c(x,y)-\phi(x)$ (resp. $x \mapsto c(x,y)-\phi^c(y)$) is bounded by that of $c$ for all $x$ (resp. $y$), $K_c$ is uniformly equicontinuous, uniformly bounded, and, consequently, precompact in $C(\Omega)\times C(\Omega)$ by the Arzel\`{a}--Ascoli theorem~\cite[Section 1.2]{santambrogio15}. Additionally, since $c$-transform is continuous under the uniform convergence, $K_c$ is compact in $C(\Omega)\times C(\Omega)$, and the existence of Kantorovich potentials in $K_c$ is guaranteed \cite[Proposition 1.11]{santambrogio15}.

\subsection{The Differentibility of the Transport Cost in the Parameter Space}

Here, we heavily rely on the Kantorovich duality \eqref{eq:K-duality} and the regularity theory of optimal value functions \cite[Chapter 4]{bonnans00}. Recall that $f$ is directionally differentiable at $\theta_0 \in \Theta$ if
\[
    \lim \limits_{t \to 0+} \frac{f(\theta_0+t \Delta \theta)-f(\theta_0)}{t}=f'(\theta_0,\Delta \theta)
\]
for all $\Delta \theta \in \R^m$~\cite[Section 2.2]{bonnans00}. Furthermore, if $\Delta \theta \mapsto f'(\theta_0,\Delta \theta)$ is linear, we say that $f$ is G\^{a}teaux differentiable at $\theta_0$ and denote by $\nabla f(\theta_0)$ the generator of this linear map.

Next, denote by $\mathcal{S}(\theta) \subset K_c$ the set of Kantorovoch potentials for the optimal transportation from $\rho(\cdot,\theta)$ to $\rho^*$.

\begin{proposition}\label{prp:dir_der_f}
Assume that \textbf{A1-A3} hold. 
\begin{itemize}
    \item[(i)] $f$ is everywhere directionally differentiable, and
\begin{equation}\label{eq:dir_der_f}
    f'(\theta_0,\Delta \theta)=\sup_{(\phi,\phi^c)\in \mathcal{S}(\theta_0)} \int_\Omega \phi(x) \nabla_\theta \rho(x,\theta_0)dx \cdot \Delta \theta 
\end{equation}
for all $\theta_0 \in \Theta$, and $\Delta \theta \in \R^m$.
    \item[(ii)] $f$ is G\^{a}teaux differentiable at $\theta_0 \in \Theta$ if and only if
    \begin{equation}\label{eq:invariance}
    \int_\Omega \phi_1(x) \nabla_{\theta} \rho(x,\theta_0) dx=\int_\Omega \phi_2(x) \nabla_{\theta} \rho(x,\theta_0) dx
\end{equation}
for all $(\phi_1,\phi_1^c),(\phi_2,\phi_2^c) \in \mathcal{S}(\theta_0)$. In this case, we have that
\begin{equation}\label{eq:grad_f}
    \nabla f(\theta_0) = \int_\Omega \phi(x) \nabla_\theta \rho(x,\theta_0) dx 
\end{equation}
for an arbitrary pair of Kantorovich potentials $(\phi,\psi) \in \Phi_c(\rho(\cdot,\theta_0),\rho^*)$.
\end{itemize}
\end{proposition}
The proof is in~\Cref{sec:dir_der_f}.

\Cref{prp:dir_der_f} asserts that $f$ is directionally differentiable at all points and that its directional derivative is a one-homogeneous closed convex function. Since we are interested in descent directions of $f$, we focus on cases when the directional derivative is a linear function and thus provides a descent direction in the form of the negative gradient. In what follows, we prove that $f$ is generically differentiable even without \eqref{eq:invariance}. Furthermore, we find sufficient structural conditions on the optimal transport plans between $\rho(\cdot,\theta_0)$ and $\rho^*$ to guarantee \eqref{eq:invariance}.

\begin{theorem}\label{thm:generic_diff_f}
Assume that \textbf{A1-A3} hold. Then $f$ is locally Lipschitz continuous, and \eqref{eq:grad_f} holds a.e.. Additionally, if \textbf{A4} holds, then $f$ is locally semiconvex, and \eqref{eq:grad_f} holds up to a set of Hausdorff dimension $d-1$.
\end{theorem}
The proof can be found in~\Cref{sec:generic_diff_f}.

There is a natural degree of freedom for Kantorovich potentials given by the addition of constants; that is, $(\phi,\phi^c)$ is a pair of Kantorovich potentials if and only if $(\phi+\lambda,\phi^c-\lambda)$ is such for an arbitrary constant $\lambda$. As a corollary of \Cref{prp:dir_der_f} we obtain that the G\^{a}teaux differentiability of $f$ is guaranteed if the addition of constants is the only degree of freedom for Kantorovich potentials.
\begin{corollary}\label{cor:diff_f}
Assume that \textbf{A1-A3} hold, and $\theta_0 \in \Theta$ is such that $\phi_2-\phi_1$ is constant $\rho(\cdot,\theta_0)$ a.e.\ for all pairs of Kantorovich potentials $(\phi_1,\psi_1)$, $(\phi_2,\psi_2)$. Then $f$ is G\^{a}teaux differentiable at $\theta_0$, and \eqref{eq:grad_f} holds.
\end{corollary}

In general, Kantorovich potentials are not unique up to constants. In what follows, we provide a sufficient condition for such uniqueness. Essentially, the optimal transportation should not amount to transportation between disjoint parts of $\operatorname{supp}(\rho(\cdot,\theta_0))$ and $\operatorname{supp}(\rho^*)$. 

More formally, assume that $\rho,\rho^* \in \Pp(\Omega)$ are such that 
$\operatorname{int}(\operatorname{supp}(\rho))\neq \emptyset$. Furthermore, denote by $\Gamma_0(\rho,\rho^*)$ the set of optimal transport plans; that is, minimizers in \eqref{eq:T_c}. We have that
\begin{equation}\label{eq:O_k}
    \operatorname{int}(\operatorname{supp}(\rho))=\cup_k O_k,
\end{equation}
where $O_k$ are disjoint open and connected sets. Next, denote by
\begin{equation}\label{eq:E_k}
\begin{split}
E_k
=&\operatorname{cl}\left(\{y:(x,y) \in \operatorname{supp}(\pi)~\mbox{for some}~x\in \operatorname{cl}(O_k),~\pi \in \Gamma_0(\rho,\rho^*)\}\right).
\end{split}
\end{equation}
In other words, $E_k$ is the set where the mass from $\operatorname{cl}(O_k)$ is transported to. 

\begin{definition}\label{def:linked}
We say that $\operatorname{cl}(O_k)$ and $\operatorname{cl}(O_l)$ are linked in the optimal transportation from $\rho$ to $\rho^*$ with a transport cost $c$, if there exist $\{i_j\}_{j=1}^m$ such that $k=i_1, l=i_m$, and $E_{i_j} \cap E_{i_{j+1}} \neq \emptyset,~1\leq j \leq m$.
\end{definition}

\begin{theorem}\label{thm:K_pot_uniq}
Assume that $c \in C^1(\Omega^2)$, $\rho,\rho^* \in \Pp(\Omega)$, and
\begin{equation}\label{eq:cl_int}
\operatorname{supp}(\rho)=\operatorname{cl}(\operatorname{int}(\operatorname{supp}(\rho))).
\end{equation}
Furthermore, suppose that $\{O_k\}$ and $\{E_k\}$ are defined as in \eqref{eq:O_k} and \eqref{eq:E_k}, respectively. Assume that all $\{\operatorname{cl}(O_k)\}$ are mutually linked. Then $\phi_2-\phi_1$ is constant $\rho$-a.e.\ for all pairs of Kantorovich potentials $(\phi_1,\psi_1)$, $(\phi_2,\psi_2)$.
\end{theorem}
The proof is presented in~\Cref{sec:K_pot_uniq}. \Cref{thm:K_pot_uniq} and \Cref{cor:diff_f} yield the following corollary.
\begin{corollary}\label{cor:f_Gat_suff}
Assume that \textbf{A1-A3} hold, and $\rho=\rho(\cdot,\theta_0)$ satisfies the hypotheses in \Cref{thm:K_pot_uniq}. Then $f$ is G\^{a}teaux differentiable at $\theta_0$.

In particular, if $\rho(\cdot,\theta_0)$ is supported on a closure of an open connected set, then $f$ is G\^{a}teaux differentiable at $\theta_0$.
\end{corollary}

The following proposition illustrates the sharpness of \Cref{cor:f_Gat_suff}. Incidentally, the same example illustrates that a smooth dependence on $\theta$ with respect to the flat $L^2$ metric does not guarantee smooth dependence on $\theta$ with respect to the Wasserstein metric.

\begin{proposition}\label{prp:non_dif_non_ac}
Assume that $\Omega=[0,4]$ and $c(x,y)=|x-y|^p$ for some $p>1$ (so that $\T_c=W_p^p$). Consider
\begin{equation*}
\begin{split}
\rho(x,\theta)=&\left(0.5+\theta\right) \chi_{[0,1]}(x)+\left(0.5-\theta\right) \chi_{[2,3]}(x),~ |\theta|<0.5,\\
\rho^*(y)=&0.5 \chi_{[1,2]}(y)+0.5 \chi_{[3,4]}(y),
\end{split}
\end{equation*}
where $\chi_A$ is the characteristic function of set $A \subset\mathbb{R}$. Then we have that
\begin{enumerate}
    \item $\{\rho(\cdot,\theta)\}$ satisfies \textbf{A1-A3}.
    \item $\{\rho(\cdot,\theta)\}$ is not absolutely continuous in $\Pp_p(\Omega)$.
    \item $\rho \mapsto W_p^p(\rho,\rho^*)$ is not G\^{a}teaux differentiable at $\rho(\cdot,\theta)$ for all $|\theta|<0.5$.
    \item $[0,1]$ and $[2,3]$ are linked in the optimal transportation from $\rho(\cdot,\theta)$ to $\rho^*$ for all $|\theta|<0.5$ except $\theta=0$.
    \item $\theta \mapsto W_p^p(\rho(\cdot,\theta),\rho^*)$ is differentiable for all $|\theta|<0.5$ except $\theta=0$.
\end{enumerate}
\end{proposition}
The proof can be found in~\Cref{sec:non_dif_non_ac}.

\subsection{Qualitative Error Analysis for the Gradient}
In this subsection, we prove that the almost-optimal solutions of Kantorovich's dual problem would provide accurate approximations of $\nabla f$.
\begin{proposition}\label{prp:grad_error}
Assume that \textbf{A1-A3} hold, and $f$ is G\^{a}teaux differentiable at $\theta_0 \in \Theta$. For every $\epsilon>0$ there exists a $\delta>0$ such that for all $(\phi,\psi) \in \Phi_c(\rho(\cdot,\theta_0),\rho^*)$ satisfying $I(\phi,\psi,\theta_0)>f(\theta_0)-\delta$ one has that
\[
\left|\nabla_\theta f(\theta_0)-\int_\Omega \phi^{cc}(x) \nabla_\theta \rho(x,\theta_0)dx\right|<\epsilon.
\]
\end{proposition}
The proof is presented in~\Cref{sec:grad_error}.

\begin{remark}
\Cref{prp:grad_error} asserts that one needs to calculate $c$-transforms of suboptimal $\phi$ for accurate gradients. This can be done very efficiently for costs of the form $c(x,y)=\sum_{i=1}^d h_i(x_i-y_i)$, where $h_i$ are even and strictly convex functions~\cite[Section 4.1]{jacobs20}. For OT algorithms that produce $c$-concave iterates, such as in \cite{jacobs20}, no further considerations are necessary.
\end{remark}

\section{Gradient Calculation} \label{gradient_descent}
Our parameter-dependent synthetic data obtained through the forward model is given by a finite-volume approximation
\begin{equation}\label{eq:rho_finite_volume}
    \rho(x,\theta)=\sum_{i=1}^n \rho_i(\theta) \frac{\chi_{C_i}(x)}{|C_i|},
\end{equation}
where $n = n_x n_y n_z$ is the total grid size, each $C_i$ is the finite volume cell, the parameter $\theta\in\Theta \subset\mathbb R^m$, and $\rho(\theta)=(\rho_i(\theta))_{i=1}^n$ is the solution to~\eqref{eq:main_direct} for some fixed $c,\epsilon>0$. Furthermore, after discretization, our reference data is given by
\[
\rho^*(y)=\sum_{i=1}^n \rho^*_i \frac{\chi_{C_i}(y)}{|C_i|}.
\]
By slightly abusing the notation we denote by $\rho^*=(\rho_i^*)_{i=1}^n$. Our goal is to solve
\begin{equation}\label{eq:f}
    \min_{\theta} f(\theta)=\T_c(\rho(\cdot,\theta),\rho^*)
\end{equation}
by gradient-based algorithms, where $\T_c$ is the optimal transport cost defined in \eqref{eq:T_c}. To apply \Cref{cor:f_Gat_suff}, which will guarantee the differentiability of $f$, we need to verify \textbf{A2} for \eqref{eq:rho_finite_volume} and that the connected components of $\operatorname{supp} \rho(\cdot,\theta)$ are linked according to \Cref{def:linked}. Since in all our experiments in~\Cref{results}, $\operatorname{supp} \rho(\cdot,\theta)=\cup_{i:\rho_i(\theta)>0} C_i$ is connected, the latter condition is satisfied. Therefore, we just need to verify \textbf{A2}, which is equivalent to the differentiability of $\theta \mapsto \rho(\theta)$. This verification is part of~\Cref{sec:implicit}.

Once all assumptions are verified, we have that
\[
    \nabla_\theta \rho(x,\theta)=\sum_{i=1}^n \nabla_\theta \rho_i(\theta) \frac{\chi_{C_i}(x)}{|C_i|}.
\]
Therefore, 
\begin{equation}\label{eq:grad_f_finite_volume}
    \nabla f(\theta )= \sum_{i=1}^n \nabla_\theta \rho_i(\theta) \phi_i(\theta),\quad \text{where }     \phi_i(\theta)=\frac{\int_{C_i} \phi(x,\theta) dx }{|C_i|}.    
\end{equation}
Here, $\phi(\cdot,\theta)$ is a Kantorovich potential for an OT from $\rho(\cdot,\theta)$ to $\rho^*$. Kantorovich potentials can be calculated by one of many available OT solvers such as~\cite{jacobs20,flamary2021pot}. Hence, we focus on calculating $\nabla_\theta \rho_i(\theta)$.

\subsection{Gradient Descent via Implicit Function Theorem} \label{sec:implicit}
First, we verify \textbf{A2}; that is, the differentiability of $\theta \mapsto \rho(\theta)$.
\begin{lemma}\label{lem:rho_diff}
Assume that $\theta \mapsto A(\theta),~\theta \in \Theta$ is a $C^1$ matrix valued function such that $A(\theta)$ is column stochastic with strictly positive entries for all $\theta \in \Theta$. Then the system of equations
\begin{equation}\label{eq:Arho=rho}
    A(\theta) \rho = \rho,\quad \rho \cdot \mathbf{1}=1,
\end{equation}
has a unique solution $\rho=\rho(\theta)$ for all $\theta \in \Theta$. Moreover, $\theta \mapsto \rho(\theta)$ is continuously differentiable with $\zeta_k(\theta)=\partial_{\theta_k} \rho(\theta)$ being the unique solution of
\begin{equation}\label{eq:(A-I)zeta}
    (A(\theta)-I) \zeta_k= -\partial_{\theta_k} A(\theta) \rho(\theta),\quad \zeta_k \cdot \mathbf{1}=0,
\end{equation}
where $\theta=(\theta_1,\theta_2,\cdots,\theta_m)$.
\end{lemma}
\begin{proof}
The existence and uniqueness of $\rho(\theta)$ is a consequence of the Perron--Frobenius Theorem as explained in \Cref{sec:finding_stationary}. Denote by $B(\theta)$ the matrix obtained from $A(\theta)-I$ by adding a $(n+1)$-st row vector $\mathbf{1}^\top$. Then we have that $\operatorname{ker}(B(\theta))=\{\mathbf{0}\}$, and so $\operatorname{rank}(B(\theta))=n$, and $n$ rows of $B(\theta)$ are linearly independent. Moreover, since $\operatorname{ker}(A(\theta)-I)=\operatorname{span}\{\rho(\theta)\}$, we have that $\operatorname{rank}(A(\theta)-I)=n-1$. Thus, the first $n$ rows of $B(\theta)$ are linearly dependent, and any list of $n$ independent rows must contain the last row $\mathbf{1}^\top$. Since $\theta \mapsto A(\theta)$ is continuous, linearly independent vectors stay so in a neighborhood of each $\theta$. Hence, we fix $\theta$ and without loss of generality assume that the rows of $B(\theta)$ from $2$ to $n+1$ are linearly independent in a neighborhood of $\theta$. 

Denote by
\[
    F(\theta,\rho)=\widetilde{B(\theta)} \rho-e_n,
\]
where $\widetilde{B(\theta)}$ is the matrix obtained from $B(\theta)$ by dropping the first row and $e_n$ is the $n$-th standard basis vector. Then we have that $\rho(\theta)$ is the unique solution of $F(\theta,\rho)=0$, and $D_\rho F(\theta,\rho)=\widetilde{B(\theta)}$ is non-degenerate. Thus, the Implicit Function Theorem applies and we obtain that $\theta \mapsto \rho(\theta)$ is continuously differentiable. Therefore, we can differentiate \eqref{eq:Arho=rho} and obtain \eqref{eq:(A-I)zeta}. Moreover, $\operatorname{ker}(B(\theta))=\{\mathbf{0}\}$ yields that the solution of \eqref{eq:(A-I)zeta} is unique.
\end{proof}

Applying \Cref{lem:rho_diff} to $A(\theta)=M_\epsilon(\theta)$ we obtain that the solution of \eqref{eq:main_direct} is differentiable and \eqref{eq:grad_f_finite_volume} holds. Thus, we can devise a gradient descent algorithm as follows:
\begin{equation}\label{eq:algo_discrete}
\begin{cases}
    M_{\epsilon}(\theta^l) \rho^l = \rho^l,\quad \rho^l \cdot \mathbf{1}=1,\\
    \\
    (M_{\epsilon}(\theta^l)-I) \zeta^l_k = -\partial_{\theta_k} M_{\epsilon}(\theta^l)\rho^l, \quad \zeta^l_k \cdot \mathbf{1}=0,\quad 1\leq k \leq m,\\
    \\
    (\phi^l,\psi^l) \in \underset{\phi_i+\psi_j \leq c(x_i,x_j)}{\operatorname{argmax}} [ \phi \cdot \rho^l+\psi \cdot \rho^* ],\\
    \\
    \theta_k^{l+1}=\theta_k^l-\tau^l~{\phi^l}\cdot \zeta_k^l, \quad 1\leq k \leq m.
\end{cases}
\end{equation}
where $\tau^l>0$ is a proper step size to for the gradient descent algorithm.

\subsection{Gradient Descent via Adjoint Method}
Here we discuss an alternative approach to calculate the gradient \eqref{eq:grad_f_finite_volume} via the adjoint-state method.

\begin{lemma}\label{lem:adjoint}
Assume that $\theta \mapsto A(\theta)$ satisfies the hypotheses in \Cref{lem:rho_diff}, $\rho(\theta)$ is the solution of \eqref{eq:Arho=rho}, and $\phi \in \R^n$ is an arbitrary vector. Then the linear system
\begin{equation}\label{eq:adjoint}
    (A(\theta)^\top - I) \lambda = -\phi+\phi \cdot \rho(\theta)~\mathbf{1}
\end{equation}
is consistent with a one-dimensional solution set. Moreover, for any solution $\lambda$ one has that
\[
    \partial_{\theta_k} (\phi \cdot \rho(\theta))=\lambda \cdot \partial_{\theta_k} A(\theta) \rho(\theta).
\]
\end{lemma}
\begin{proof}
Since $\operatorname{im}(A(\theta)^\top-I)=\operatorname{ker}(A(\theta)-I)^\perp$, we have to show that
\[
    -\phi+\phi \cdot \rho(\theta)~\mathbf{1} \in \operatorname{ker}(A(\theta)-I)^\perp=\operatorname{span}\{\rho(\theta)\}^\perp.
\]
A simple calculation yields the result:
\[
     (-\phi+\phi \cdot \rho(\theta)~\mathbf{1} )\cdot \rho(\theta)=-\phi \cdot \rho(\theta)+\phi \cdot \rho(\theta)~\mathbf{1}\cdot \rho(\theta)=0.
\]
Furthermore, since $\operatorname{ker}(A(\theta)^\top-I)=\operatorname{span}\{\mathbf{1}\}$, the solution set of \eqref{eq:adjoint} is a one-dimensional coset of $\operatorname{span}\{\mathbf{1}\}$.

Finally, assume that $\lambda$ is an arbitrary solution of \eqref{eq:adjoint}. Then applying \eqref{eq:(A-I)zeta} we obtain that
\begin{equation*}
    \begin{split}
        \partial_{\theta_k} (\phi \cdot \rho(\theta))=& \phi \cdot \zeta_k = (\phi \cdot \rho(\theta)~\mathbf{1}-(A(\theta)^\top -I)\lambda ) \cdot \zeta_k\\
        =& \phi \cdot \rho(\theta)~\mathbf{1}\cdot \zeta_k-\lambda \cdot (A(\theta)-I)\zeta_k=\lambda \cdot \partial_{\theta_k} A(\theta) \rho(\theta)
    \end{split}
\end{equation*}
\end{proof}

Applying \Cref{lem:adjoint} to $A(\theta)=M_\epsilon(\theta)$, we obtain an alternative, but equivalent, gradient descent algorithm: 
\begin{equation}\label{eq:algo_adjoint_discrete}
\begin{cases}
    M_{\epsilon}(\theta^l) \rho^l = \rho^l,\quad \rho^l \cdot \mathbf{1}=1,\\
    \\
    (\phi^l,\psi^l) \in \underset{\phi_i+\psi_j \leq c(x_i,x_j)}{\operatorname{argmax}} [ \phi \cdot \rho^l+\psi \cdot \rho^* ],\\
    \\
    (M_{\epsilon}(\theta^l)^\top-I) \lambda^l = -\phi^l+\phi^l \cdot \rho^l~\mathbf{1},\quad \lambda^l \cdot \mathbf{1}=0, \\
    \\
    \theta_k^{l+1}=\theta_k^l-\tau^l~\lambda^l \cdot \partial_{\theta_k} M_{\epsilon}(\theta^l) \rho^l,\quad 1\leq k \leq m.
\end{cases}
\end{equation}
Here, $\tau^l>0$ is a chosen step size to guarantee enough decrease in the objective function. Note that we add a condition $\lambda^l \cdot \mathbf{1}$ to ensure the uniqueness of $\lambda^l$.

We present a numerical scheme for efficiently solving systems of equations~\eqref{eq:algo_discrete} and \eqref{eq:algo_adjoint_discrete} in~\Cref{sec:numerical_scheme_1}.

\subsection{The Gradient of the $M_\epsilon(\theta)$}
For both algorithms \eqref{eq:algo_discrete} and \eqref{eq:algo_adjoint_discrete} we need to evaluate $\partial_{\theta_i} M_{\epsilon}(\theta)$. Denote by $H(x)=\frac{dx^+}{dx}$ the Heaviside function. We then have 
\[
    \partial_{\theta_i} v^+=H(v) \partial_{\theta_i} v,\quad \partial_{\theta_i} v^-=(1-H(v)) \partial_{\theta_i} v.
\]
We can also consider smoothed versions of $H$ such as
\[
    H_k(x)=\frac{d}{dx} k \log(1+e^{\frac{x}{k}})=\frac{e^{\frac{x}{k}}}{1+e^{\frac{x}{k}}}.
\]
It is not hard to show that $H_k$ is smooth and $\lim_{k\to 0^+} H_k(x)=H(x)$. Based on \eqref{eq:M_eps}, we derive that
\begin{equation*}
\begin{split}
    \partial_{\theta_i} M_{\epsilon}=(1-\epsilon)c \cdot \partial_{\theta_i} K_{mat}=(1-\epsilon)c \cdot \sum_{i_d=1}^d \frac{\Delta t}{\Delta x_{(i_d)}} \partial_{\theta_i} K_{(i_d)}(\theta)
\end{split}
\end{equation*}
where each matrix $\partial_{\theta_i} K_{(i_d)}(\theta)$ has three nonzero diagonals for each pair of $(i, i_d)$ where $1\leq i \leq m,\, 1\leq i_d \leq d$, while the offsets of the diagonals depend on $i_d$, as we have discussed earlier regarding~\eqref{eq:Kd}. We emphasize that $\partial_{\theta_i} K_{(i_d)}(\theta) $  shares the same tridiagonal structure with $K_{(i_d)}(\theta)$ for each $i_d$ as illustrated below.
\[
\partial_{\theta_i} K_{(i_d)}(\theta) = %
 {\scriptsize \begin{bmatrix}
 \ddots & & &  & \\
  &\ddots & &  & \\
 \ddots & & -\left(1-H_k(v_{(i_d-\frac{1}{2})})\right) \partial_{\theta_i} v_{(i_d-\frac{1}{2})}  &  & \\
 & \ddots &\vdots& \ddots &\\
  \ddots  & &\left(1-H_k(v_{(i_d-\frac{1}{2})})\right) \partial_{\theta_i} v_{(i_d-\frac{1}{2})} - H_k(v_{(i_d+\frac{1}{2})}) \partial_{\theta_i} v_{(i_d+\frac{1}{2})} &   & \ddots\\
 &   \ddots & \vdots  & \ddots &\\
 &  & H_k(v_{(i_d+\frac{1}{2})}) \partial_{\theta_i} v_{(i_d+\frac{1}{2})}  &   & \ddots\\
 & & & \ddots&  \\
 & & & & \ddots
\end{bmatrix}. }
\]
One can also compute $\partial_{\theta_i} K_{(i_d)}(\theta)$ through automatic differentiation; see~\Cref{sec:AD} for details of implementation and performance comparison.

\section{Numerical Results}\label{results}
In this section, we show several numerical results on dynamical system parameter identification, following the methodology described in the earlier sections. The forward problem is to solve for the steady state of the corresponding PDE~\eqref{eq:continuity} rather than the ODE system~\eqref{eq: dynamicalSystem}. The objective function that compares the observed and the synthetic invariant measures is the quadratic Wasserstein metric ($W_2$) from optimal transportation. The optimization algorithm implemented for all inversion tests is the gradient descent method with backtracking line search to control the step size~\cite{nocedal2006numerical}.

\subsection{Chaotic System Examples} 
\begin{figure}
\centering
    \includegraphics[width=0.24\textwidth]{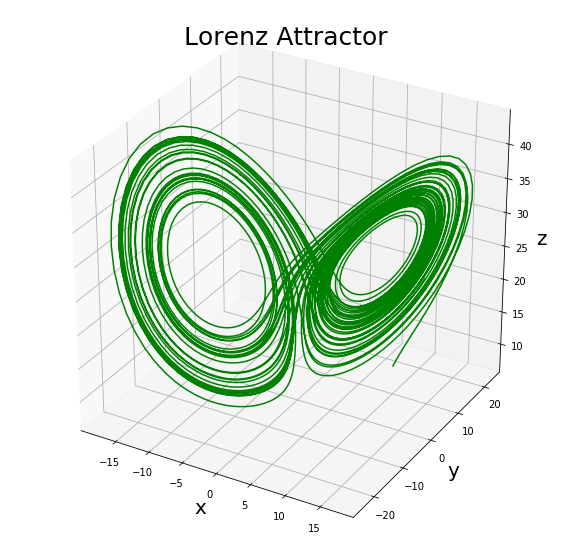}
    \includegraphics[width=0.24\textwidth]{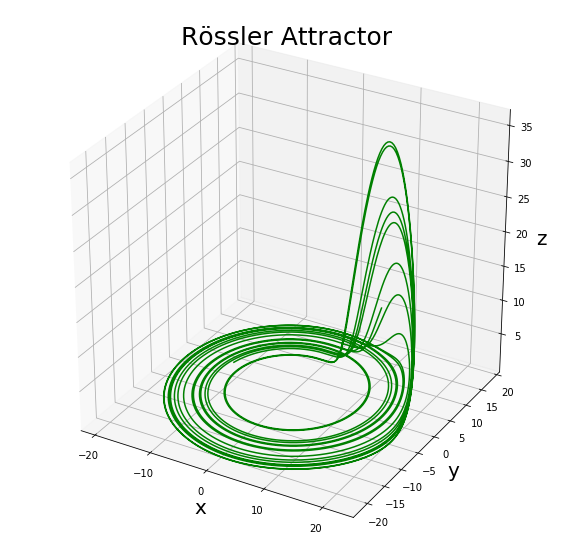} 
    \includegraphics[width=0.24\textwidth]{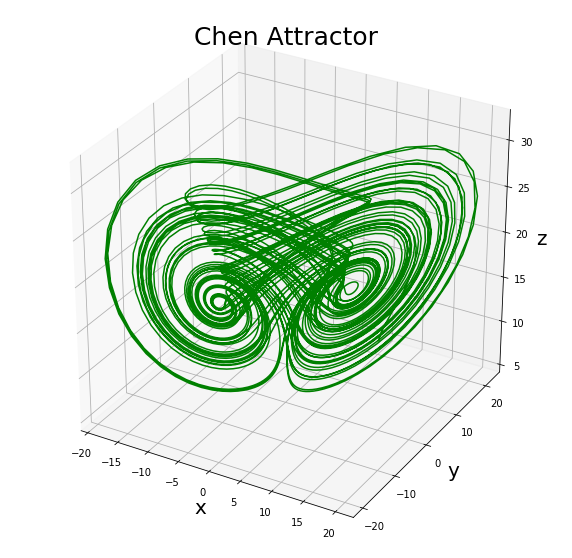} 
   \includegraphics[width=0.24\textwidth]{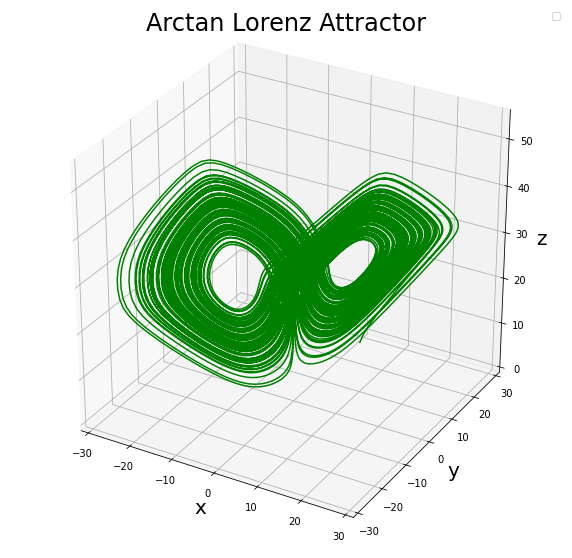} 
    \caption{From left to right: the Lorenz, R\"ossler, Chen, and Arctan Lorenz attractors.}\label{fig:attractors}
\end{figure}

We test our proposed method on three classic chaotic systems: the Lorenz, R\"ossler, and Chen systems. These models are widely used benchmarks that illustrate typical features of dynamical systems with instabilities and nonlinearities that give rise to deterministic chaos. We also perform an inversion test on a modified Arctan Lorenz system in which the unknown parameters are nonlinear with respect to the flow velocity in terms of monomial basis. The true parameters are selected such that the dynamical systems exhibit chaotic behaviors; see the illustration through the time trajectories in~\Cref{fig:attractors}.

\subsubsection{Lorenz System}
Consider the following Lorenz system.
\begin{equation}\label{eq:Lorenz}
\begin{cases}
    \dot{x} &= \sigma (y-x),\\
    \dot{y} &= x(\rho -z)-y,\\
    \dot{z} &= x y -\beta z. 
\end{cases}
\end{equation}
The equations form a simplified mathematical model for atmospheric convection, where $x,y,z$ denote variables proportional to convective intensity, horizontal and vertical temperature differences. The parameters $\sigma, \beta, \rho$ are proportional to the Prandtl number, Rayleigh number, and a geometric factor.
The true parameter values that we will try to infer are $\sigma = 10,\; \beta=8/3,\;\rho=28$. These are well-known parameter values for which Lorenz system shows a chaotic behaviour. 

\subsubsection{R\"ossler System}
Consider the following R\"ossler System.
\begin{equation}\label{eq:Rossler}
\begin{cases}
    \dot{x} &= -y - z, \\
    \dot{y} &= x+a y, \\
    \dot{z} &= b + z(x-c). 
\end{cases}
\end{equation}
Here $x, y, z$ denote variables, while $a, b, c$ are the parameters we want to infer.
The system exhibits continuous-time chaos and is described by the above three coupled ODEs. The Rössler attractor behaves similarly to the Lorenz attractor, but it is easier to analyze qualitatively since it generates a chaotic attractor having a single lobe rather than two. The true parameters that we try to infer are $a = 0.1,\; b=0.1,\; c=14$. 

\subsubsection{Chen System}
Consider the following Chen System~\cite{chen1999yet}.
\begin{equation} \label{eq:chen}
\begin{cases}
    \dot{x} &= a(y-x), \\
    \dot{y} &= (c-a)x-xz+cy, \\
    \dot{z} &= xy -bz. 
\end{cases}
\end{equation}
Again, $x, y, z$ are variables and $a, b, c$ are parameters we will infer. 
The system has a double-scroll chaotic attractor, which is often observed from a physical, electronic chaotic circuit. The true parameters that we will infer are $a = 40,\; b=3,\; c=28$.

\subsubsection{Arctan Lorenz System}
The parameters in the earlier examples are all coefficients of the monomial basis. Here, we modify the right-hand side of the Lorenz system~\eqref{eq:Lorenz} to create a new dynamical system such that the particle flow velocity is nonlinear with respect to the monomial basis.
\begin{equation} \label{eq:arctan_lorenz}
\begin{cases}
    \dot{x} &= 50 \arctan \left( \sigma(y-x)/50 \right) , \\
    \dot{y} &= 50 \arctan \left(x(\rho-z)/50 - y/50 \right), \\
    \dot{z} &= 50 \arctan \left( (xy-\beta z)/50 \right). 
\end{cases}
\end{equation}
Again, $x, y, z$ are variables, and $\sigma, \rho, \beta$ are parameters we want to infer. The reference values are set to be $(10,28, 8/3)$, the same as the original Lorenz system.

\subsection{The Invariant Measures} \label{subsec:Invariant}
Here, we follow the numerical scheme described in~\Cref{sec:finding_stationary} and approximate the invariant measure through the regularized PDE surrogate model, represented by the corresponding probability density function (PDF), for the three dynamical systems at the given sets of parameters.

We compare PDFs obtained through the steady-state solution to~\eqref{eq:continuity} with the histogram accumulated from long-time trajectories from Direct Numerical Simulation (DNS). That is, we solve systems~\eqref{eq:Lorenz}--\eqref{eq:chen} forward in time using the explicit Euler scheme with time step $\Delta t$ from $t=0$ to its final time $t=T$. We then compute the physical invariant measure following~\eqref{eq:our_data}. Moreover, we use time trajectories that are enforced with either the intrinsic or the extrinsic noises.

\begin{figure}
    \centering
\subfloat[Steady-state solution to~\eqref{eq:continuity}]{
       \includegraphics[width = 0.3\textwidth]{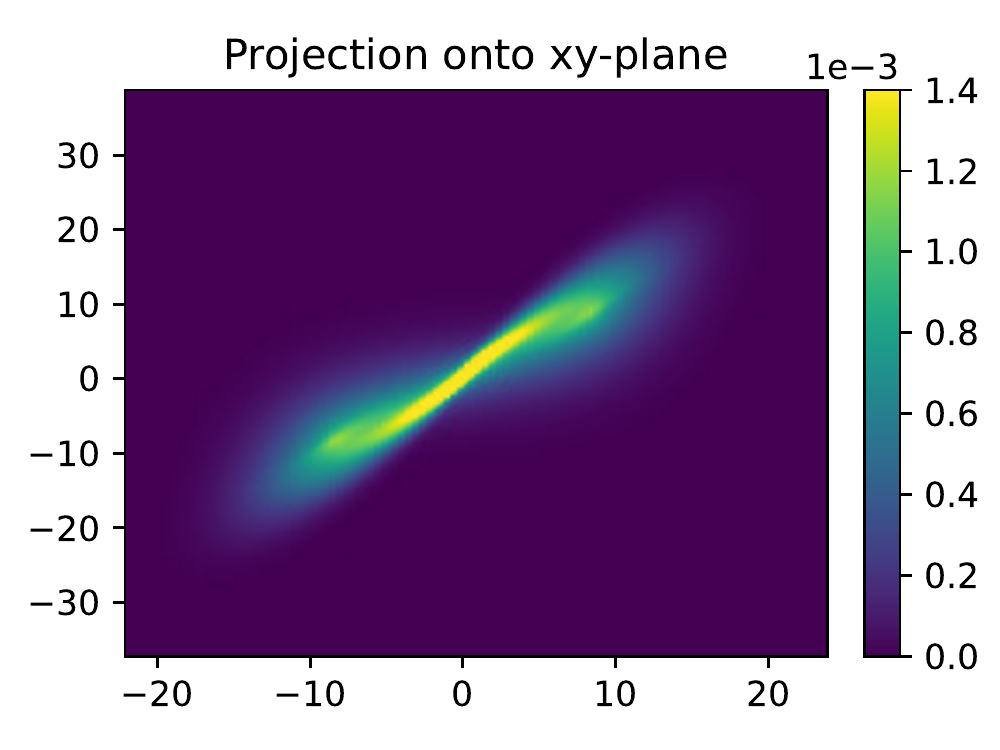}
       \includegraphics[width = 0.3\textwidth]{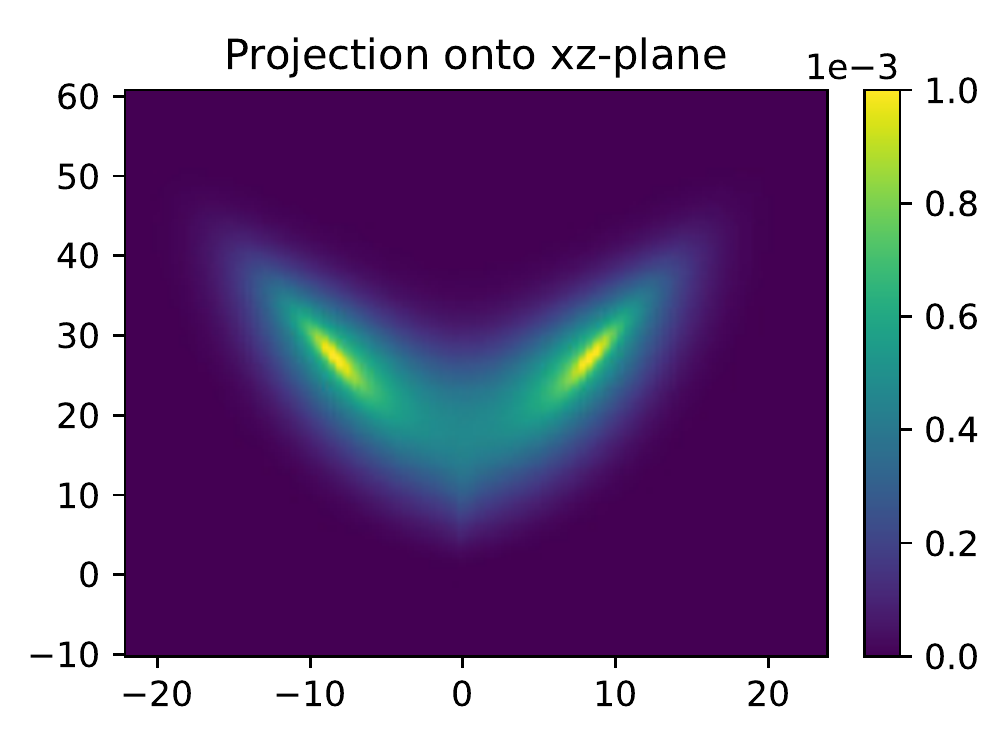}
       \includegraphics[width = 0.3\textwidth]{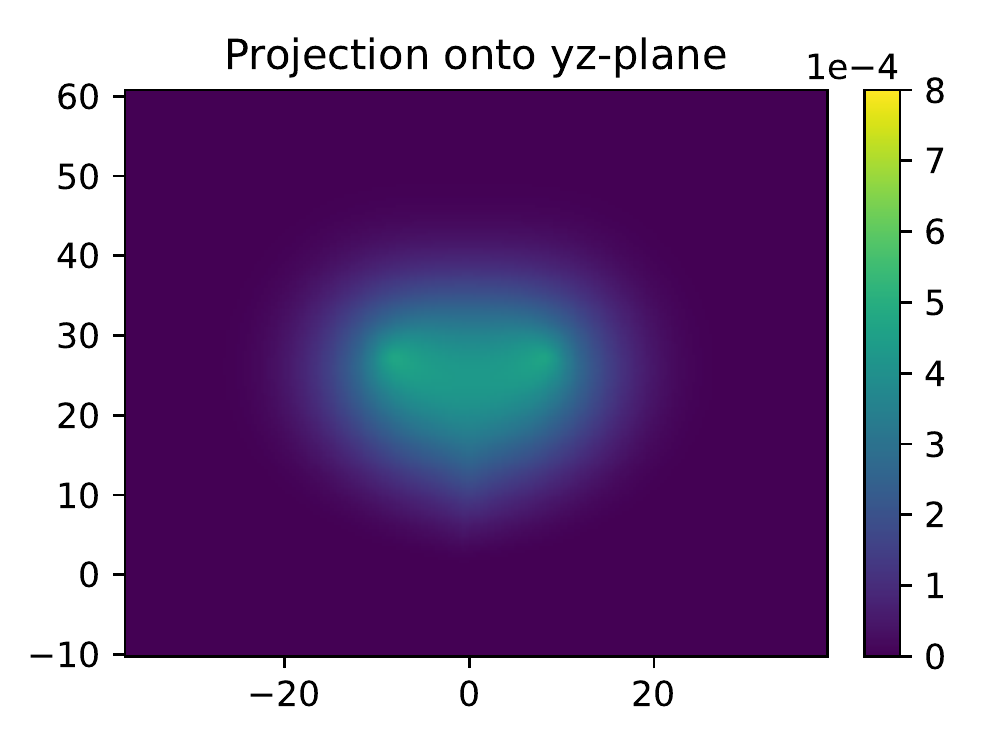} \label{fig:Lorenz_FPE_pdf}}\\
\subfloat[Histogram accumulated from noise-free Lorenz system time trajectory]{ \includegraphics[width = 0.3\textwidth]{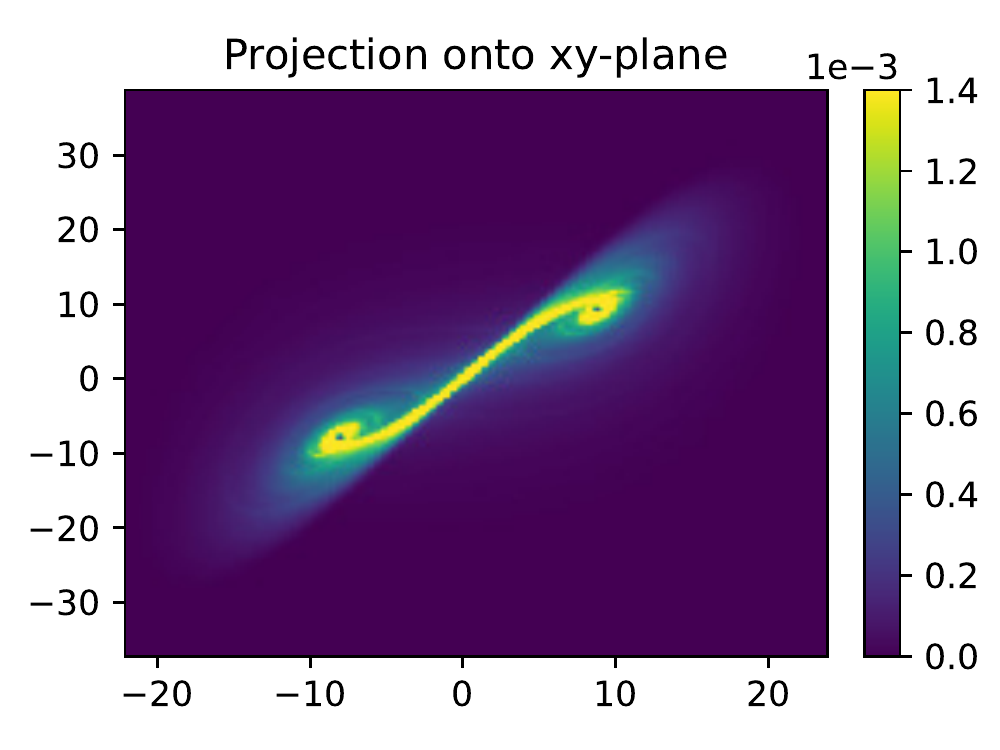}
        \includegraphics[width = 0.3\textwidth]{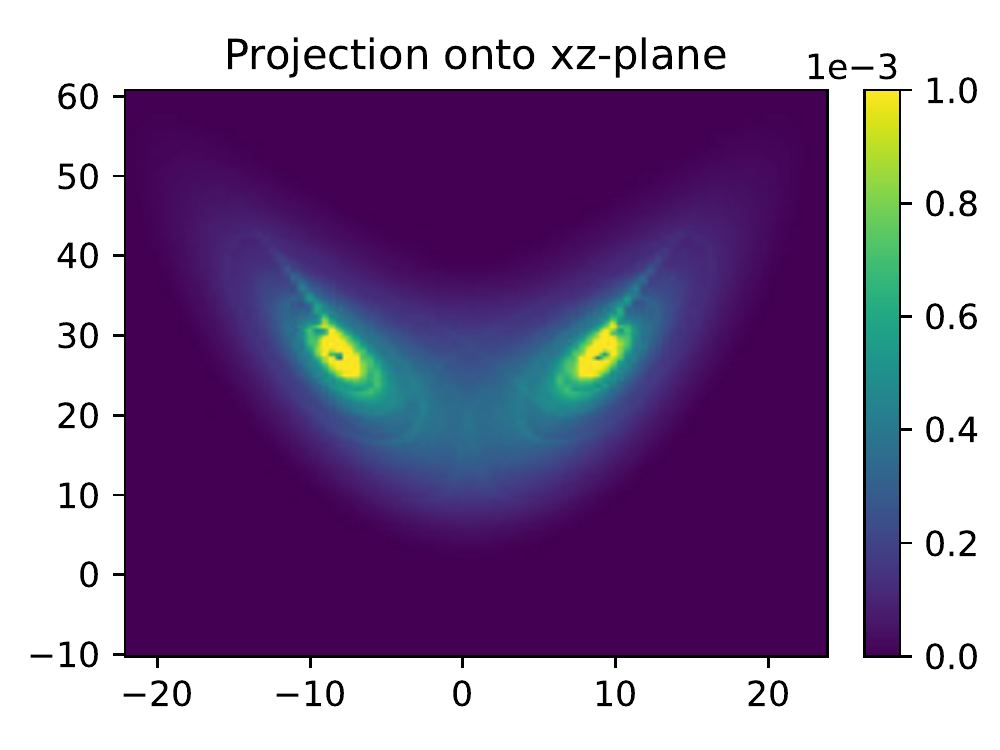}
            \includegraphics[width = 0.3\textwidth]{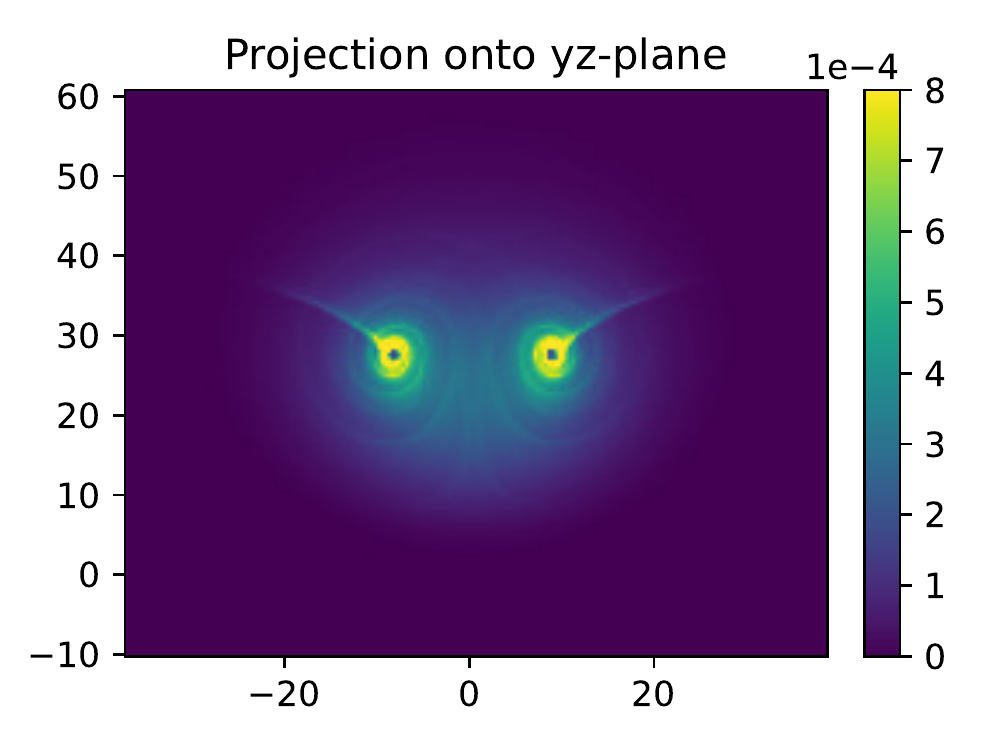} \label{fig:Lorenz_DNS_Clear_pdf}}\\
\subfloat[Histogram accumulated from Lorenz system time trajectory with \textit{intrinsic} noise]{
                \includegraphics[width = 0.3\textwidth]{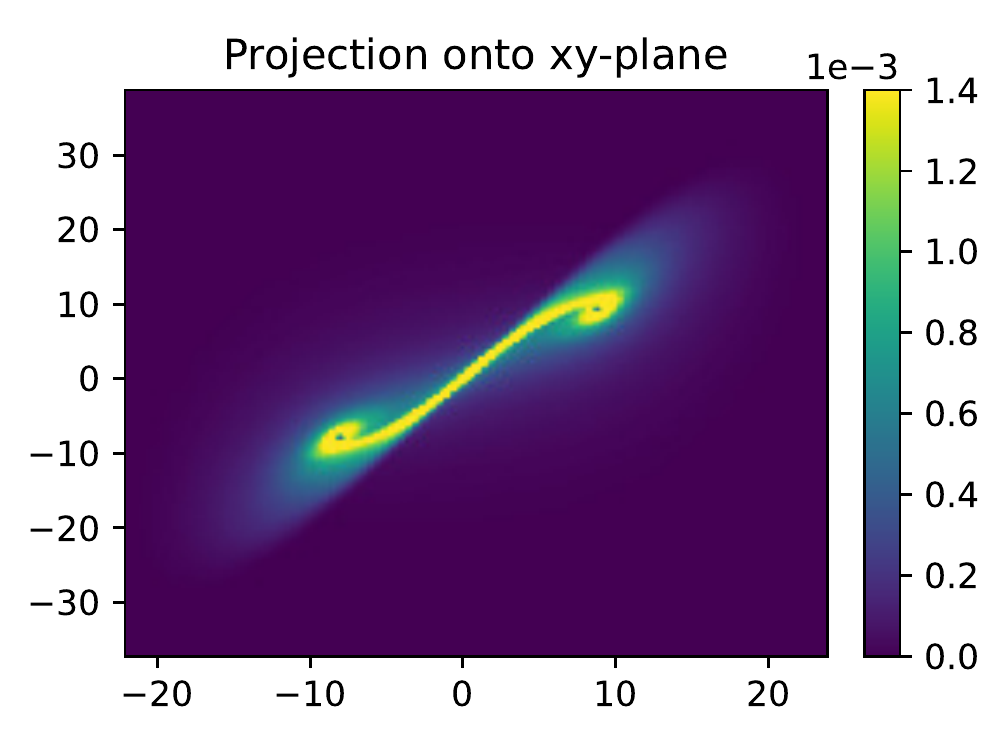}
        \includegraphics[width = 0.3\textwidth]{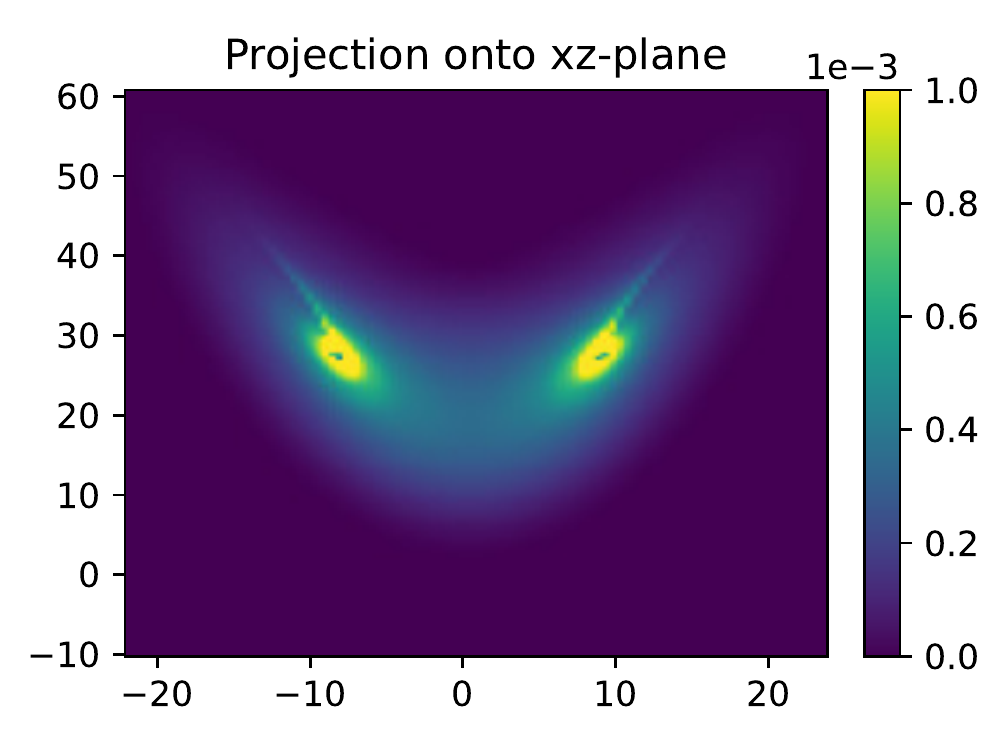}
            \includegraphics[width = 0.3\textwidth]{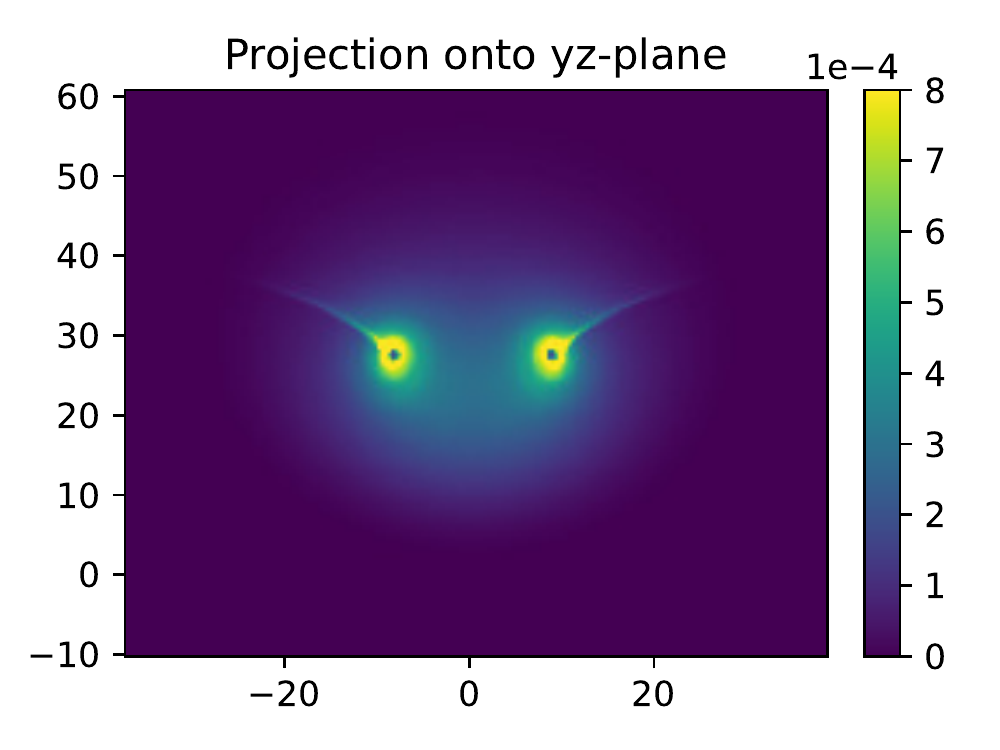} \label{fig:Lorenz_DNS_intrinsic_pdf}}\\
\subfloat[Histogram accumulated from Lorenz system time trajectory with \textit{extrinsic} noise]{
                \includegraphics[width = 0.3\textwidth]{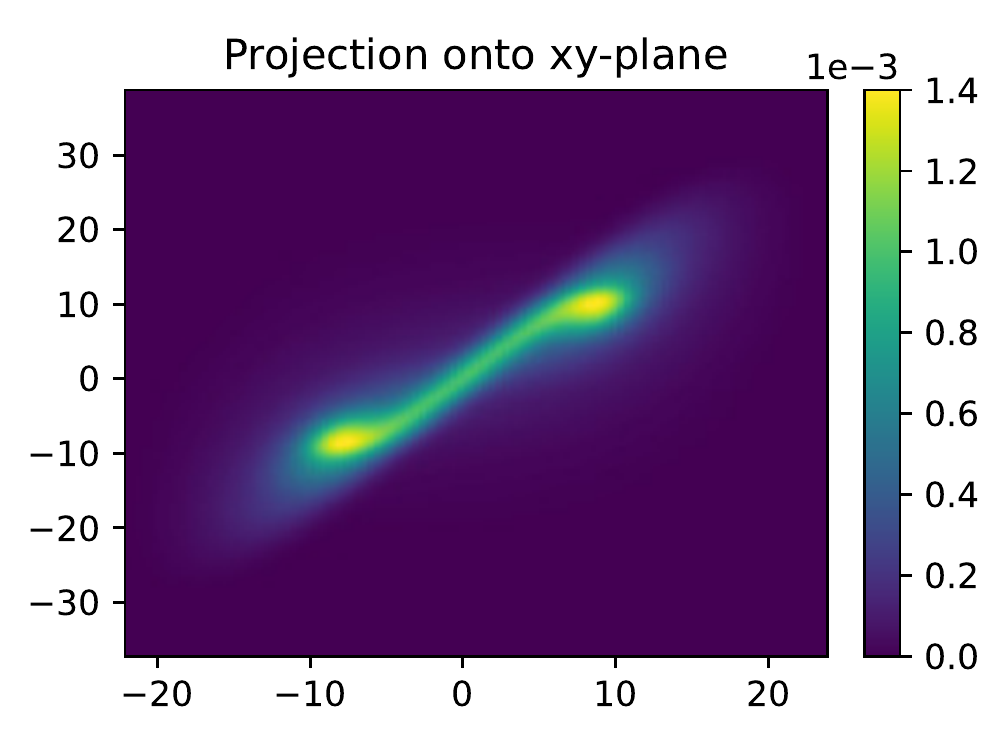}
        \includegraphics[width = 0.3\textwidth]{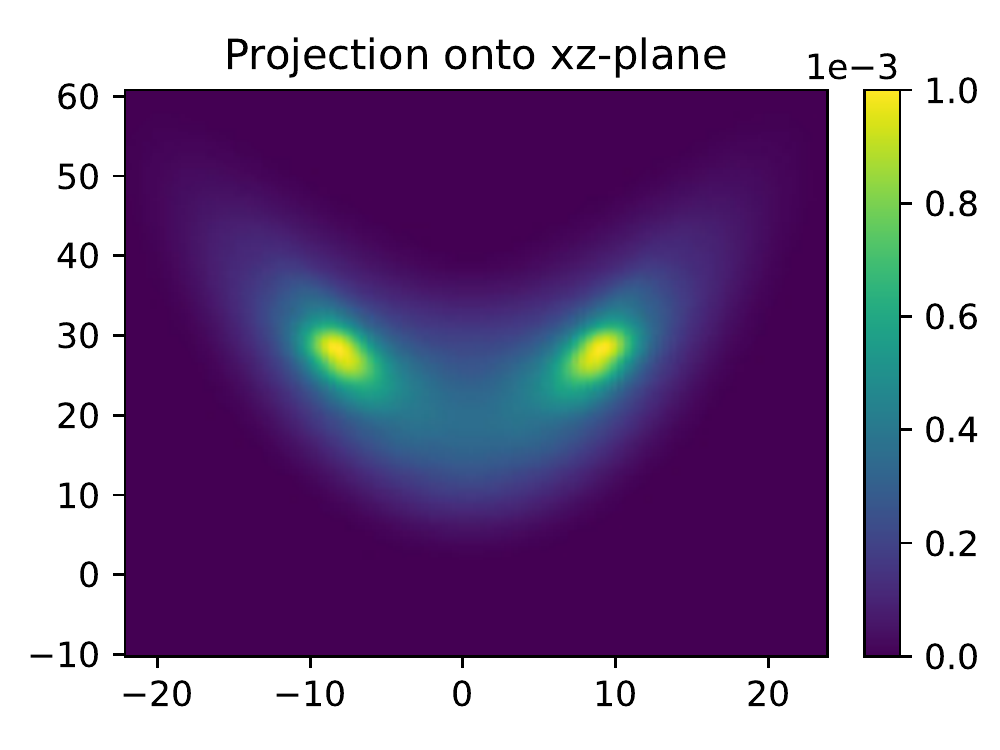}
            \includegraphics[width = 0.3\textwidth]{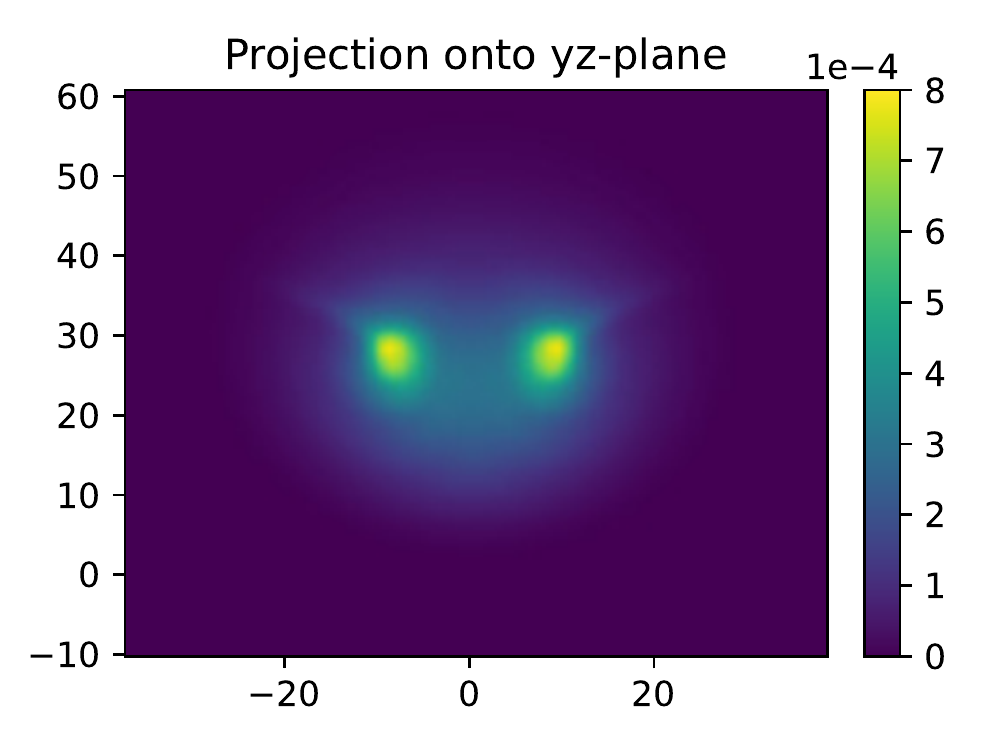}\label{fig:Lorenz_DNS_extrinsic_pdf}}
    \caption{Lorenz system. Top row: the steady state on the grid size $93\times 153\times 143$ by solving~\eqref{eq:continuity}. The teleportation parameter is $\epsilon = 10^{-6}$. Second row: projections of physical invariant measure from noise-free time trajectory for $T = 2\times 10^6$. Third row: projections of physical invariant measure from time trajectory with \textit{intrinsic} noise ${\omega}\sim \mathcal N(0,\mathbf I)$. Last row: projections of physical invariant measure from time trajectory with \textit{extrinsic} noise $\gamma\sim \mathcal N(0,\mathbf I)$.}
    \label{fig:Lorenz_pdf}
\end{figure}

\subsubsection{Numerical Illustrations}\label{sec:invariant_measure}
Comparisons for the Lorenz system~\eqref{eq:Lorenz} are displayed in \Cref{fig:Lorenz_pdf}. The three plots in the top row show the $x$--$y$, $x$--$z$, and $y$--$z$ projections of the dominant eigenvector of the Markov matrix $M_\epsilon$. The grid size for the finite volume discretization of~\eqref{eq:continuity} is $93\times 153\times 143$. The teleportation parameter is $\epsilon = 10^{-6}$. In the second row, we see the corresponding three projections of the physical invariant measure from \textit{noise-free} time trajectory for total time $T = 2\times 10^6$. The third row and the  bottom row show three projections of the physical invariant measure from time trajectories of the same total time $T$ but with \textit{intrinsic} noise ${\omega}\sim \mathcal N(0,\mathbf I)$ (the noise occurs on the right-hand side of the dynamical system as $\dot{\bf{x}}=v(\bf x)+ \omega$) and \textit{extrinsic} noise $\gamma\sim \mathcal N(0,\mathbf I)$ (the observation of the time trajectory suffers from noise as ${\bf x_\gamma}={\bf x}+{\gamma}$), respectively. The bin size for all three histograms is a cube of volume $0.5^3$.

\begin{figure}
\centering
\subfloat[Steady-state solution to~\eqref{eq:continuity}]{
\includegraphics[width = 0.3\textwidth]{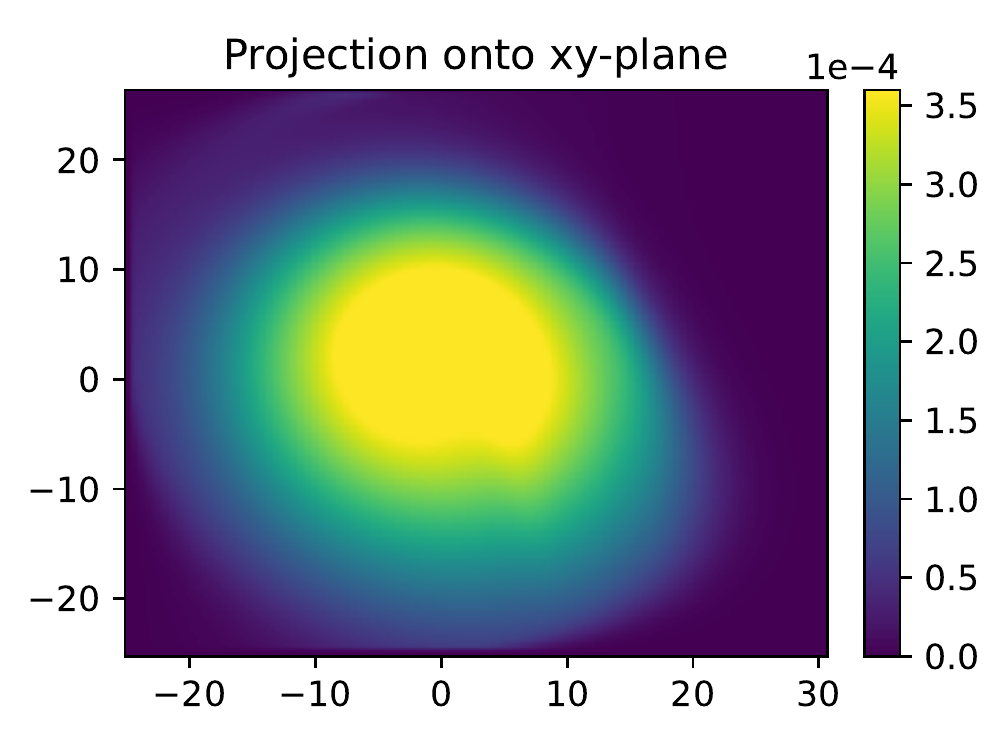}
\includegraphics[width = 0.3\textwidth]{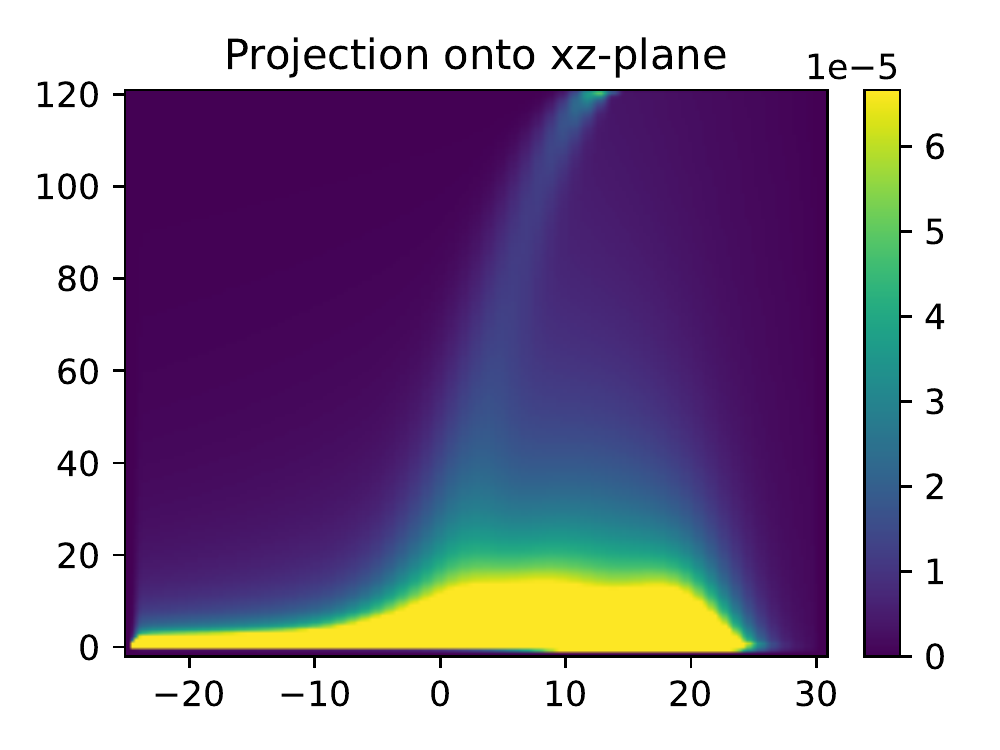}
\includegraphics[width = 0.3\textwidth]{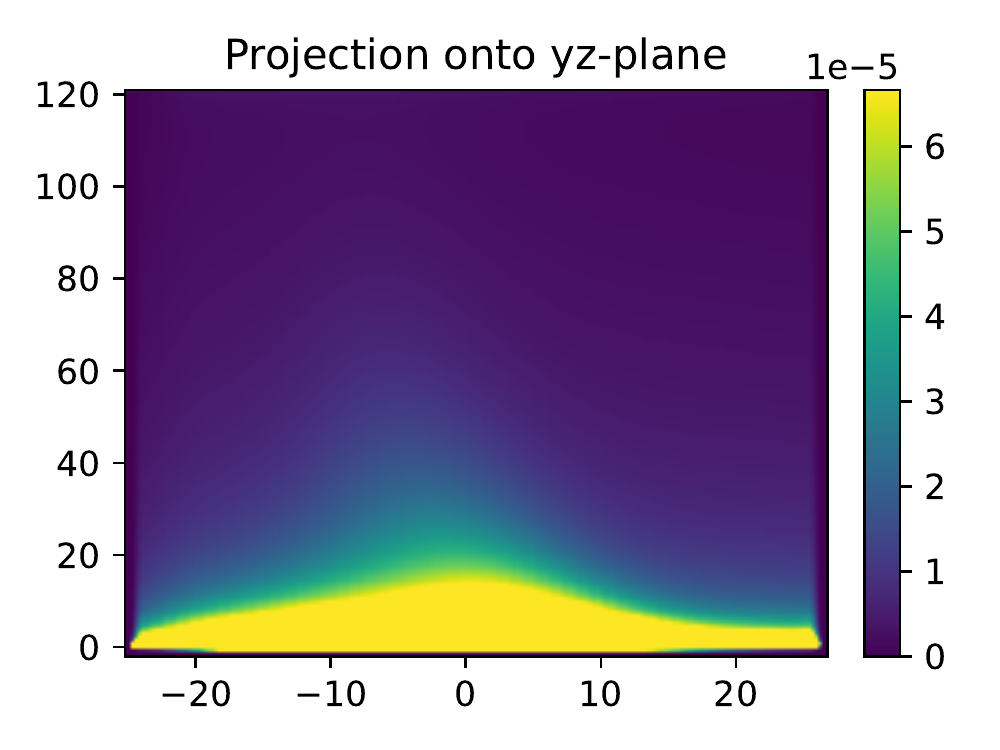} \label{fig:Rossler_FPE_pdf}}\\
\subfloat[Histogram accumulated from R\"ossler system time trajectory with intrinsic noise]{
\includegraphics[width = 0.3\textwidth]{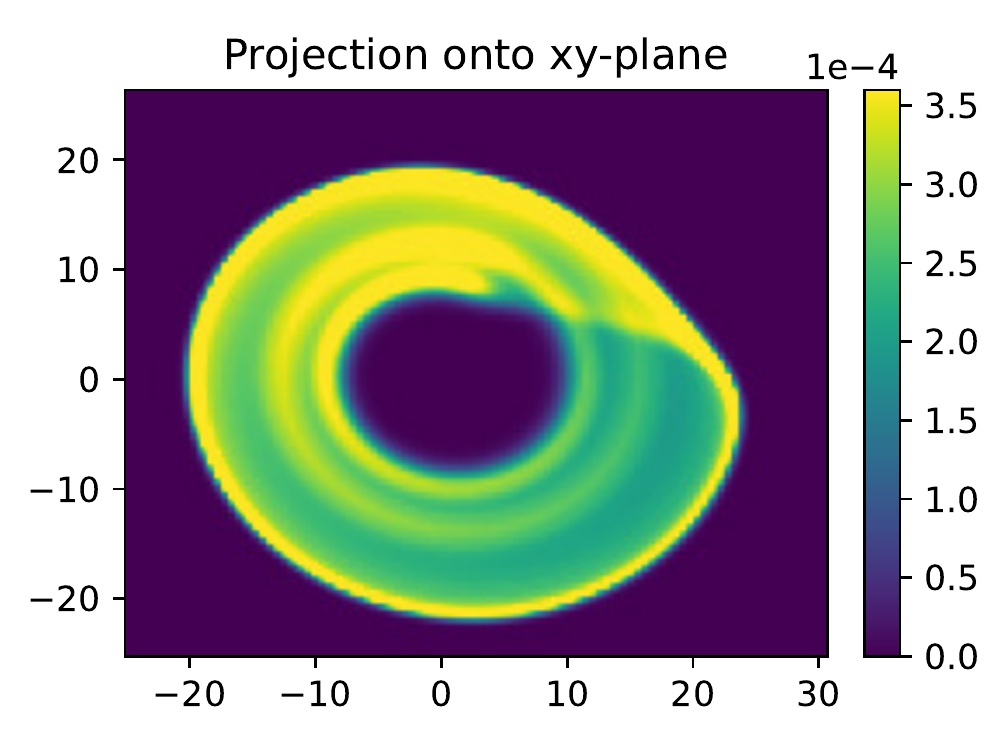}
\includegraphics[width = 0.3\textwidth]{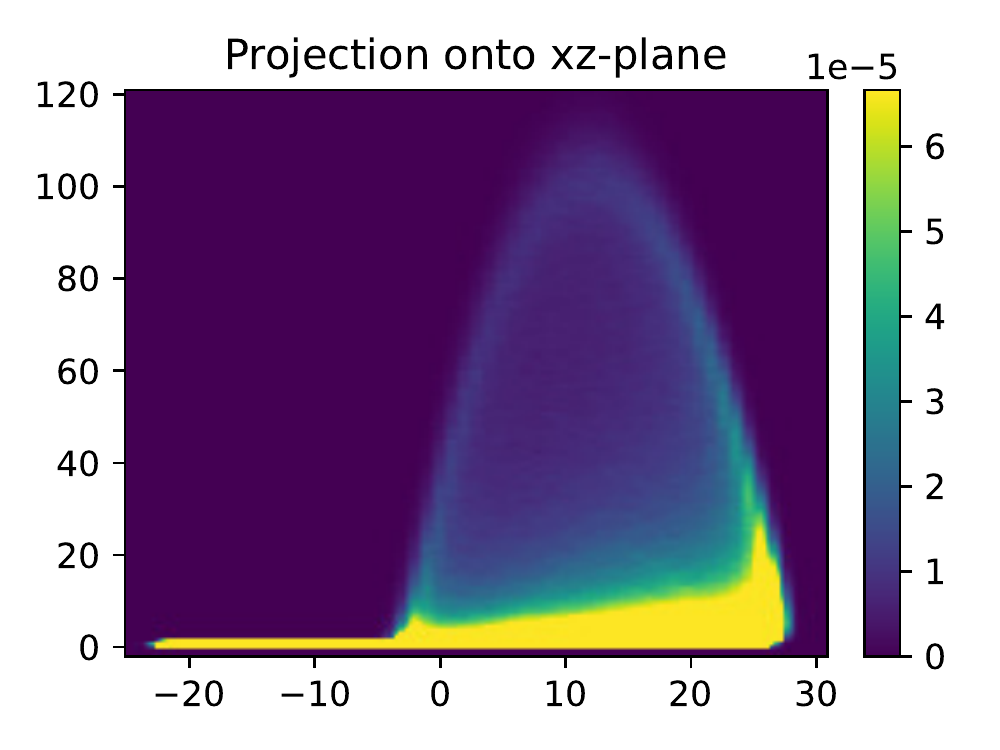}
\includegraphics[width = 0.3\textwidth]{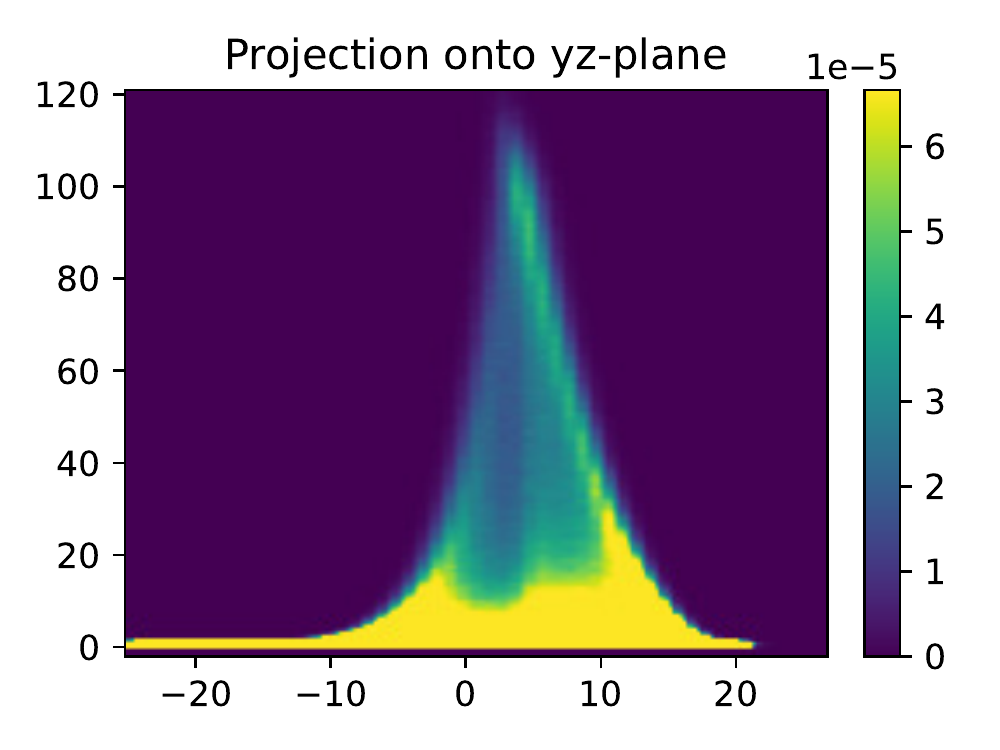} }
\caption{R\"ossler system. Top row: the steady-state solution to~\eqref{eq:continuity} on the grid size $94\times 87 \times 106$. The teleportation parameter is $\epsilon = 10^{-6}$. Bottom row: the  histogram accumulated from R\"ossler system time trajectory for total time $T = 1\times 10^6$ with intrinsic noise $\omega \sim \mathcal N(0,0.2 \mathbf I)$.}
\label{fig:Rossler_pdf}
\end{figure}

Similar plots for the R\"ossler system~\eqref{eq:Rossler} are presented in~\Cref{fig:Rossler_pdf}. Top row shows the steady-state solution to~\eqref{eq:continuity} computed on a grid size is $94\times 87 \times 106$. The teleportation parameter is $\epsilon = 10^{-6}$. For the bottom row, the R\"ossler system time trajectory runs for a total time $T = 1\times 10^6$ with an \text{intrinsic} noise ${\omega}\sim \mathcal N(0,0.2 \mathbf I)$. The bin size for the histogram is a cube of volume $0.6^3$.

\begin{figure}
 \centering
 \subfloat[Steady-state solution to~\eqref{eq:continuity}]{
\includegraphics[width = 0.3\textwidth]{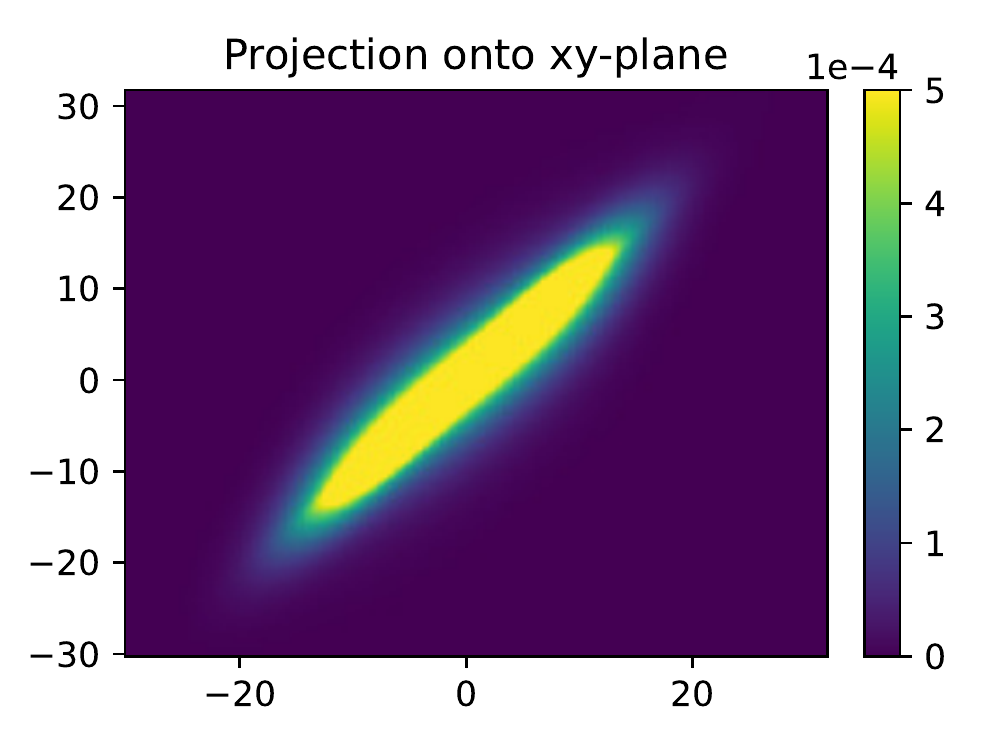}
\includegraphics[width = 0.3\textwidth]{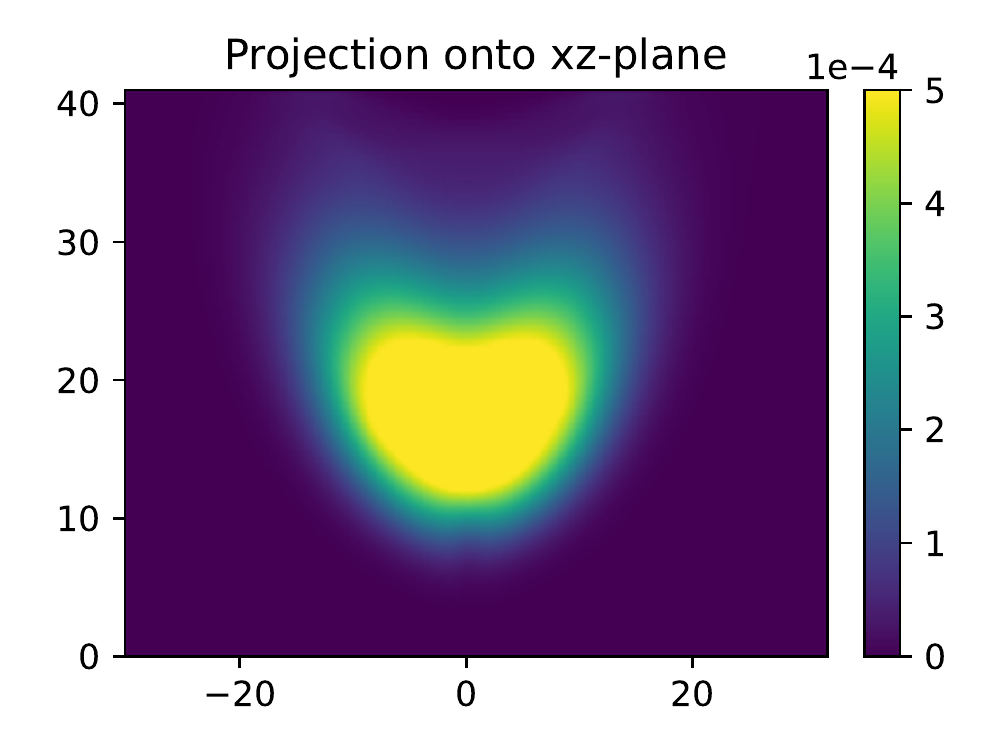}
\includegraphics[width = 0.3\textwidth]{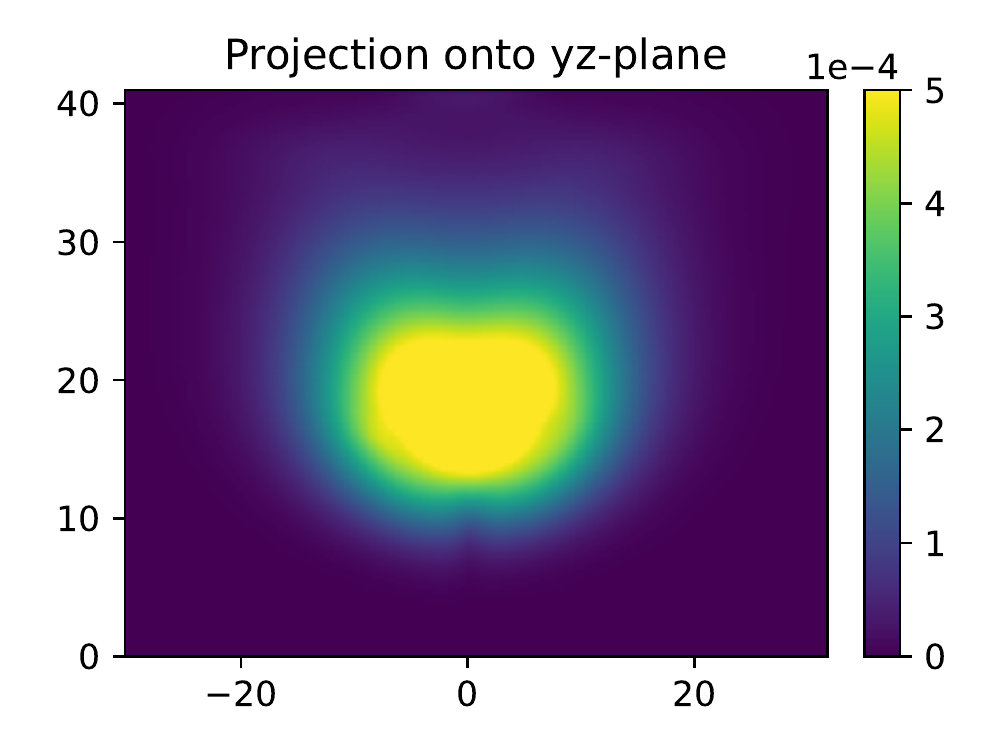}
\label{fig:Chen_FPE_pdf}}\\
\subfloat[Histogram accumulated from Chen system time trajectory with intrinsic noise]{
 \includegraphics[width = 0.3\textwidth]{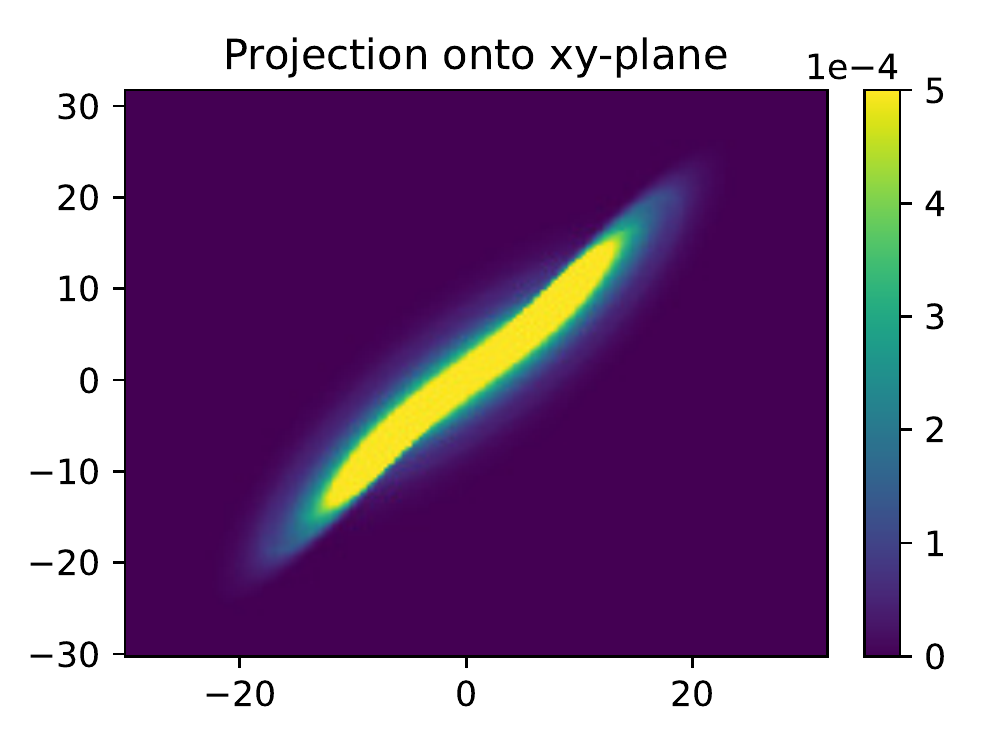}
\includegraphics[width = 0.3\textwidth]{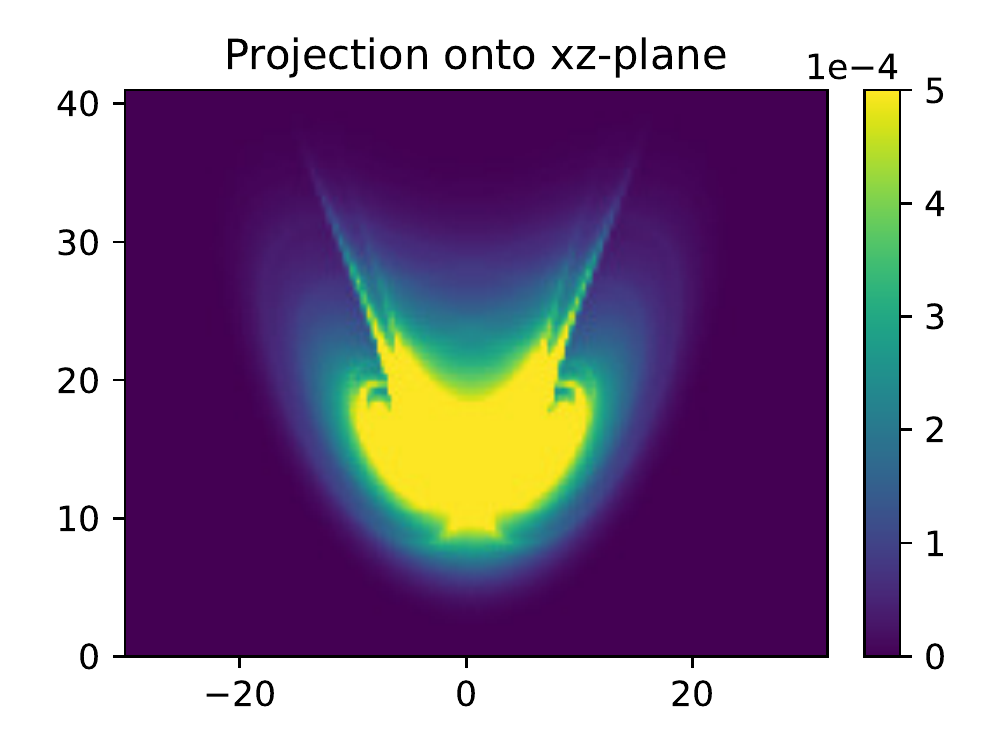}
\includegraphics[width = 0.3\textwidth]{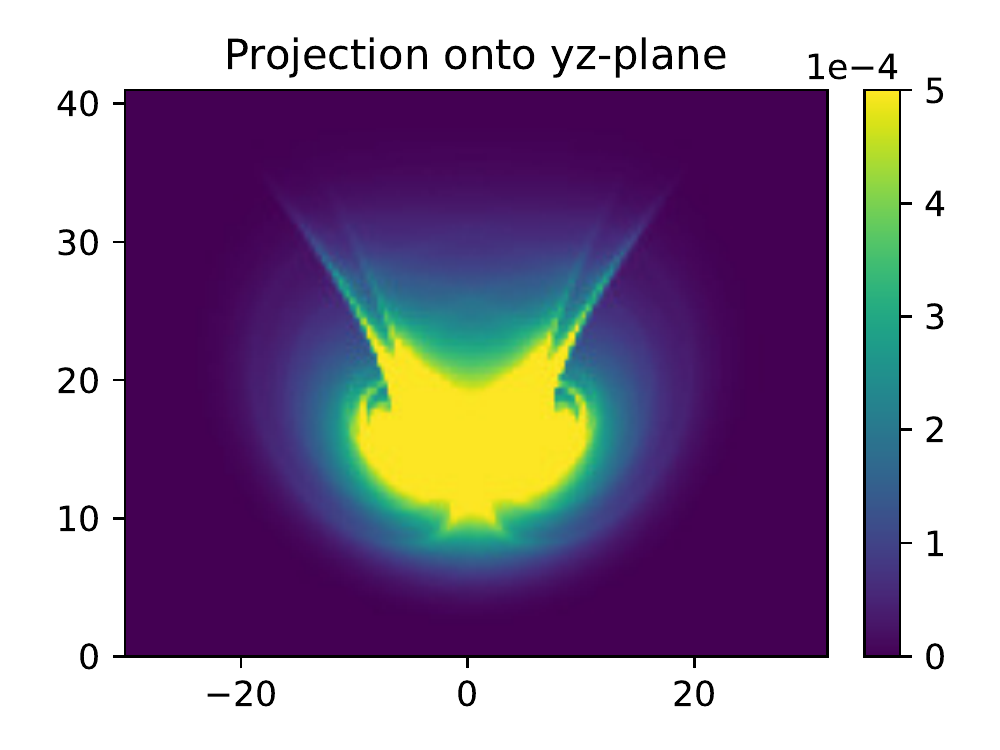}}
\caption{Chen system. Top row: the steady-state solution to~\eqref{eq:continuity} on the grid size $125\times 125 \times 83$. The teleportation parameter is $\epsilon = 10^{-6}$. Bottom row: the histogram accumulated from Chen system time trajectory with $T = 5\times 10^5$ and intrinsic noise $\omega\sim \mathcal N(0,0.2 \mathbf I)$.}
    \label{fig:Chen_pdf}
\end{figure}

\Cref{fig:Chen_pdf} shows the comparisons for the Chen system~\eqref{eq:chen}. The first row displays the three projections of the steady-state solution to~\eqref{eq:continuity} on a $104\times 104 \times 69$ grid. The teleportation parameter is $\epsilon = 10^{-6}$. The bottom row shows the projections of the physical invariant measure accumulated from time trajectory with intrinsic noise for a total time $T = 5\times 10^5$. The bin size for the histogram is a cube of volume $0.5^3$. The intrinsic noise ${\omega}\sim\mathcal N(0,0.2 \mathbf I)$.

\subsubsection{The Effect of Noise}
\label{sec:invariant_measure_noise}
It is important to understand the fundamental limitations and challenges of converging the low-order solver for~\eqref{eq:continuity}, particularly the role that the addition of the \textit{extrinsic} and \textit{intrinsic} noises play here as an approximation of the diffusive errors expected in the PDE solver.  

After the ODE is solved, the \textit{extrinsic} noise applied to the trajectory corresponds to an effective Gaussian blur of the DNS results. In the limit of long time DNS simulation, the true density is the result of taking every point on the invariant measure, represented by a delta function in state space based on the DNS solution, and then replacing it with a Gaussian ball of equal integral mass with width defined by the standard deviation of the noise. This process is equivalent to the Gaussian blur common in image processing.

The \textit{intrinsic} noise case is more complicated. Since the three examples we have all admit non-trivial basins of attraction, the accumulation of energy resulting from the addition of noise is balanced by the dissipation inherent to the dynamics for directions that are orthogonal to the attractor manifold.  While the \textit{extrinsic} noise corresponds to a spatially uniform low pass filter, the blurring resulting from the \textit{intrinsic} noise depends more on the local stability of the attractor in state space.



\subsubsection{The Effect of Mesh Size and Numerical Diffusion}\label{sec:invariant_measure_dx}
\begin{figure}
    \centering
    \includegraphics[width=0.95\textwidth]{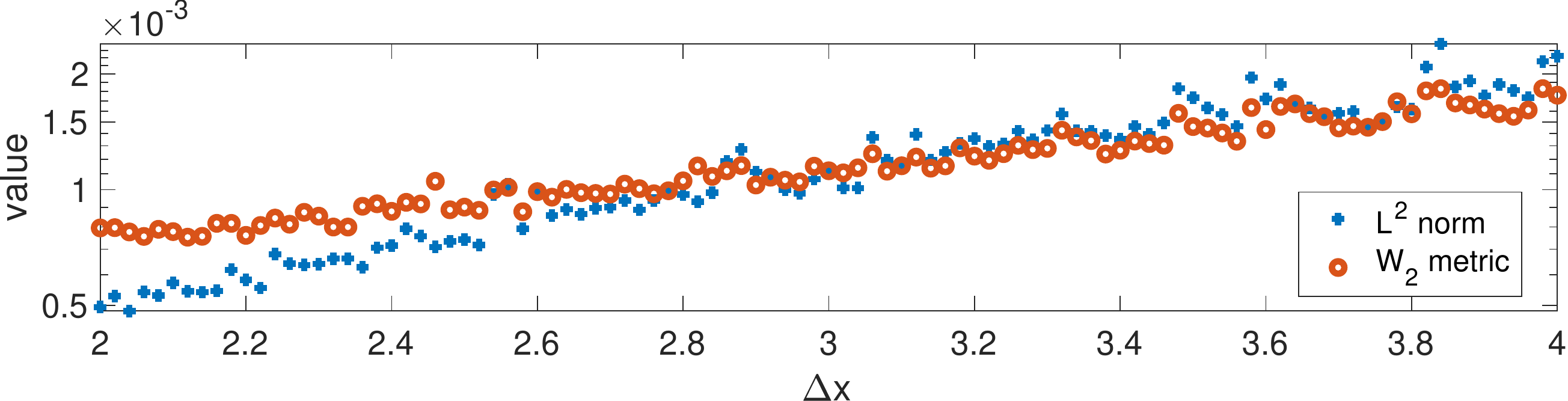}
    \caption{The $W_2$ metric and the $L^2$ difference between the PDF accumulated from DNS with bin volume $(\Delta x)^3$ and the PDF solved as the steady-state solution to~\eqref{eq:continuity} with spatial spacing $\Delta x$. The PDFs are for the Lorenz system at the true parameters.} \label{fig:model-discrepancy}
\end{figure}

While of a form dominated by diffusion, numerical errors of the PDE solver have a dependence on the flow velocity $\propto v^2\Delta t$, as described in~\cite{bewley2012efficient}. This is the well-known numerical diffusion that motivates running computational fluid dynamics solvers with a Courant--Freidrich--Lewy (CFL) condition  number as close to $1$ as possible for low-order methods to minimize the numerical diffusivity.  While in this work, we seek a steady-state solution, the time step of the forward operator has effectively been selected to comply with this CFL restriction in the act of ensuring that the forward operator is at least positive semi-definite in~\eqref{eq:CFL}. Substituting the CFL restriction, $\Delta t=\Delta x/v_{max}$, into the expression for the numerical diffusivity, it can be seen that numerical diffusion in the PDE solver is effectively $\propto v^2\Delta x/ v_{max}$, which is bounded by $v_{max}\Delta x$, suggesting first-order convergence with $\Delta x$ if $v_{max}$ is bounded. More detailed numerical analysis for the convergence and numerical errors can be found in~\cite{leveque2002finite}. The linear convergence is also seen in~\Cref{fig:model-discrepancy}, where we compare the differences between the PDF accumulated from the Lorenz system DNS with again $T = 2\times 10^6$ and the steady-state solution to~\eqref{eq:continuity}, both evaluated at the true parameters for the Lorenz system. The histogram bin size changes as we use different $\Delta x$ in the finite volume discretization. 

However, as a steady-state problem, the error due to this numerical diffusion, similar to the \textit{intrinsic} noise added to the DNS solution, accumulates until balanced by the dissipation of the attractor dynamics is balanced.  This, too, depends on how dissipative the basin is.

We remark that all the inversion tests in this paper use $\Delta x = 3$. It is for demonstration only and thus far from being optimal. The size of the Markov matrix $M$ grows $\propto \Delta x^{-3}$ as $\Delta x$ decreases, making it very expensive to compute the steady state at a fine mesh. Mesh-refinement strategies could help provide better parameter estimates while saving computational costs of the forward solve. This, along with more efficient numerical implementations, will be left to future work.

\subsubsection{The Effect of Random Samples}\label{sec:invariant_measure_sample}
\begin{figure}
    \centering
\subfloat[$10^4$ samples of the trajectory]{
\includegraphics[width = 0.48\textwidth]{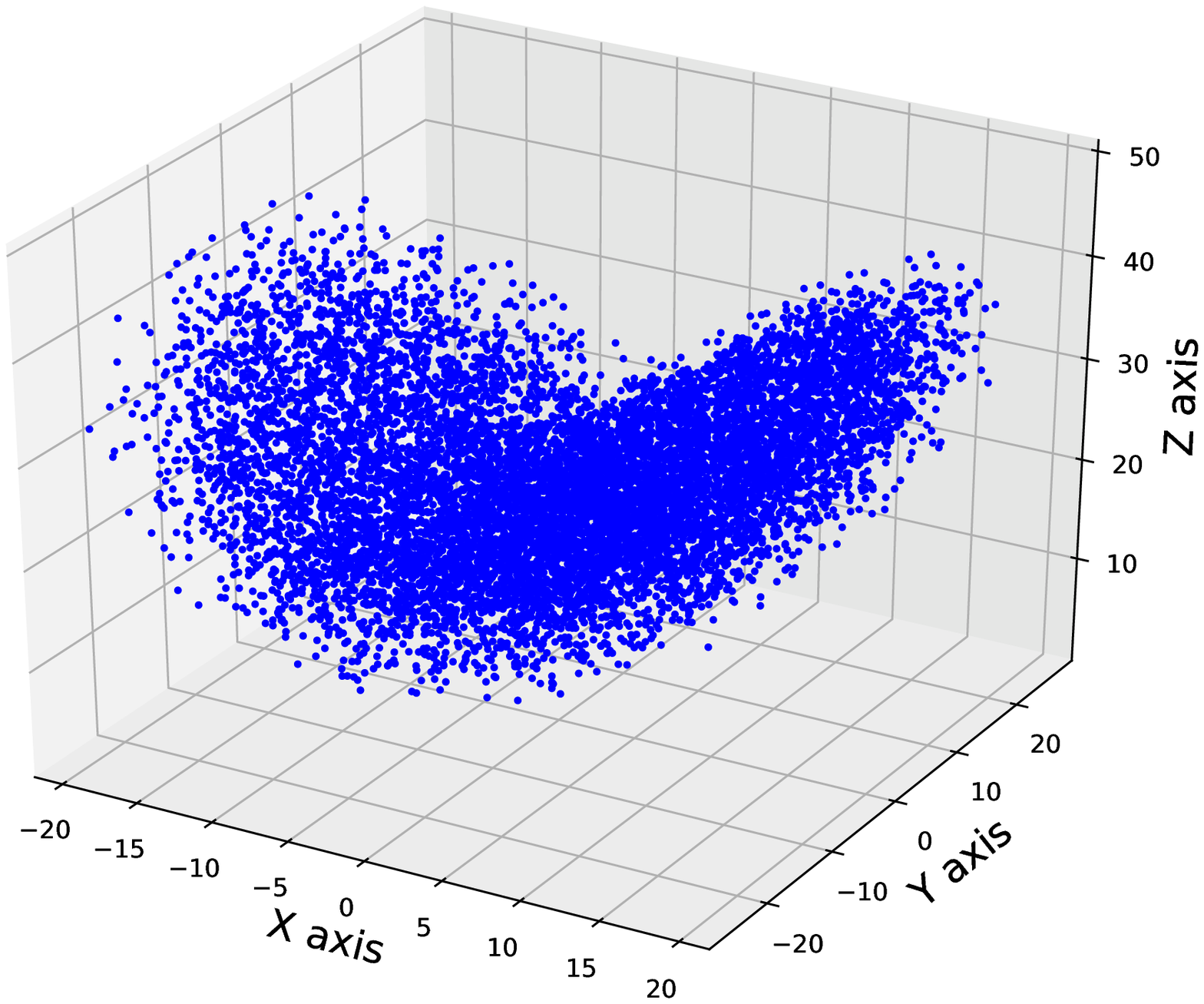} \label{fig:subsample_all}}
\subfloat[Zoom-in view of the first $10$ samples]{\includegraphics[width = 0.48\textwidth]{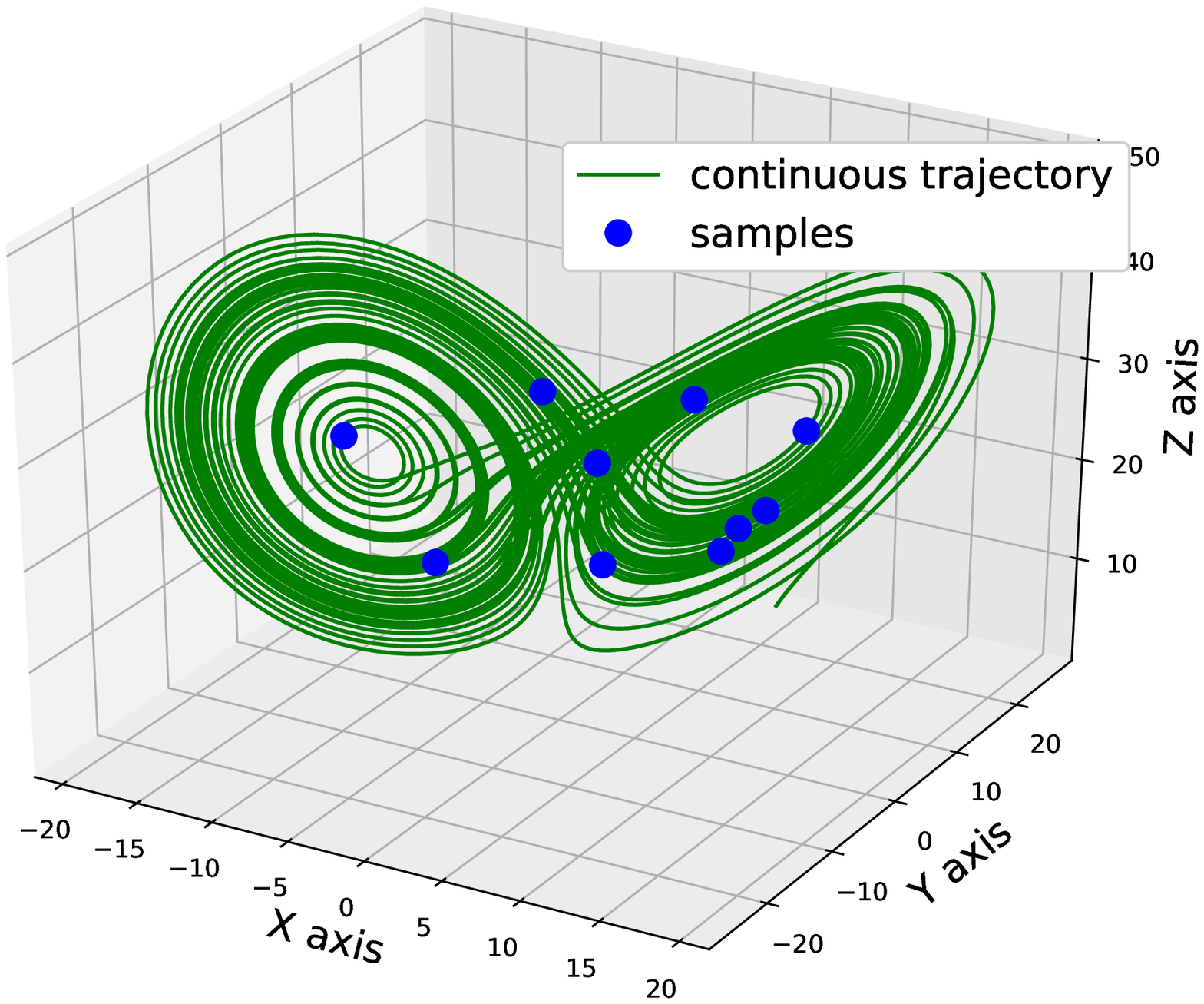} \label{fig:subsample_sparse}}
    \caption{Left: $10^4$ random samples of the Lorenz trajectory; Right: illustration of the first $10$ samples of~\Cref{fig:subsample_all} compared with the continuous trajectory.\label{fig:Lorenz_subsample}}
\end{figure}

One main advantage of the proposed framework is that we allow the trajectory data to be ``slowly'' sampled, in which case we do not have access to the state-space velocity or velocity estimates, i.e., the $\dot{\bf x}$. In~\Cref{fig:subsample_all}, we illustrate the total samples of the trajectory that will be used in the parameter inference, while~\Cref{fig:subsample_sparse} displays the relationship of the first $10$ samples in the time series with the continuous trajectory in the corresponding time window. One can observe that our random samples of state-space positions are ``sparse'' and could not accurately estimate the state-space velocity. Later in~\Cref{sec:real_inversion}, we use the reference measure constructed from such slowly sampled and completely randomized state measurements to perform parameter identification.

\begin{figure}
    \centering
\subfloat[DNS vs.~DNS]{
\includegraphics[width = 0.48\textwidth]{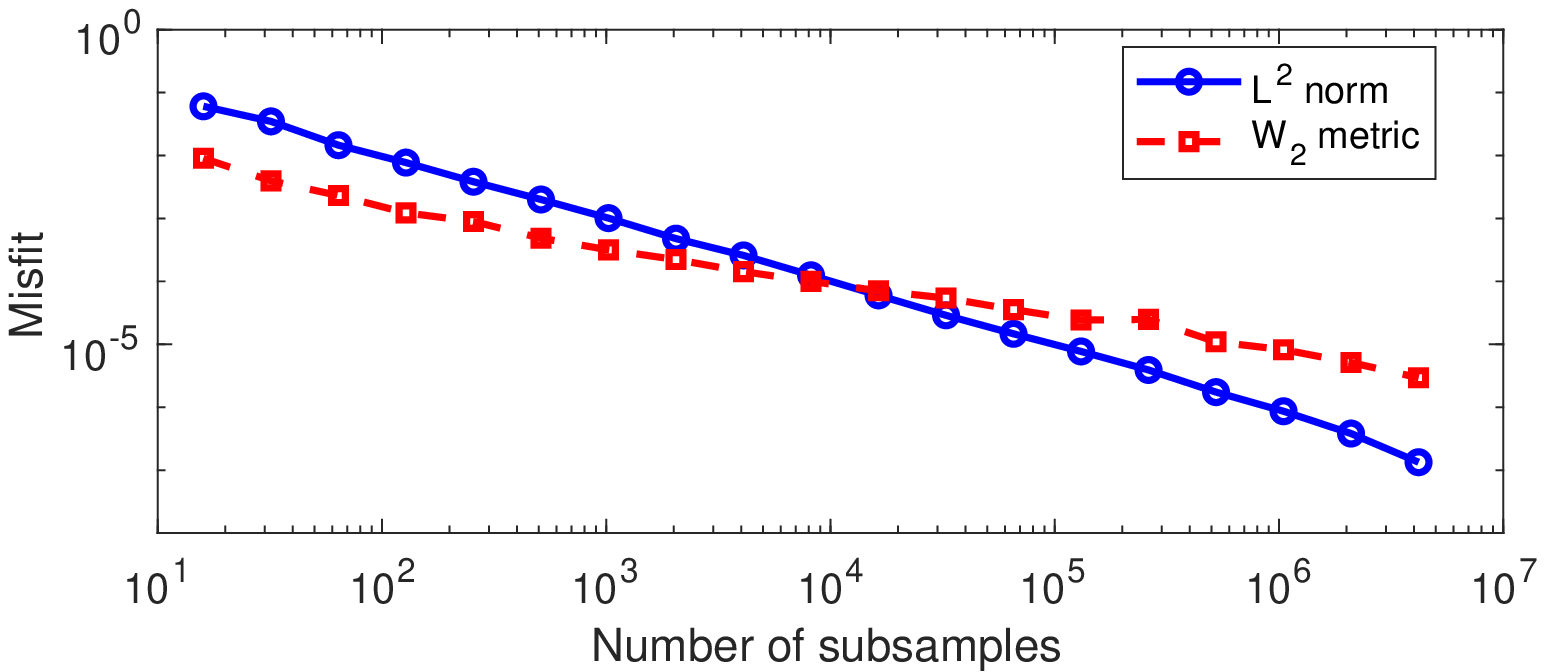} \label{fig:subsample_DNSvsDNS}}
\subfloat[DNS vs.~PDE steady state]{
\includegraphics[width = 0.48\textwidth]{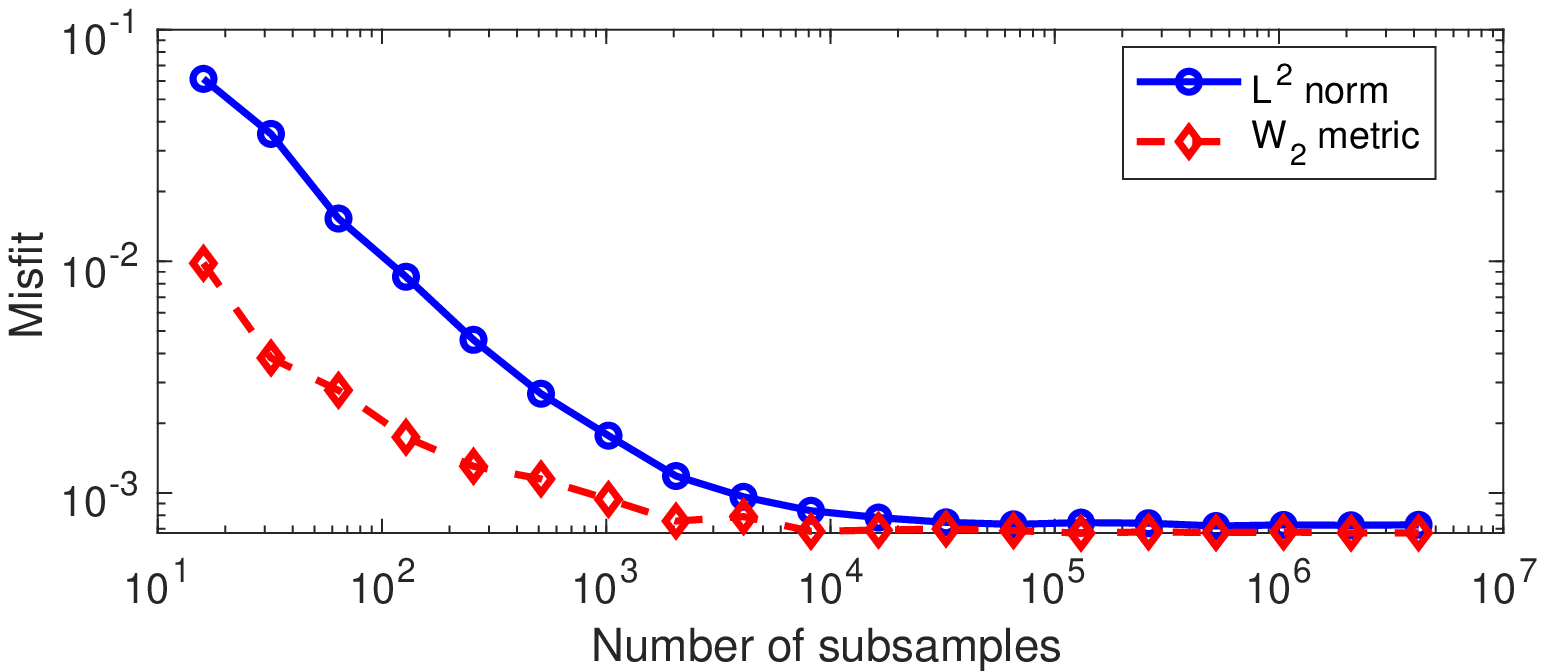} \label{fig:subsample_DNSvsFPE}}
\caption{Left: the misfit between the density accumulated from sub-sampled data and the one from the entire trajectory; Right: the misfit between the density accumulated from sub-sampled data and the steady-state solution from the FPE solver.\label{fig:Lorenz_subsample_compare}}
\end{figure}

In~\Cref{fig:Lorenz_subsample_compare}, we numerically investigate the relationship between the amount of state-space position samples and the approximation error for the invariant measure. In~\Cref{fig:subsample_DNSvsDNS}, we set the reference density to be the histogram accumulated from $10^8$ samples and compare it with the histogram accumulated from much fewer samples. We observe the classical Monte Carlo error, $\mathcal{O}(1/\sqrt{N})$, where $N$ is the number of samples. In~\Cref{fig:subsample_DNSvsFPE}, we change the reference density to the steady-state solution from the FPE solver. The error plateaus for large $N$ since the modeling error, mainly due to the numerical diffusion discussed in~\Cref{sec:invariant_measure_dx}, becomes the dominant factor of the mismatch when $N$ is large enough. It also indicates that we do not need too many trajectory samples to perform parameter identification.

\subsubsection{The Effect of Teleportation Parameter}\label{sec:invariant_measure_tele}
To obtain the steady-state solution, we used the so-called teleportation trick to regularize the Markov matrix; see~\Cref{subsec:finite_volume} for details. Here, we numerically investigate the impact of the teleportation parameter $\epsilon$ on the obtained steady-state solution.

\begin{figure}
    \centering
\subfloat[PDE steady state vs.~PDE steady state]{
\includegraphics[width = 0.48\textwidth]{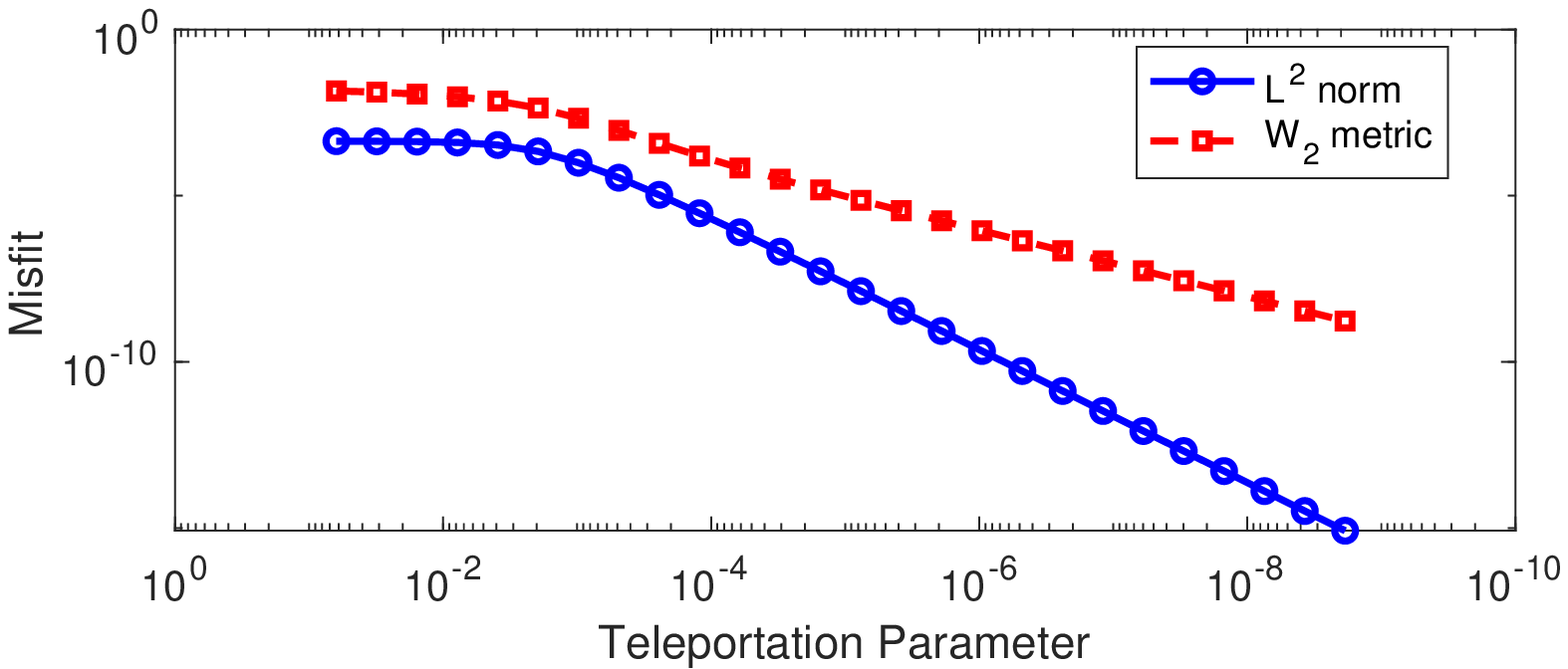} \label{fig:Tele_FPEvsFPE}}
\subfloat[PDE steady state vs.~DNS]{
\includegraphics[width = 0.48\textwidth]{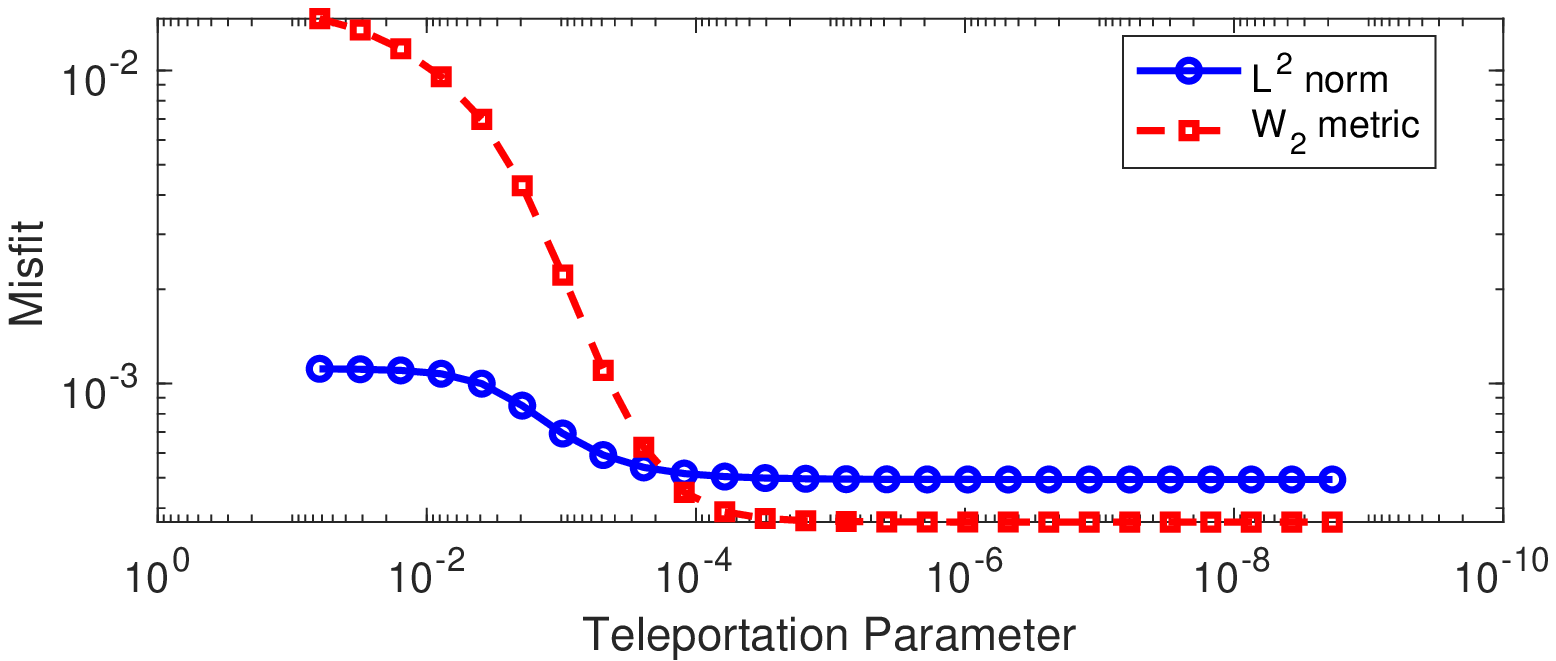} \label{fig:Tele_FPEvsDNS}}
\caption{The $L^2$ norm and $W_2$ metric when the steady-state solution of various teleportation parameters is compared with the steady state without teleportation (left), and  with a fixed invariant measure obtained from the trajectory samples (right).\label{fig:Lorenz_teleportation_compare}}
\end{figure}

In~\Cref{fig:Tele_FPEvsFPE}, we use the steady-state density in which the teleportation parameter $\epsilon=0$ as the reference data. We then compare it with those generated with a nonzero $\epsilon$ in terms of the $L^2$ norm and $W_2$ metric. The misfit monotonically decreases to zero as $\epsilon\rightarrow 0$. When the reference density is replaced by the histogram accumulated from trajectory samples, the misfit again plateaued when $\epsilon$ becomes small since the modeling error, mainly the numerical diffusion from the finite volume solver, becomes the dominant factor of their difference. As discussed in~\Cref{sec:invariant_measure_dx}, the error from numerical diffusion could be effectively reduced as the mesh is refined, i.e., $\Delta x\rightarrow 0$.

\subsection{Parameter Inference} \label{sec:single_parameter}
One main goal of this work is to perform parameter identification using the invariant measure, a macroscopic statistical quantity, as the data, rather than inferring the parameter directly through the time trajectories. All steady-state distributions in this section are solved on a mesh with spacing $\Delta x = 3$.

\subsubsection{Single Parameter Inference}
\begin{figure}
    \centering
   \subfloat[Single-parameter inversion]{    \includegraphics[width = 0.3\textwidth]{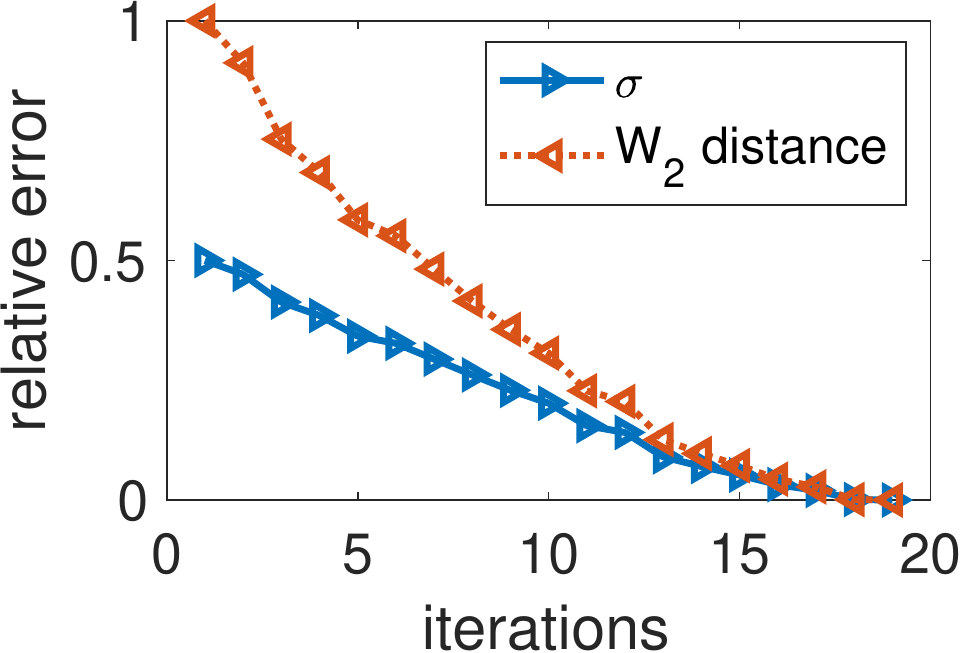}
    \includegraphics[width = 0.3\textwidth]{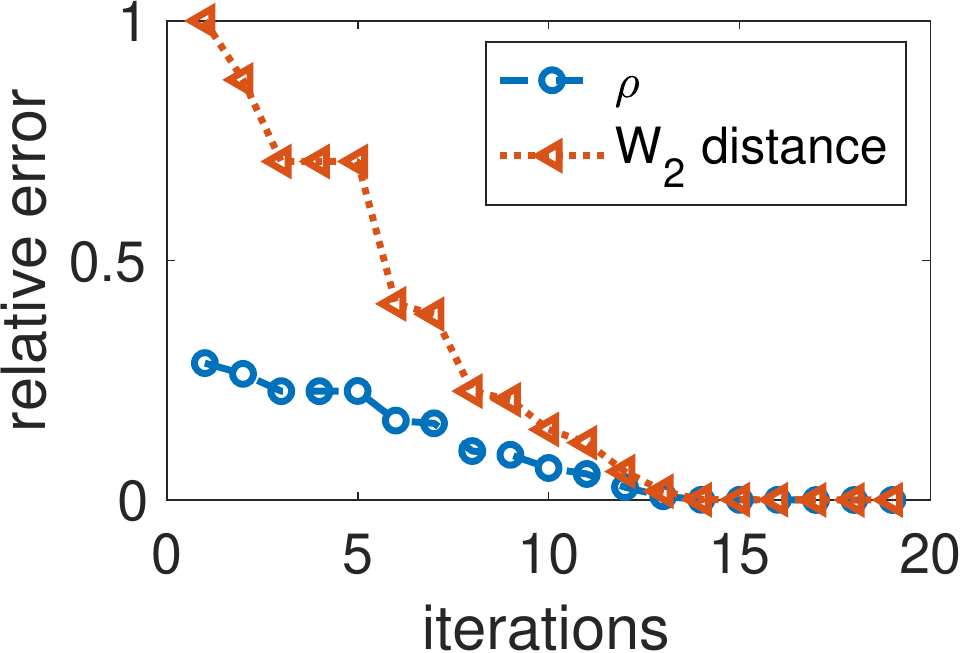}
    \includegraphics[width = 0.3\textwidth]{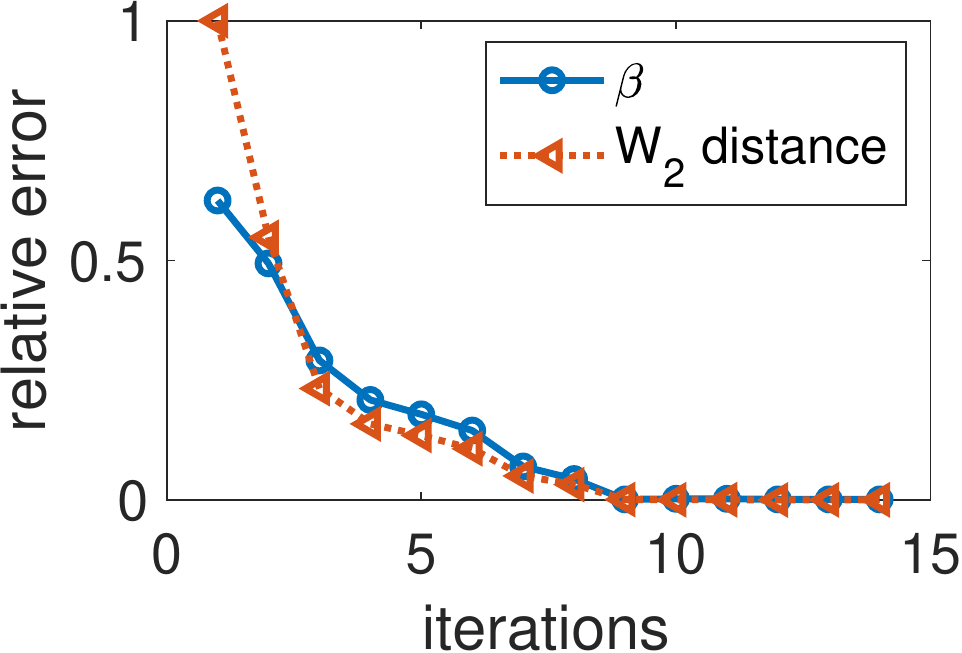}\label{fig:lorenz_W2crime_inv_single}}\\
   \subfloat[Multi-parameter inversion]{\includegraphics[height = 2.5cm]{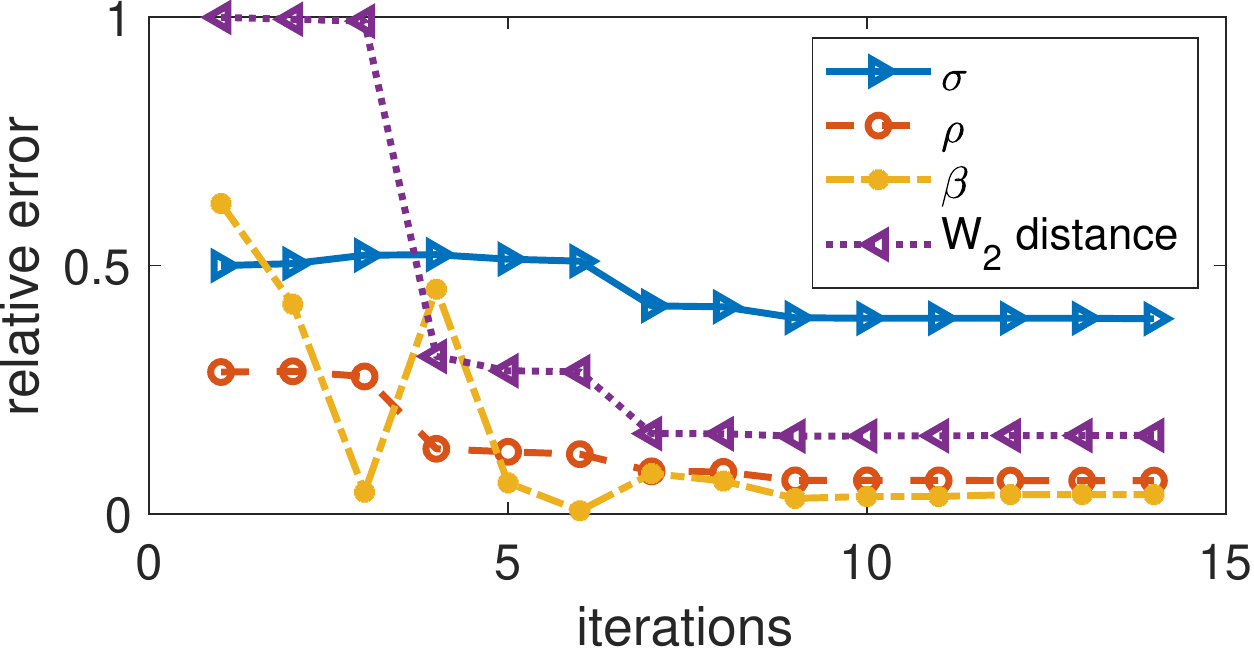}     \includegraphics[height = 2.5cm]{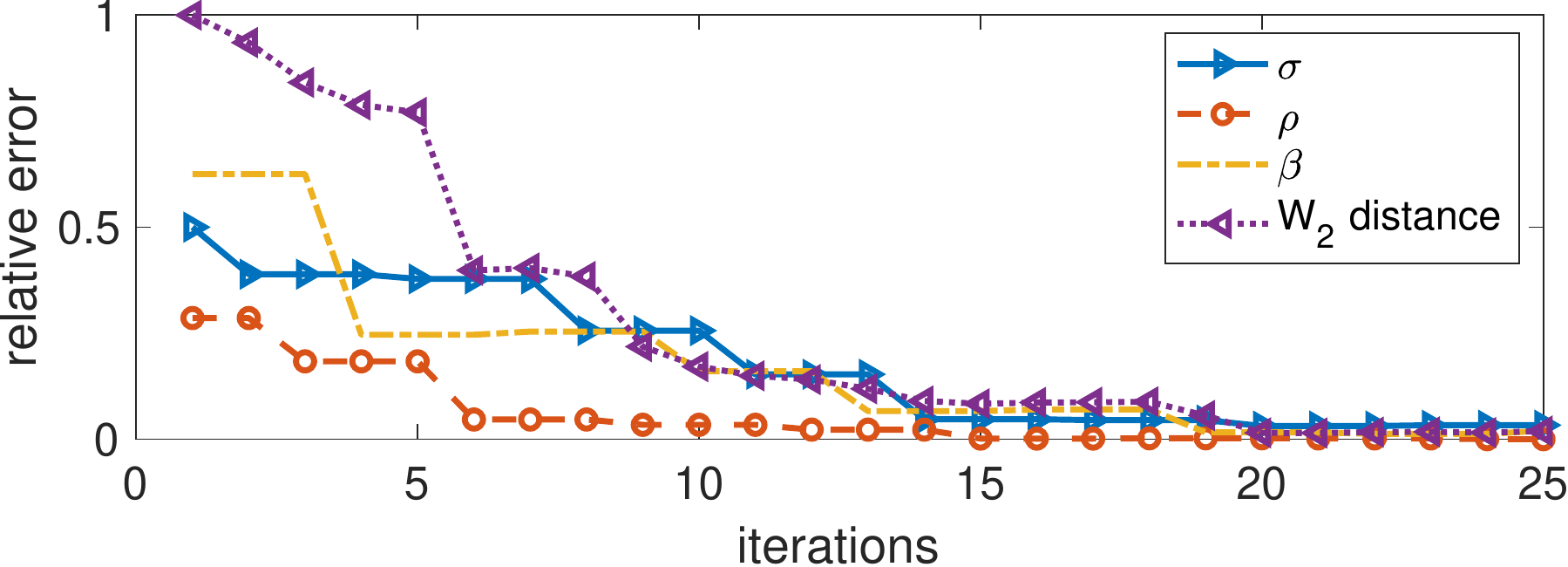}\label{fig:lorenz_W2crime_inv_mul}}
    \caption{Top row: Lorenz system single-parameter inference starting with $\sigma = 5$ (left), $\rho=20$ (middle), $\beta=1$ (right), respectively. Bottom row: multi-parameter inference by updating three parameters simultaneously (bottom left) and using coordinate gradient descent (bottom right) with initial guess $(\sigma, \rho, \beta) = (5, 20, 1)$. The reference PDF is generated through the same numerical solver producing the synthetic PDF.}
    \label{fig:Lorenz_inverse_crime_W2}
\end{figure}

We first focus on the single-parameter reconstruction by assuming that the other parameters in the dynamical systems are accurately known. \Cref{fig:lorenz_W2crime_inv_single} shows the single-parameter inversions of the Lorenz system where the ones for the R\"ossler and Chen systems can be found in~\Cref{sec:RosslerChen}. All experiments use the squared $W_2$ metric as the objective function; see~\eqref{eq:f}. One can see that both the objective function that measures the data mismatch and the relative error of the reconstructed parameters decay to zero rapidly.

We remark that in these tests, the target invariant measure (our reference data) is simulated as the steady-state solution to~\eqref{eq:continuity} at the true parameters, using the same PDE solver that produces the synthetic data. Later, to mimic the realistic scenarios, we will show numerical inversion tests where the reference data directly comes from time trajectories and thus contains both noise and model discrepancy.

\subsubsection{Multi-Parameter Inference via Coordinate Gradient Descent}

For numerical tests we consider here, all dynamical systems have three parameters, while our observation is the invariant measure $\rho(\theta_1,\theta_2,\theta_3)$. Under certain assumptions for the continuous dependency on the parameters, the first-order variation gives 
\[
    \delta\rho = \rho_{\theta_1}\delta\theta_1 + \rho_{\theta_2}\delta\theta_2  + \rho_{\theta_3}\delta\theta_3,
\]
which discloses the issue of multi-parameter inversion. In the forward problem, a small perturbation in each parameter causes a corresponding perturbation in the data $\rho$, but in the inverse problem, the observed misfit in $\rho$ could be contributed from any of the parameters, causing nonzero and possibly wrong gradient updates.

Numerical strategies exist to reduce the inter-parameter trade-off. One may mitigate the inter-parameter dependency either from the formulation of the optimization problem or through the optimization algorithm. Here, we take the second pathway: separate the parameters in the optimization algorithm by using the coordinate gradient descent by only updating one parameter at one iteration. 

\Cref{fig:lorenz_W2crime_inv_mul} shows the Lorenz system multi-parameter inversion. We remark again that the reference data in these tests are produced by the same PDE solver that produces the synthetic data and thus contains no modeling discrepancy. The left plot in \Cref{fig:lorenz_W2crime_inv_mul} shows the convergence history of simultaneously updating all three parameters, but the iterates get stuck at an incorrect set of values with no feasible descent direction. On the other hand, the right plot shows the convergence result using coordinate gradient descent. The gradient descent algorithm quickly converges to the true value $(\sigma, \rho, \beta) = (10,28,8/3)$ starting from $(5, 20, 1)$. The different convergence behaviors of the two plots in~\Cref{fig:lorenz_W2crime_inv_mul} demonstrate that the reconstruction process is affected by the inter-parameter interaction.

\subsubsection{Parameter Inference for Chaotic Systems with Noise}\label{sec:real_inversion}

\begin{figure}
    \centering
\subfloat[Lorenz with extrinsic noise]{\includegraphics[width = 0.49\textwidth]{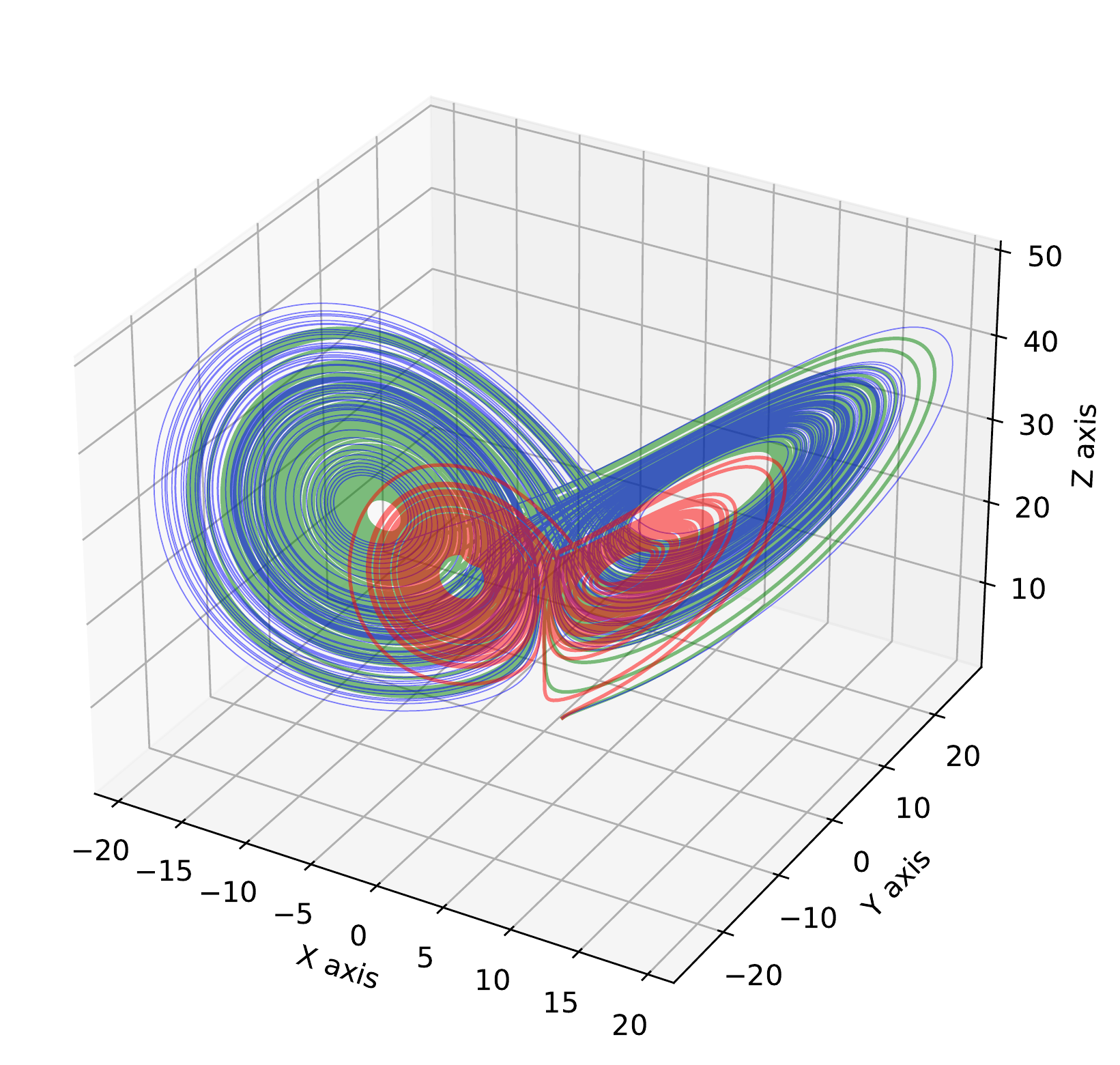}}
\subfloat[Arctan Lorenz with intrinsic noise]{\includegraphics[width = 0.49\textwidth]{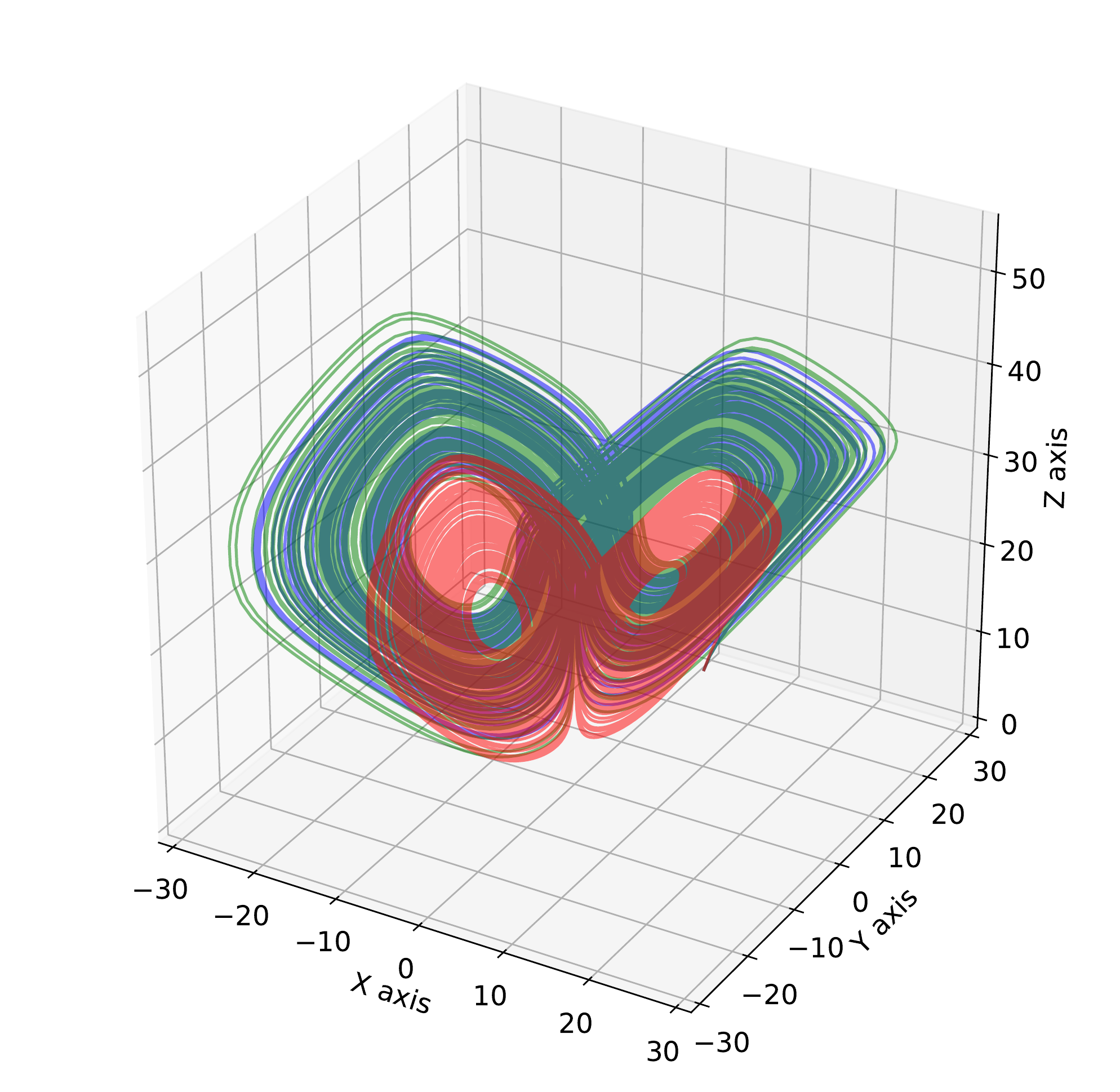}}
    \caption{Comparison among the dynamics produced by the initial parameter (red); true parameter (green); reconstructed parameters (blue) for two examples.}
    \label{fig:Lorenz_dym_compare}
\end{figure}

\begin{figure}
    \centering
\includegraphics[width = 0.9\textwidth]{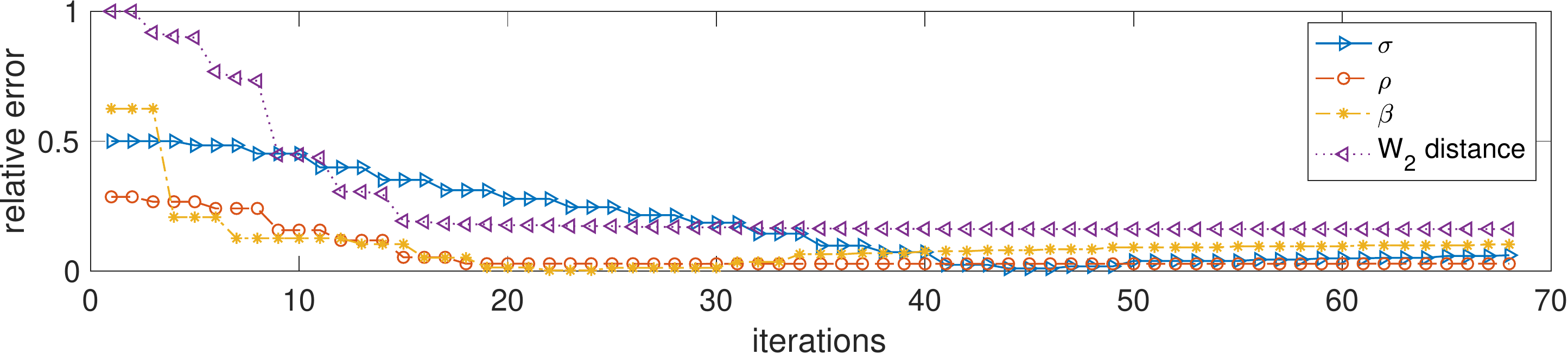}
\caption{
Lorenz system~\eqref{eq:Lorenz}: Multi-parameter inference using coordinate gradient descent with initial guess $(\sigma, \rho, \beta) = (5, 20, 1)$. The reference PDF is the histogram from the time trajectory with \textit{extrinsic} noise.}
\label{fig:Lorenz_extrinsic_noise_W2}
\end{figure}

In this work, we formulate an inverse problem into a nonlinear regression problem, usually subject to at least three sources of errors: model discrepancy, data noise, and optimization error. As discussed earlier, the almost perfect reconstructions in the previous section are achieved under the so-called ``inverse crime'' regime and thus are immune to the first two types of errors. Here, we set up tests to avoid the ``inverse crime'' regime. We first solve the dynamical system forward in time with a fixed time step $\Delta t$ from $t=0$ to $T = 2\times 10^6$, achieving the DNS solution. We then randomly subsample $10^4$ state-space positions; see~\Cref{fig:Lorenz_subsample} for their illustrations. The reference data, i.e., the target estimated invariant measure, is obtained from the histogram that results from binning the subsampled data into cubic boxes in $\mathbb{R}^3$. Moreover, we also use time trajectories affected by intrinsic and extrinsic noises. Starting from the initial guess $(5,20,1)$, the multi-parameter inversion for the Lorenz system~\eqref{eq:Lorenz} with the extrinsic noise converges to $(\sigma, \rho, \beta) = (10.63, 28.82, 3.04)$, and the test with the intrinsic noise converges to $(10.50, 28.41, 2.89)$. For the Arctan Lorenz system~\eqref{eq:arctan_lorenz}, the reconstruction converges to $(11.37, 27.64, 2.35)$ starting from $(5,20,2)$, where the reference data is polluted by the intrinsic noise. We demonstrate the reconstructed dynamics in~\Cref{fig:Lorenz_dym_compare}. Plot for the convergence history of the Lorenz example is shown in~\Cref{fig:Lorenz_extrinsic_noise_W2}. More numerical results can be found in~\Cref{sec:inference_with_noise}.


Earlier in~\Cref{sec:invariant_measure_dx}, we have analyzed the numerical error between the synthetic steady-state solution using the first-order finite volume method. It is shown both in~\Cref{fig:model-discrepancy} and by numerical analysis that the error grows linearly with $\Delta x$. It is also a good characterization of the model discrepancy and could be utilized to design specific stopping criteria to avoid parameter overfitting. For example, \Cref{fig:model-discrepancy} could serve as the baseline: whenever the objective function ($W_2$ metric in our case) is minimized to a value smaller than the model discrepancy, one should execute early stopping: terminate the iterative parameter reconstruction to avoid overfitting the noise. In machine learning, early stopping is designed to monitor the generalization error of one model and stop training when generalization error begins to degrade, which is quite similar to the situation we encounter here.

\subsection{Discussions and Future Directions}\label{subsec:future}
We presented several preliminary numerical results to illustrate the feasibility of our proposed method. There are quite a few future directions that we wish to pursue to improve efficiency and accuracy. 

\subsubsection{The Choice of the Objective Function}
We used the squared $W_2$ metric as the objective function to measure the discrepancy between the reference data and the synthetic stationary distribution. We were motivated by the advantageous properties of the $W_2$ metric, such as the differentiability (see~\Cref{sec:f_diff}), the geometric feature, and the robustness to noises and small perturbations~\cite{engquist2020quadratic,dunlop2020stability}. Nevertheless, it would be interesting to investigate further other choices of objective functions, including the commonly used $1$-Wasserstein metric as well as many other families of probability metrics~\cite{gibbs2002choosing}. Intuitively, we know that the TV and Hellinger distances will not reflect the geometric differences between two delta functions as $d_{\text{TV}}(\delta_x,\delta_y) = d_{\text{Hell}}(\delta_x,\delta_y) = 1$, $\forall x,y\in\R^3$, which potentially causes local minima trapping. Moreover, the $\chi^2$ and Kullback--Leibler divergences are not suitable to compare measures with compact and singular supports. Note that the supports of the synthetic and reference measures in our application may not overlap and are on the low-dimensional manifolds, which inevitably causes $\rho=0$ in the denominator of such divergences. A complete study is needed along this direction.

\subsubsection{The Computational Cost and the Modeling Discrepancy}
Similar to all computational inverse problems solved as PDE-constrained optimization problems, the major bottleneck in memory and computational cost is solving the forward problem repetitively throughout the gradient- or Hessian-based optimization algorithms. It takes from a few hours to a few days on a single computer to produce the stationary distribution once on a fine grid as in~\Cref{fig:Lorenz_FPE_pdf},~\Cref{fig:Rossler_FPE_pdf} and~\Cref{fig:Chen_FPE_pdf}, which is unrealistic for solving inverse problems. In the tests for parameter estimation, we use a much more coarse grid for the PDE solver, which gives us a much smaller Markov matrix so that it is feasible to compute the stationary distribution repetitively. The corresponding histograms of the time trajectory are also accumulated on the coarsened grid. It is essential to understand the error incurred in parameter estimation by producing the synthetic data on a coarsened grid. One may expect some balance between the computation time of solving the inverse problem and the error contributed by the numerical solver.

Besides a proper choice of the grid size, we only use a first-order finite volume discretization for~\eqref{eq:continuity} in this paper. It directly affects the resulting Markov matrix and the computed steady-state solution to~\eqref{eq:continuity}. It also directly contributes to the model discrepancy, a major source of error in the parameter inversion, as discussed earlier. A more accurate discretization of~\eqref{eq:continuity}, such as those including the corner transport and second-order terms~\cite{bewley2012efficient,leveque2002finite}, could reduce the model discrepancy and mitigate the overfitting phenomenon. Beyond higher-order discretizations, exploring different approaches such as adaptive cell-based approximations or ``SRB"-based methods (see~\Cref{sec:FPE_FWD} for more literature review) to approximate the invariant measure and  exploit the inherent sparsity of the problem will be particularly advantageous or even necessary for higher-dimensional state spaces. 


However, these numerical issues highlight a particular challenge associated with the PDE framework with respect to the applicability of this approach to even moderately high-dimensional state spaces.  This has been a perennial challenge for solving problems involving high dimensional flows in state space as encountered in the solutions to Vlasov, Fokker--Planck, and Boltzmann equations.  While solutions that exploit natural sparsity~\cite{bewley2012efficient} offer the potential to reduce the complexity from that of the dimension of the ambient state-space down to that of the attractor, adaptive models such as moving meshes and arbitrary Lagrangian-Eulerian (ALE)  methods~\cite{masud2005,TAITANO2021107861} and hierarchically adaptive methods~\cite{arslanbekov2013} will likely be necessary to extend this approach beyond very low dimensional state-space. However, extending this approach to very high-dimensional problems will likely require significantly larger modifications using recently developed approaches to high-dimensional problems such as rank adaptive tensor~\cite{dektor2021} or machine learning methods~\cite{Line2024713118}. 

\subsubsection{The Data Requirements}
The approach, as explored here, assumes access to a parameterized and differentiable representation of the state-space dynamics valid near and especially on the invariant measure of the system. Access to fully resolved state-space trajectory data in the model state-space coordinates, transformed here into a reference measure $\rho^*$, could be challenging for realistic problems. When both these conditions are met, existing methods such as SINDy~\cite{schmid2010dynamic,brunton2016discovering,schaeffer2018extracting} and related sparse regression frameworks are likely to be considerably more data-efficient than the proposed approach. However, an advantage to this approach is that it does not require that the reference data be sampled at timescales commensurate with the inherent dynamics of the autonomous system. Assuming that the observed states are simply samples from some invariant measure, data sampled slowly with respect to timescales for which prediction would be well-posed due to the chaotic divergence of individual trajectories is still well suited for use as reference data. 

The longer-term goal is to connect this framework to the fully black-box parameter estimation problem of matching the invariant measure in the time-delay-embedded coordinates that can be constructed from available observable data 
as investigated in~\cite{greve2019data}.  This will instead require differentiable surrogate models of the flow in time-delay observable coordinates.  While this poses a different set of challenges left to future investigations, the purpose of this work has been to explore mathematical foundations for the empirically motivated parameter estimation problem described in the prior work~\cite{greve2019data} related to the conditions under which the steady-state distribution can be expected to exist and the differentiability of the Wasserstein metric in the parameter space (see~\Cref{sec:f_diff}). Exploration of the extension of these topics to black-box parameter estimation based on the system's flow as observed in delay-embedded coordinates is left to future work.

\section{Conclusion}\label{conclusion}
In this paper, we propose a data-driven approach for parameter estimation of chaotic dynamical systems. There are two significant contributions. First, we shift from an ODE forward model to the related PDE forward model through the tool of physical measure. Instead of using pure time trajectories as the inference data, we treat statistics accumulated from the direct numerical simulation as the observable, whose continuous analog is the steady-state solution to~\eqref{eq:continuity}. As a result, the original parameter identification problem is translated into a data-fitting, PDE-constrained optimization problem. We then use an upwind scheme based on the finite volume method to discretize and solve the forward problem. Second, we use the quadratic Wasserstein metric from optimal transportation as the data fidelity term measuring the difference between the synthetic and the reference datasets. We first provide a rigorous analysis of the differentiability regarding the Wasserstein-based parameter estimation and then derive two ways of calculating the Wasserstein gradient following the discretize-then-optimize approach. In particular, the adjoint approach is efficient as the computational cost of gradient evaluation is independent of the size of the unknown parameters, making the method scalable for large-scale parameterization of the velocity fields. Finally, we show several numerical results to demonstrate the promises of this new approach for chaotic dynamical system parameter identification.

For this method, sufficient data is required to converge the histogram estimate of the reference distribution. As in any non-parametric density estimate, the amount of data is therefore dependent on the coarseness of the approximation and level of stochastic error tolerated.  In this work, knowledge of the full state is also presumed. The approximated invariant measure from the time trajectories as our reference data might be a singular probability measure with highly complex support that has fractional fractal dimension. Thus, we use the regularized forward PDE model as a surrogate in solving this inverse problem. We approximate the steady-state solution to the PDE model with first-order accuracy based on the finite-volume upwind discretization. Due to the sparsity of the Markov matrix and a coarse grid, we can evaluate the gradient of the resulting PDE-constrained optimization problem quite efficiently in terms of both memory and computation complexity. The Wasserstein metric from optimal transportation is our objective function, which can compare measures with singular and compact support and handle the fractional fractal dimension of the reference invariant measure. Future works along the lines discussed in~\Cref{subsec:future} would generalize and help improve the proposed method.

\section*{Acknowledgements}
YY gratefully acknowledges the support by National Science Foundation through grant number DMS-1913129. LN was partially supported by AFOSR MURI FA 9550 18-1-0502 grant. EN acknowledges that results in this paper were obtained in part using a high-performance computing system acquired through NSF MRI grant DMS-1337943 to WPI. RM was partially supported by AFOSR Grants FA9550-20RQCOR098 (PO: Leve) and FA9550-20RQCOR100 (PO: Fahroo).
We thank Prof.\ Adam Oberman for initiating our collaboration. We also thank Prof.\ Alex Townsend for his constructive suggestions.

This work was in part completed during the long program on High Dimensional Hamilton--Jacobi PDEs held in the Institute for Pure and Applied Mathematics (IPAM) at UCLA, March 9-June 12, 2020. The authors thank the program's organizers, IPAM scientific committee, and staff for the hospitality and stimulating research environment.

The authors are also grateful to the peer referees for their time, comments, and constructive suggestions during the review process.

\bibliographystyle{siamplain}
\bibliography{lorenz.bib}

\begin{thebibliography}{10}

\bibitem{aguirre2009modeling}
{\sc L.~A. Aguirre and C.~Letellier}, {\em Modeling nonlinear dynamics and
  chaos: a review}, Mathematical Problems in Engineering, 2009 (2009).

\bibitem{alberti99}
{\sc G.~Alberti and L.~Ambrosio}, {\em A geometrical approach to monotone
  functions in $\mathbb{R}^n$}, Mathematische Zeitschrift, 230 (1999),
  pp.~259--316.

\bibitem{allawala2016statistics}
{\sc A.~Allawala and J.~Marston}, {\em {Statistics of the stochastically forced
  Lorenz attractor by the Fokker--Planck equation and cumulant expansions}},
  Physical Review E, 94 (2016), p.~052218.

\bibitem{ags08}
{\sc L.~Ambrosio, N.~Gigli, and G.~Savar\'{e}}, {\em Gradient flows in metric
  spaces and in the space of probability measures}, Lectures in Mathematics ETH
  Z\"{u}rich, Birkh\"{a}user Verlag, Basel, second~ed., 2008.

\bibitem{arjovsky2017wasserstein}
{\sc M.~Arjovsky, S.~Chintala, and L.~Bottou}, {\em {W}asserstein generative
  adversarial networks}, in International conference on machine learning, PMLR,
  2017, pp.~214--223.

\bibitem{arslanbekov2013}
{\sc R.~R. Arslanbekov, V.~I. Kolobov, and A.~A. Frolova}, {\em Kinetic solvers
  with adaptive mesh in phase space}, Phys. Rev. E, 88 (2013), p.~063301.

\bibitem{baake1992fitting}
{\sc E.~Baake, M.~Baake, H.~Bock, and K.~Briggs}, {\em Fitting ordinary
  differential equations to chaotic data}, Physical Review A, 45 (1992),
  p.~5524.

\bibitem{bakker2000learning}
{\sc R.~Bakker, J.~C. Schouten, C.~L. Giles, F.~Takens, and C.~M. Van
  Den~Bleek}, {\em Learning chaotic attractors by neural networks}, Neural
  Computation, 12 (2000), pp.~2355--2383.

\bibitem{bernton2019parameter}
{\sc E.~Bernton, P.~E. Jacob, M.~Gerber, and C.~P. Robert}, {\em On parameter
  estimation with the {W}asserstein distance}, Information and Inference: A
  Journal of the IMA, 8 (2019), pp.~657--676.

\bibitem{bewley2012efficient}
{\sc T.~R. Bewley and A.~S. Sharma}, {\em Efficient grid-based {B}ayesian
  estimation of nonlinear low-dimensional systems with sparse non-{G}aussian
  {PDF}s}, Automatica, 48 (2012), pp.~1286--1290.

\bibitem{bezruchko2010extracting}
{\sc B.~P. Bezruchko and D.~A. Smirnov}, {\em Extracting knowledge from time
  series: An introduction to nonlinear empirical modeling}, Springer Science \&
  Business Media, 2010.

\bibitem{bonnans00}
{\sc J.~F. Bonnans and A.~Shapiro}, {\em Perturbation analysis of optimization
  problems}, Springer Series in Operations Research, Springer-Verlag, New York,
  2000.

\bibitem{bowen1975equilibrium}
{\sc R.~Bowen}, {\em Equilibrium state and the ergodic theory of anosov
  diffeomorphisms lecture notes in mathematics no. 470 springer-verlag berlin},
  Bricmont and A. Kupiainen 1994 Coupled analytic maps, to appear in
  Non-linearity,  (1975).

\bibitem{jax2018github}
{\sc J.~Bradbury, R.~Frostig, P.~Hawkins, M.~J. Johnson, C.~Leary,
  D.~Maclaurin, G.~Necula, A.~Paszke, J.~Vander{P}las, S.~Wanderman-{M}ilne,
  and Q.~Zhang}, {\em {JAX}: composable transformations of {P}ython+{N}um{P}y
  programs}, 2018, \url{http://github.com/google/jax}.

\bibitem{brunton2016discovering}
{\sc S.~L. Brunton, J.~L. Proctor, and J.~N. Kutz}, {\em Discovering governing
  equations from data by sparse identification of nonlinear dynamical systems},
  Proceedings of the national academy of sciences, 113 (2016), pp.~3932--3937.

\bibitem{chen1999yet}
{\sc G.~Chen and T.~Ueta}, {\em Yet another chaotic attractor}, International
  Journal of Bifurcation and chaos, 9 (1999), pp.~1465--1466.

\bibitem{cowieson_young_2005}
{\sc W.~Cowieson and L.-S. Young}, {\em {SRB measures as zero-noise limits}},
  Ergodic Theory and Dynamical Systems, 25 (2005), p.~1115–1138,
  \url{https://doi.org/10.1017/S0143385704000604}.

\bibitem{dektor2021}
{\sc A.~Dektor, A.~Rodgers, and D.~Venturi}, {\em Rank-adaptive tensor methods
  for high-dimensional nonlinear pdes}, Journal of Scientific Computing, 88
  (2021).

\bibitem{dellnitz2001algorithms}
{\sc M.~Dellnitz, G.~Froyland, and O.~Junge}, {\em The algorithms behind
  {GAIO}—set oriented numerical methods for dynamical systems}, in Ergodic
  theory, analysis, and efficient simulation of dynamical systems, Springer,
  2001, pp.~145--174.

\bibitem{dellnitz1997almost}
{\sc M.~Dellnitz and O.~Junge}, {\em Almost invariant sets in chua's circuit},
  International Journal of Bifurcation and Chaos, 7 (1997), pp.~2475--2485.

\bibitem{dellnitz1998adaptive}
{\sc M.~Dellnitz and O.~Junge}, {\em An adaptive subdivision technique for the
  approximation of attractors and invariant measures}, Computing and
  Visualization in Science, 1 (1998), pp.~63--68.

\bibitem{dellnitz1999approximation}
{\sc M.~Dellnitz and O.~Junge}, {\em On the approximation of complicated
  dynamical behavior}, SIAM Journal on Numerical Analysis, 36 (1999),
  pp.~491--515.

\bibitem{dellnitz2002set}
{\sc M.~Dellnitz and O.~Junge}, {\em Chapter 5 - set oriented numerical methods
  for dynamical systems}, in Handbook of Dynamical Systems, B.~Fiedler, ed.,
  vol.~2 of Handbook of Dynamical Systems, Elsevier Science, 2002,
  pp.~221--264.

\bibitem{dunlop2020stability}
{\sc M.~M. Dunlop and Y.~Yang}, {\em Stability of gibbs posteriors from the
  wasserstein loss for bayesian full waveform inversion}, SIAM/ASA Journal on
  Uncertainty Quantification, 9 (2021), pp.~1499--1526.

\bibitem{effah2018study}
{\sc S.~Effah-Poku, W.~Obeng-Denteh, and I.~Dontwi}, {\em A study of chaos in
  dynamical systems}, Journal of Mathematics, 2018 (2018).

\bibitem{eidenschink1997exploring}
{\sc M.~Eidenschink}, {\em Exploring global dynamics: A numerical algorithm
  based on the Conley Index theory}, PhD thesis, Georgia Institute of
  Technology, 1995.

\bibitem{engquist2020quadratic}
{\sc B.~Engquist, K.~Ren, and Y.~Yang}, {\em The quadratic {W}asserstein metric
  for inverse data matching}, Inverse Problems, 36 (2020), p.~055001.

\bibitem{engquist2020optimal}
{\sc B.~Engquist and Y.~Yang}, {\em Optimal transport based seismic inversion:
  Beyond cycle skipping}, Communications on Pure and Applied Mathematics,
  (2020).

\bibitem{evansgariepy92}
{\sc L.~C. Evans and R.~F. Gariepy}, {\em Measure theory and fine properties of
  functions}, Studies in Advanced Mathematics, CRC Press, Boca Raton, FL, 1992.

\bibitem{feng2008reconstruction}
{\sc J.~C. Feng}, {\em Reconstruction of chaotic signals with applications to
  chaos-based communications}, World Scientific, 2008.

\bibitem{fernandez2014generalized}
{\sc D.~C. D.~R. Fern{\'a}ndez, P.~D. Boom, and D.~W. Zingg}, {\em A
  generalized framework for nodal first derivative summation-by-parts
  operators}, Journal of Computational Physics, 266 (2014), pp.~214--239.

\bibitem{fiedler2002handbook}
{\sc B.~Fiedler}, {\em Handbook of dynamical systems}, Gulf Professional
  Publishing, 2002.

\bibitem{flamary2021pot}
{\sc R.~Flamary, N.~Courty, A.~Gramfort, M.~Z. Alaya, A.~Boisbunon, S.~Chambon,
  L.~Chapel, A.~Corenflos, K.~Fatras, N.~Fournier, et~al.}, {\em {POT: Python
  Optimal Transport}}, Journal of Machine Learning Research, 22 (2021),
  pp.~1--8.

\bibitem{fradkov2005control}
{\sc A.~L. Fradkov and R.~J. Evans}, {\em Control of chaos: Methods and
  applications in engineering}, Annual Reviews in Control, 29 (2005),
  pp.~33--56.

\bibitem{froyland2001extracting}
{\sc G.~Froyland}, {\em Extracting dynamical behavior via markov models}, in
  Nonlinear dynamics and statistics, Springer, 2001, pp.~281--321.

\bibitem{gabor2015robust}
{\sc A.~G{\'a}bor and J.~R. Banga}, {\em Robust and efficient parameter
  estimation in dynamic models of biological systems}, BMC systems biology, 9
  (2015), pp.~1--25.

\bibitem{gibbs2002choosing}
{\sc A.~L. Gibbs and F.~E. Su}, {\em On choosing and bounding probability
  metrics}, International statistical review, 70 (2002), pp.~419--435.

\bibitem{givon2004extracting}
{\sc D.~Givon, R.~Kupferman, and A.~Stuart}, {\em Extracting macroscopic
  dynamics: model problems and algorithms}, Nonlinearity, 17 (2004), p.~R55.

\bibitem{gleich2015pagerank}
{\sc D.~F. Gleich}, {\em Pagerank beyond the web}, SIAM Review, 57 (2015),
  pp.~321--363.

\bibitem{GVL12}
{\sc G.~H. Golub and C.~F. Van~Loan}, {\em Matrix {C}omputations, 4th ed.},
  Johns Hopkins University Press, Baltimore, 2013.

\bibitem{greve2019data}
{\sc C.~Greve, K.~Hara, R.~Martin, D.~Eckhardt, and J.~Koo}, {\em A data-driven
  approach to model calibration for nonlinear dynamical systems}, Journal of
  Applied Physics, 125 (2019), p.~244901.

\bibitem{griewank2008evaluating}
{\sc A.~Griewank and A.~Walther}, {\em Evaluating derivatives: principles and
  techniques of algorithmic differentiation}, SIAM, 2008.

\bibitem{haker2004optimal}
{\sc S.~Haker, L.~Zhu, A.~Tannenbaum, and S.~Angenent}, {\em Optimal mass
  transport for registration and warping}, International Journal of computer
  vision, 60 (2004), pp.~225--240.

\bibitem{huang18}
{\sc W.~Huang, M.~Ji, Z.~Liu, and Y.~Yi}, {\em Concentration and limit
  behaviors of stationary measures}, Physica D: Nonlinear Phenomena, 369
  (2018), pp.~1--17.

\bibitem{jacobs20}
{\sc M.~Jacobs and F.~L{\'e}ger}, {\em A fast approach to optimal transport:
  The back-and-forth method}, Numerische Mathematik, 146 (2020), pp.~513--544.

\bibitem{jaeger1996unbiased}
{\sc L.~Jaeger and H.~Kantz}, {\em Unbiased reconstruction of the dynamics
  underlying a noisy chaotic time series}, Chaos: An Interdisciplinary Journal
  of Nonlinear Science, 6 (1996), pp.~440--450.

\bibitem{KaiserJFM2014}
{\sc E.~Kaiser, B.~R. Noack, L.~Cordier, A.~Spohn, M.~Segond, M.~Abel,
  G.~Daviller, J.~{\"O}sth, S.~Krajnovi{\'c}, and R.~K. Niven}, {\em
  Cluster-based reduced-order modelling of a mixing layer}, Journal of Fluid
  Mechanics, 754 (2014), pp.~365--414.

\bibitem{kifer1986small}
{\sc Y.~Kifer}, {\em General random perturbations of hyperbolic and expanding
  transformations}, Journal d'Analyse Math{\'e}matique, 47 (1986),
  pp.~111--150.

\bibitem{kostelich1992problems}
{\sc E.~J. Kostelich}, {\em Problems in estimating dynamics from data}, Physica
  D: Nonlinear Phenomena, 58 (1992), pp.~138--152.

\bibitem{leveque2002finite}
{\sc R.~J. LeVeque}, {\em Finite volume methods for hyperbolic problems},
  vol.~31, Cambridge university press, 2002.

\bibitem{Line2024713118}
{\sc A.~T. Lin, S.~W. Fung, W.~Li, L.~Nurbekyan, and S.~J. Osher}, {\em
  Alternating the population and control neural networks to solve
  high-dimensional stochastic mean-field games}, Proceedings of the National
  Academy of Sciences, 118 (2021).

\bibitem{lu2021deepxde}
{\sc L.~Lu, X.~Meng, Z.~Mao, and G.~E. Karniadakis}, {\em Deepxde: A deep
  learning library for solving differential equations}, SIAM Review, 63 (2021),
  pp.~208--228.

\bibitem{masud2005}
{\sc A.~Masud and L.~Bergman}, {\em Solution of the four dimensional
  fokker--planck equation: still a challenge.}, in ICOSSAR 2005, Millpress,
  Rotterdam, 2005.

\bibitem{mcgoff2015statistical}
{\sc K.~McGoff, S.~Mukherjee, and N.~Pillai}, {\em Statistical inference for
  dynamical systems: A review}, Statistics Surveys, 9 (2015), pp.~209--252.

\bibitem{medio2001}
{\sc A.~Medio and M.~Lines}, {\em Nonlinear dynamics: A primer}, Cambridge
  University Press, 2001.

\bibitem{meyer2000matrix}
{\sc C.~D. Meyer}, {\em Matrix analysis and applied linear algebra}, vol.~71,
  {}SIAM, 2000.

\bibitem{michalik2009incremental}
{\sc C.~Michalik, R.~Hannemann, and W.~Marquardt}, {\em Incremental single
  shooting---a robust method for the estimation of parameters in dynamical
  systems}, Computers \& Chemical Engineering, 33 (2009), pp.~1298--1305.

\bibitem{nakagawa1999chaos}
{\sc M.~Nakagawa}, {\em Chaos and fractals in engineering}, World Scientific,
  1999.

\bibitem{negrini2021neural}
{\sc E.~Negrini, G.~Citti, and L.~Capogna}, {\em A neural network ensemble
  approach to system identification}, arXiv preprint arXiv:2110.08382,  (2021).

\bibitem{negrini2021system}
{\sc E.~Negrini, G.~Citti, and L.~Capogna}, {\em System identification through
  {Lipschitz} regularized deep neural networks}, Journal of Computational
  Physics, 444 (2021), p.~110549.

\bibitem{nocedal2006numerical}
{\sc J.~Nocedal and S.~Wright}, {\em Numerical optimization}, Springer Science
  \& Business Media, 2006.

\bibitem{peyre2019computational}
{\sc G.~Peyr{\'e} and M.~Cuturi}, {\em Computational Optimal Transport: With
  Applications to Data Science}, Foundations and trends in machine learning,
  Now, the essence of knowledge., 2019.

\bibitem{rippl16}
{\sc T.~Rippl, A.~Munk, and A.~Sturm}, {\em Limit laws of the empirical
  {W}asserstein distance: Gaussian distributions}, Journal of Multivariate
  Analysis, 151 (2016), pp.~90--109.

\bibitem{robin17}
{\sc Y.~Robin, P.~Yiou, and P.~Naveau}, {\em Detecting changes in forced
  climate attractors with wasserstein distance}, Nonlinear Processes in
  Geophysics, 24 (2017), pp.~393--405.

\bibitem{rodriguez2006novel}
{\sc M.~Rodriguez-Fernandez, J.~A. Egea, and J.~R. Banga}, {\em Novel
  metaheuristic for parameter estimation in nonlinear dynamic biological
  systems}, BMC bioinformatics, 7 (2006), pp.~1--18.

\bibitem{ruan2003chaotic}
{\sc H.~Ruan, T.~Zhai, and E.~E. Yaz}, {\em A chaotic secure communication
  scheme with extended {Kalman} filter based parameter estimation}, in
  Proceedings of 2003 IEEE Conference on Control Applications, 2003. CCA 2003.,
  vol.~1, IEEE, 2003, pp.~404--408.

\bibitem{santambrogio15}
{\sc F.~Santambrogio}, {\em Optimal transport for applied mathematicians},
  vol.~87, Springer, 2015.

\bibitem{schaeffer2018extracting}
{\sc H.~Schaeffer, G.~Tran, and R.~Ward}, {\em Extracting sparse
  high-dimensional dynamics from limited data}, SIAM Journal on Applied
  Mathematics, 78 (2018), pp.~3279--3295.

\bibitem{schmid2010dynamic}
{\sc P.~J. Schmid}, {\em Dynamic mode decomposition of numerical and
  experimental data}, Journal of fluid mechanics, 656 (2010), pp.~5--28.

\bibitem{siegmund2006approximation}
{\sc S.~Siegmund and P.~Taraba}, {\em Approximation of box dimension of
  attractors using the subdivision algorithm}, Dynamical Systems, 21 (2006),
  pp.~1--24.

\bibitem{sommerfeld18}
{\sc M.~Sommerfeld and A.~Munk}, {\em Inference for empirical {W}asserstein
  distances on finite spaces}, Journal of the Royal Statistical Society Series
  B, 80 (2018), pp.~219--238.

\bibitem{TAITANO2021107861}
{\sc W.~Taitano, B.~Keenan, L.~Chacón, S.~Anderson, H.~Hammer, and
  A.~Simakov}, {\em An eulerian vlasov-fokker–planck algorithm for spherical
  implosion simulations of inertial confinement fusion capsules}, Computer
  Physics Communications, 263 (2021), p.~107861.

\bibitem{tucker99}
{\sc W.~Tucker}, {\em {The Lorenz attractor exists}}, Comptes Rendus de
  l'Académie des Sciences - Series I - Mathematics, 328 (1999),
  pp.~1197--1202.

\bibitem{tucker02}
{\sc W.~Tucker}, {\em A rigorous {ODE} solver and {S}male's $14$th problem},
  Foundations of Computational Mathematics, 2 (2002), pp.~53--117,
  \url{https://doi.org/10.1007/s002080010018}.

\bibitem{villani03}
{\sc C.~Villani}, {\em Topics in optimal transportation}, vol.~58, American
  Mathematical Soc., 2003.

\bibitem{young2002srb}
{\sc L.-S. Young}, {\em What are {SRB} measures, and which dynamical systems
  have them?}, Journal of Statistical Physics, 108 (2002), pp.~733--754.

\end{thebibliography}

\begin{appendix}
\section{Proofs From~\Cref{sec:f_diff}}\label{sec:proofs}
\subsection{Proof of~\Cref{prp:dir_der_f}}\label{sec:dir_der_f}
\begin{proof}
We fix $\theta_0 \in \Theta$ and firstly prove that (i) implies (ii). Note that \eqref{eq:invariance} follows immediately from \eqref{eq:dir_der_f}. Furthermore, assume that $(\phi,\psi) \in \Phi_c(\rho(\cdot,\theta_0),\rho^*)$ is an arbitrary pair of Kantorovich potentials. Note that $(\phi,\psi)$ are not necessarily from $\mathcal{S}(\theta_0)$. Since $\int_\Omega  \nabla_\theta \rho(x,\theta_0) dx=0$, we can add an arbitrary constant to $\phi$ and assume that $\sup \phi= \|c\|_\infty $. In that case, we obtain that $(\phi^{cc},\phi^c) \in \mathcal{S}(\theta_0)$, and
\[
    \phi^{cc}(x)=\phi(x),~x\in \operatorname{supp}(\rho(\cdot,\theta_0)), \quad \mbox{and}\quad  \phi^c(y)=\psi(y),~y\in \operatorname{supp}(\rho^*).
\]
Next, we have that $\operatorname{supp}(\nabla_\theta \rho(\cdot,\theta_0) ) \subset \operatorname{supp}(\rho(\cdot,\theta_0))$. Therefore, we have that
\[
    \int_\Omega \phi(x) \nabla_\theta \rho(x,\theta_0)dx=\int_\Omega \phi^{cc}(x) \nabla_\theta \rho(x,\theta_0)dx,
\]
and \eqref{eq:grad_f} follows from~\eqref{eq:dir_der_f} and~\eqref{eq:invariance}.

Next, we prove (i). We apply \cite[Proposition 4.12]{bonnans00} with $U=\Theta$, $X=C(\Omega)\times C(\Omega)$, $\Phi=C=K_c$, and an objective function given by
\[    I(\phi,\psi,\theta)=\int_\Omega \phi(x) \rho(x,\theta) dx+ \int_\Omega \psi(y) \rho^*(y) dy.
\]
For $\theta_1,\theta_2 \in \Theta$ such that $[\theta_1,\theta_2] \subset \Theta$, we have that 
\begin{equation}\label{eq:obj_Lip}
    |I(\phi_2,\psi_2,\theta_2)-I(\phi_1,\psi_1,\theta_1)|\leq  \|\phi_2-\phi_1\|_\infty + \|\psi_2-\psi_1\|_\infty +\|\phi_1\|_\infty \|\eta\|_1|\theta_2-\theta_1|,
\end{equation}
and so $I$ is continuous. Since $K_c$ is compact, the sup-compactness condition holds. Furthermore, \textbf{A2} and the dominated convergence theorem yield the directional differentiability of $I(\phi,\psi,\cdot)$ with
\[
    I'(\phi,\psi,\theta_0,\Delta \theta)= \int_\Omega \phi(x) \nabla_\theta \rho (x,\theta_0) dx \cdot \Delta \theta.
\]
Finally, assume that $t_n \to 0+$, $(\phi_n,\psi_n) \in K_c$, $\Delta \theta \in \R^m$, and $(\phi_n,\psi_n)\to (\phi,\psi) \in K_c$. Then by the dominated convergence theorem we have that
\begin{equation*}
    \begin{split}
        &\lim \limits_{n\to \infty} \frac{I(\phi_n,\psi_n,\theta_0+t_n \Delta \theta)-I(\phi_n,\psi_n,\theta_0 \theta)}{t_n}\\
        =&\lim \limits_{n\to \infty} \int_\Omega \phi_n(x) \frac{\rho(x,\theta_0+t_n \Delta \theta)-\rho(x,\theta_0)}{t_n} dx 
        =I'(\phi,\psi,\theta_0,\Delta \theta).
    \end{split}
\end{equation*}
Thus, all conditions in \cite[Proposition 4.12]{bonnans00} are satisfied and \eqref{eq:dir_der_f} follows.
\end{proof}

\subsection{Proof of~\Cref{thm:generic_diff_f}}\label{sec:generic_diff_f}
\begin{proof}
Assume that \textbf{A1-A3} hold. Then \eqref{eq:obj_Lip} yields that $\theta \mapsto I(\phi,\psi,\theta)$ is locally Lipschitz for all $(\phi,\psi)\in C(\Omega)\times C(\Omega)$. Invoking \eqref{eq:K-duality}, we conclude that $f$ is locally Lipschitz and a.e.\ differentiable by Rademacher's theorem \cite[Section 3.1]{evansgariepy92}. 

Next, assume that \textbf{A4} also holds and denote by $C_0=\|c\|_\infty\|h\|_1$. For arbitrary $(\phi,\phi^c) \in K_c$ we have that
\begin{equation*}
\begin{split}
    I(\phi,\phi^c,\theta)+\frac{C_0|\theta|^2}{2}=&\int_\Omega \phi(x) \left(\rho(x,\theta)+\frac{h(x)|\theta|^2}{2}\right) dx+\int_\Omega \phi^c(y) \rho^*(y) dy\\
    &+\left(\|c\|_{\infty} \|h\|_1-\int_\Omega \phi(x) h(x) dx \right)\frac{|\theta|^2}{2}.
\end{split}
\end{equation*}

Since $0\leq \phi \leq \|c\|_\infty$, and $\theta \mapsto \rho(x,\theta) + \frac{h(x)|\theta|^2}{2}$ is convex for a.e.\ $x$, we obtain that $\theta \mapsto I(\phi,\phi^c,\theta)+\frac{C_0|\theta|^2}{2}$ is convex. Invoking Kantorovich duality again, we obtain that
\[
    f(\theta)+\frac{C_0|\theta|^2}{2}=\sup_{(\phi,\phi^c) \in K_c} I(\phi,\phi^c,\theta)+\frac{C_0|\theta|^2}{2}
\]
is convex. Thus, by a theorem of Anderson and Klee \cite{alberti99} $f$ is differentiable up to a set of Hausdorff dimension $d-1$.
\end{proof}

\subsection{Proof of~\Cref{thm:K_pot_uniq}}\label{sec:K_pot_uniq}
\begin{proof}
Fix an arbitrary pair of Kantorovich potentials $(\phi_1,\psi_1)$, $(\phi_2,\psi_2)$. Note that \eqref{eq:cl_int} guarantees that $\operatorname{int}(\operatorname{supp}(\rho)) \neq \emptyset$, and $\{O_k\},~\{E_k\}$ are well defined.

First, we prove that $\phi_2-\phi_1$ is constant on $\operatorname{cl}(O_k)$ for all $k$. Fix an optimal plan $\pi_0 \in \Gamma_0(\rho,\rho^*)$. For all $x \in \operatorname{supp}(\rho)$ there exists $y \in \Omega$ such that $(x,y) \in \operatorname{supp}(\pi_0)$. Therefore $\phi_i(x)+\psi_i(y)=c(x,y)$, and so $\phi_i(x)=\psi_i^c(x)$ for $x\in \operatorname{supp}(\rho)$. Furthermore, since $c \in C^1(\Omega^2)$ is locally Lipschitz continuous, $\phi_i$ are locally Lipschitz continuous in $O_k$. Thus, by Rademacher's theorem we have that $\phi_i$ are a.e.\ differentiable in $O_k$, and by \cite[Proposition 1.15]{santambrogio15} we obtain that $\nabla \phi_2=\nabla \phi_1$ a.e.\ in $O_k$. Since $O_k$ are connected and $\phi_i$ are continuous, we obtain that $\phi_2-\phi_1=\lambda_k$ in $\operatorname{cl}(O_k)$ for some constants $\lambda_k$.

Next, we show that $\lambda_k=\lambda_l$ for all $k,l$. We start with a claim that
\begin{equation}\label{eq:E_k_O_k}
    \psi_i(y)=\inf_{x\in \operatorname{cl}(O_k)} \{ c(x,y)-\phi_i(x) \},\quad y \in E_k.
\end{equation}
Indeed, we have that $y=\lim_{n\to \infty} y_n$ where $y_n$ are such that $(x_n,y_n) \in \operatorname{supp}(\pi_n)$ for some $\pi_n \in \Gamma_0(\rho,\rho^*)$, and $x_n \in \operatorname{cl}(O_k)$. Therefore, for all $n$ we have that $\phi_i(x_n)+\psi_i(y_n)=c(x_n,y_n)$, and so
\[
    \psi_i(y_n)=\inf_{x\in \operatorname{cl}(O_k)} \{ c(x,y_n)-\phi_i(x) \}.
\]
Since both $\psi_i$ and $y \mapsto \inf_{x\in \operatorname{cl}(O_k)} \{ c(x,y)-\phi_i(x) \}$ are continuous, we deduce \eqref{eq:E_k_O_k}. Next, $\phi_2-\phi_1=\lambda_k$ in $\operatorname{cl}(O_k)$, and \eqref{eq:E_k_O_k} yields that $\psi_2-\psi_1=-\lambda_k$ in $E_k$.

Now fix arbitrary $k,l$. Since $\operatorname{cl}(O_k),\operatorname{cl}(O_l)$ are linked, there exist $\{i_j\}_{j=1}^m$ such that $k=i_1, l=i_m$, and $E_{i_j} \cap E_{i_{j+1}} \neq \emptyset,~1\leq j \leq m$. Since $\psi_2-\psi_1=-\lambda_{i_j}$ in $E_{i_j}$, and $\psi_2-\psi_1=-\lambda_{i_{j+1}}$ in $E_{i_{j+1}}$, we obtain that $\lambda_{i_j}=\lambda_{i_{j+1}}$ for all $j$. Thus, $\lambda_k=\lambda_l$, and, consequently, $\phi_2-\phi_1=\lambda$ in $\operatorname{int}(\operatorname{supp}(\rho))=\cup_k O_k$. Finally, \eqref{eq:cl_int} and the continuity of $\phi_i$ yield that $\phi_2-\phi_1=\lambda$ in $\operatorname{supp}(\rho)$.
\end{proof}

\subsection{Proof of \Cref{prp:non_dif_non_ac}} \label{sec:non_dif_non_ac}
\begin{proof}
The proof is based on the following points.
\begin{enumerate}
    \item We have that $|\partial_\theta \rho(x,\theta)|=|\chi_{[0,1]}(x)-\chi_{[2,3]}(x)| \leq 1$, for all $x \in \Omega$.
    \item Assume that $-0.5<\theta_1<\theta_2<0.5$. In $\R$, OT maps are precisely the order-preserving ones \cite[Section 2.2]{villani03}.  The total mass of $[0,1]$ with respect to $\rho(\cdot,\theta_1)$ and $\rho(\cdot,\theta_2)$ is $0.5+\theta_1$ and $0.5+\theta_2$, respectively. Since $0.5+\theta_1<0.5+\theta_2$, all of the mass of $\rho(\cdot,\theta_1)$ from $[0,1]$ has to be transported to $[0,1]$ with a linear transport map $T(x)=\frac{0.5+\theta_1}{0.5+\theta_2}x$. Meanwhile, the excess mass of $\rho(\cdot,\theta_2)$ in $[0,1]$, supported on $\left[\frac{0.5+\theta_1}{0.5+\theta_2}, 1\right]$, has to be transported from $[2,3]$, and therefore has to travel a distance $\geq 1$. Since the excess mass of $\rho(\cdot,\theta_2)$ left in $[0,1]$ is equal to $0.5+\theta_2-(0.5+\theta_1)=\theta_2-\theta_1$, we obtain that the transport cost is at least $(\theta_2-\theta_1)\cdot 1^p$. Thus,
\[
        W_p(\rho(\cdot,\theta_1),\rho(\cdot,\theta_2)) \geq |\theta_2-\theta_1|^\frac{1}{p},\quad \forall \theta_1,\theta_2 \in (-0.5,0.5),
\]
    which means that $\theta \mapsto \rho(\cdot,\theta)$ is not absolutely continuous with respect to $W_p$ metric.
    \item Fix an arbitrary $|\theta|<0.5$. We only use the fact that $\operatorname{supp}(\rho(\cdot,\theta)) \subsetneq [0,4]$. Assume by contradiction that $\rho \mapsto W_p^p(\rho,\rho^*)$ is G\^{a}teaux differentiable at $\rho(\cdot,\theta)$ in the sense of \cite[Definition 7.12]{santambrogio15}; that is, there exists a measurable function $g$ such that
\[
        \frac{d}{d\epsilon} W_p^p(\rho(\cdot,\theta)+\epsilon (\widetilde{\rho}-\rho(\cdot,\theta)),\rho^*)\big|_{\epsilon=0+}=\int_0^4 g(x) (\widetilde{\rho}(x)-\rho(x,\theta))dx 
\]
    for all $\widetilde{\rho}\in \Pp(\Omega) \cap L^\infty(\Omega)$. Let $\phi \in C([0,4])$ be an arbitrary Kantorovich potential. From \cite[Proposition 7.17]{santambrogio15} we have that $\phi$ is in the subdifferential of $\rho \mapsto W_p^p(\rho,\rho^*)$ at $\rho(\cdot,\theta)$, and so
\[
        \frac{d}{d\epsilon} W_p^p(\rho(\cdot,\theta)+\epsilon (\widetilde{\rho}-\rho(\cdot,\theta)),\rho^*)\big|_{\epsilon=0+} \geq \int_0^4 \phi(x) (\widetilde{\rho}(x)-\rho(x,\theta))dx 
\]
    Combining this inequality with the preceding equality, we obtain
\[
        \int_0^4 (\phi(x)-g(x)) \rho(x,\theta)dx  \geq \int_0^4 (\phi(x)-g(x)) \widetilde{\rho}(x) dx 
\]
    for all $\widetilde{\rho}\in \Pp(\Omega) \cap L^\infty(\Omega)$ and Kantorovich potentials $\phi$. Fix an arbitrary potential $\phi_0$ and take $\widetilde{\rho}(x)=\chi_{(1,2)}(x)$. Furthermore, for every $\lambda \in \mathbb{R}$ consider
\[
        \phi_\lambda(x)=\phi_0(x)+\lambda (x-1)(2-x) \chi_{(1,2)}(x)
\]
    Note that $\phi_\lambda$ is continuous and $\phi_\lambda=\phi_0$ in $\operatorname{supp}\rho(\cdot,\theta)$. Thus, if $(\phi_0,\psi_0)$ is a pair of Kantorovich potentials, then $(\phi_\lambda,\psi_0)$ is also a pair of Kantorovich potentials. Plugging in $\phi=\phi_\lambda$ in the inequality above we obtain
    \begin{equation*}
    \begin{split}
        \int_0^4 (\phi_0(x)-g(x)) \rho(x,\theta)dx  \geq& \int_1^2 (\phi_0(x)-g(x)) dx + \lambda \int_1^2 (x-1)(2-x)dx\\
        =&\int_1^2 (\phi_0(x)-g(x)) dx+\frac{\lambda}{6}
    \end{split}
    \end{equation*}
    for all $\lambda \in \R$, which is a contradiction.
    \item For this and the following item, we need an explicit characterization of the OT map, $T_\theta$, from $\rho(\cdot,\theta)$ to $\rho^*$. For $\theta=0$, we have that $\rho^*$ is a translation of $\rho(\cdot,0)$. Thus, we have that
\[
    T_0(x)=x+1,\quad W_p^p(\rho(\cdot,0),\rho^*)=1.    
\]
    Next, for $\theta>0$ we have that $\rho([0,1],\theta)=0.5+\theta>0.5=\rho^*([1,2])$. Therefore,
    \begin{equation}\label{eq:T_theta>0}
        T_\theta(x)=\begin{cases}
            1+\frac{0.5+\theta}{0.5}x,\quad x\in [0,\frac{0.5}{0.5+\theta}]\\
            3+\frac{0.5+\theta}{0.5}(x-\frac{0.5}{0.5+\theta}),\quad x\in [\frac{0.5}{0.5+\theta},1]\\
            3+\frac{\theta}{0.5}+\frac{0.5-\theta}{0.5} (x-2),\quad x\in [2,3]
        \end{cases}
    \end{equation}
    For $\theta<0$ we have that $\rho([0,1],\theta)=0.5+\theta<0.5=\rho^*([1,2])$. Therefore,
    \begin{equation}\label{eq:T_theta<0}
        T_\theta(x)=\begin{cases}
            1+\frac{0.5+\theta}{0.5}x,\quad x\in [0,1]\\
            1+\frac{0.5+\theta}{0.5}+\frac{0.5-\theta}{0.5} (x-2),\quad x\in [2,2-\frac{\theta}{0.5-\theta}]\\
            3+\frac{0.5-\theta}{0.5}(x-2+\frac{\theta}{0.5-\theta}), \quad x\in [2-\frac{\theta}{0.5-\theta},3]\\
        \end{cases}
    \end{equation}
    For all $\theta$, the connected components of $\operatorname{int}(\operatorname{supp}(\rho(\cdot,\theta)))$ are
\[
        O_1=(0,1),\quad O_2=(2,3)
\]
    Furthermore, using the definition \eqref{eq:E_k} and invoking \eqref{eq:T_theta>0}, \eqref{eq:T_theta<0} we obtain
    \begin{equation*}
E_1=\begin{cases}
                [1,2],\,\theta=0\\
                [1,2] \cup [3,3+\frac{\theta}{0.5}],\, \theta>0\\
                [1,1+\frac{0.5+\theta}{0.5}],\, \theta<0
            \end{cases},\,
            E_2=\begin{cases}
                [3,4],\, \theta=0\\
                [3+\frac{\theta}{0.5},4],\,\theta>0\\
                [1+\frac{0.5+\theta}{0.5},2] \cup [3,4],\, \theta<0
            \end{cases}.
    \end{equation*}
    Thus, we have that
    \begin{equation*}
        E_1 \cap E_2=\begin{cases}
            \emptyset,\quad \theta=0\\
            \{3+\frac{\theta}{0.5}\},\quad \theta>0\\
            \{1+\frac{0.5+\theta}{0.5}\},\quad \theta<0
        \end{cases}
    \end{equation*}
    which means that $\operatorname{cl}(O_1),\operatorname{cl}(O_2)$ are linked for all $|\theta|<0.5$ except $\theta=0$.
    \item The differentiability of $\theta \mapsto W_p^p(\rho(\cdot,\theta),\rho^*)$ at $\theta\neq 0$ follows from \Cref{cor:f_Gat_suff}, and Item 4 above.
    
    Recall that $W_p^p(\rho(\cdot,0),\rho^*)=1$. Next, assume that $\theta>0$. From \eqref{eq:T_theta>0},
    \begin{equation}\label{eq:W_p^p_theta>0}
            W_p^p(\rho(\cdot,\theta),\rho^*)=\sum_{k=1}^3\int_{I_k} |T_\theta(x)-x|^p \rho(x,\theta) dx,
    \end{equation}
\[
      \text{where\quad }  I_1=[0,\frac{0.5}{0.5+\theta}],\quad I_2=[\frac{0.5}{0.5+\theta},1],\quad I_3=[2,3].
\]
    For $x\in I_1 \cup I_3$ we use the elementary inequality
\[
        |T_\theta(x)-x|^p \geq 1+p(T_\theta(x)-x-1)
\]
    For $x\in I_2$, we have that 
\[
        |T_\theta(x)-x|^p \geq 2^p
\]
    Plugging these inequalities in \eqref{eq:W_p^p_theta>0} and using \eqref{eq:T_theta>0} for evaluating elementary integrals, we obtain
\[
        W_p^p(\rho(\cdot,\theta),\rho^*)\geq 1+(2^p+p-1) \theta - p\theta^2,\quad 0<\theta < 0.5,
\]
    and so
\[
    \liminf_{\theta\to 0+} \frac{W_p^p(\rho(\cdot,\theta),\rho^*)-W_p^p(\rho(\cdot,0),\rho^*)}{\theta} \geq 2^p+p-1
\]
    For $\theta<0$, we have that
    \begin{equation}\label{eq:W_p^p_theta<0}
            W_p^p(\rho(\cdot,\theta),\rho^*)=\sum_{k=1}^3\int_{J_k} |T_\theta(x)-x|^p \rho(x,\theta) dx,
    \end{equation}
    where $J_1=[0,1]$, $J_2=[2,2-\frac{\theta}{0.5-\theta}]$,  and $J_3=[2-\frac{\theta}{0.5-\theta},3]$. Furthermore,
    \begin{equation*}
    \begin{split}
        |T_\theta(x)-x|^p \geq & 1+p(T_\theta(x)-x-1),\quad x\in J_1 \cup J_3,\\
        |T_\theta(x)-x|^p \geq & 0,\quad x\in J_2.
    \end{split}
    \end{equation*}
    Plugging these inequalities in \eqref{eq:W_p^p_theta<0}, we obtain
\[
        W_p^p(\rho(\cdot,\theta),\rho^*)\geq 1+(p+1) \theta + p\theta^2,\quad -0.5<\theta < 0,
\]
    and so
\[
    \limsup_{\theta\to 0-} \frac{W_p^p(\rho(\cdot,\theta),\rho^*)-W_p^p(\rho(\cdot,0),\rho^*)}{\theta} \leq p+1
\]
    Since $2^p+p-1>p+1$ for $p>1$, we obtain that $\theta \mapsto W_p^p(\rho(\cdot,\theta),\rho^*)$ is not differentiable at $\theta=0$. 
\end{enumerate}
\end{proof}

\subsection{Proof of~\Cref{prp:grad_error}}\label{sec:grad_error}
\begin{proof}
Assume by contradiction that there exist $(\phi_n,\psi_n) \in \Phi_c(\rho(\cdot,\theta_0),\rho^*)$ and $\epsilon_0>0$ such that $I(\phi_n,\psi_n)>f(\theta_0)-\frac{1}{n}$ and
\begin{equation}\label{eq:eps_0_error}
    \left|\nabla_\theta f(\theta_0)-\int_\Omega \phi_n^{cc}(x) \nabla_\theta \rho(x,\theta_0)dx\right|\geq \epsilon_0.
\end{equation}
Note that by adding a suitable constant to $\phi_n$, we can assume that $\sup \phi_n=\|c\|_\infty$. Thus, $(\phi_n^{cc},\phi_n^c) \in K_c$ and
\[
    f(\theta_0) \geq I(\phi_n^{cc},\phi_n^c,\theta_0) \geq I(\phi_n,\psi_n,\theta_0)>f(\theta_0)-\frac{1}{n}.
\]
Since $K_c$ is compact, we have that $(\phi_n^{cc},\phi_n^c) \to (\phi,\phi^c) \in K_c$ at least through a subsequence. Thus,
\[
    I(\phi,\phi^c,\theta_0)=\lim \limits_{n \to \infty} I(\phi_n^{cc},\phi_n^c,\theta_0)=f(\theta_0),
\]
and so $(\phi,\phi^c)\in \mathcal{S}(\theta_0)$. Hence, from \Cref{prp:dir_der_f} we have that
\begin{equation*}
\begin{split}
    &\left|\nabla_\theta f(\theta_0)-\int_\Omega \phi_n^{cc}(x) \nabla_\theta \rho(x,\theta_0)dx\right|\\
    =&\left|\int_\Omega \phi(x) \nabla_\theta \rho(x,\theta_0)dx-\int_\Omega \phi_n^{cc}(x) \nabla_\theta \rho(x,\theta_0)dx\right|\leq \|\phi-\phi_n^{cc}\|_\infty \|\eta\|_1,
\end{split}
\end{equation*}
which contradicts to \eqref{eq:eps_0_error} and finishes the proof.
\end{proof}

\section{Numerical Schemes for Computing the Gradient}\label{sec:numerical_scheme}
\subsection{Numerical Scheme for \eqref{eq:algo_discrete} and \eqref{eq:algo_adjoint_discrete}} \label{sec:numerical_scheme_1}
We remark that the first equation in both~\eqref{eq:algo_discrete} and \eqref{eq:algo_adjoint_discrete} are the same, which corresponds to solving the forward problem~\eqref{eq:main_direct} given the current iterate of the unknown parameter $\theta^l$. There are at least three ways to solve the linear system: (1) the power method, (2) the Richardson iteration, and (3) the sparse linear solve. We refer the readers to~\cite{gleich2015pagerank} for more details about the first two approaches and explain (3) in more detail.

In~\eqref{eq:main_direct}, we are interested in finding the solution $\rho$ to the linear system
\begin{equation}\label{eq:fwd_lins_sys}
    M_\epsilon \rho = (1-\epsilon) M \rho +\frac{\epsilon}{n} \mathbf{1}\mathbf{1}^\top \rho = \rho
\end{equation}
where $\mathbf{1}=[1,1,\ldots,1]^\top$, $M$ is defined in~\eqref{eq:main_direct-1} and $\epsilon$ is our teleportation (regularization) parameter. Thus, we can rewrite the linear system as
\[
     \big( (1-\epsilon) M  - I \big) \rho = - \frac{\epsilon}{n} \mathbf{1} \mathbf{1}^\top \rho.
\]
Since the biggest eigenvalue of $(1-\epsilon) M$ is $1-\epsilon < 1$, the matrix on the left-hand side is invertible, and the solution is unique. We have
\begin{equation}\label{eq:fwd_sol}
    \rho^* = \frac{\rho}{\mathbf{1}^\top \rho} =  - \big( (1-\epsilon) M  - I \big)^{-1}  \frac{\epsilon}{n} \mathbf{1}, 
\end{equation}
where $\rho^*$ is our solution that we seek as $\mathbf{1}^\top \rho^* = 1$. 

Regarding the second equation of~\eqref{eq:algo_discrete}, we consider the general linear system as below to solve for $\zeta$ given the right-hand side $y$ where
\begin{equation}\label{eq:disc_lins_sys}
     (M_{\epsilon} - I) \zeta = y.
\end{equation}
Based on~\Cref{lem:rho_diff}, we know the right-hand side of~\eqref{eq:algo_discrete}, which we denote as $y$, satisfies $y\cdot \mathbf{1} = 0$, and $M_{\epsilon} - I$ has a one-dimensional null space with generator $\rho^*$. We seek a unique solution $\zeta^*$ where $\mathbf{1}^\top \zeta^* = 0$.
Note that~\eqref{eq:disc_lins_sys} is equivalent to
\[
    \big( (1-\epsilon) M  - I \big) \zeta = y - \frac{\epsilon}{n} \mathbf{1}\mathbf{1}^\top \zeta.
\]
Since $\mathbf{1}^\top \zeta^*=0$, we obtain that $\zeta^*$ must satisfy $\big( (1-\epsilon) M  - I \big) \zeta^* = y$. As above, $(1-\epsilon) M  - I$ is invertible, and this system has a unique solution. Therefore,
\begin{equation}
    \zeta^*= \big( (1-\epsilon) M - I \big)^{-1} y.
\end{equation}

Regarding the third equation of~\eqref{eq:algo_adjoint_discrete}, we consider the general linear system as below to solve for $\zeta$ given the right-hand side $b$ where
\begin{equation}\label{eq:disc_lins_sys2}
     (M^\top_{\epsilon} - I) \zeta = b.
\end{equation}
Based on~\Cref{lem:adjoint}, we know the right-hand side of~\eqref{eq:algo_adjoint_discrete}, which we denote as $b$, satisfies $b\cdot \rho^* = 0$, and $M^\top_{\epsilon} - I$ has a one-dimensional null space with generator $\mathbf{1}$. We seek a unique solution $\zeta^*$ where $\mathbf{1}^\top \zeta^* = 0$.
Note that~\eqref{eq:disc_lins_sys2} is equivalent to
\[
    \big( (1-\epsilon) M^\top  - I \big) \zeta = b - \frac{\epsilon}{n} \mathbf{1}\mathbf{1}^\top \zeta.
\]
Since $\mathbf{1}^\top \zeta^*=0$, we obtain that $\zeta^*$ must satisfy $\big( (1-\epsilon) M^\top  - I \big) \zeta = b$. As above, $(1-\epsilon) M^\top  - I$ is invertible, and this system has a unique solution. Therefore
\begin{equation}
    \zeta^*= \big( (1-\epsilon) M^\top  - I \big)^{-1} b.
\end{equation}

Note that both the matrix $(1-\epsilon)M-I$ and its transpose are sparse. Therefore, it is relatively efficient to solve the linear systems that are essential for gradient calculation. The other components in~\eqref{eq:algo_discrete} and \eqref{eq:algo_adjoint_discrete} are rather straightforward once we solve ~\eqref{eq:fwd_lins_sys} and~\eqref{eq:disc_lins_sys} (or~\eqref{eq:disc_lins_sys2}). 

\subsection{Automatic Differentiation}\label{sec:AD}
In the previous sections, we explained how to directly compute $\nabla_{\theta} K_{mat} (\theta)$ (or equivalently $\nabla_{\theta} M_{\epsilon} (\theta)$), which is necessary to calculate $\nabla_\theta \rho(\theta)$. However, if the numerical scheme for the forward problem changes, the structure of $K_{mat} (\theta)$ changes, and consequently, one has to re-derive the explicit form of $\nabla_{\theta} K_{mat} (\theta)$. Such situations occur when using a higher-order finite volume method or switching to other standard numerical schemes such as the discontinuous Galerkin method. In order to make our code more flexible, we also implemented an automatic differentiation version using the Python library JAX \cite{jax2018github}.

Automatic differentiation techniques have been used since the 1990s for optimization, parameter identification, nonlinear equation solving, the numerical integration of differential equations, and combinations of these (see for example \cite{nocedal2006numerical, griewank2008evaluating}). However, it is only after the advent of high-speed computers and modern deep learning algorithms that automatic differentiation became extremely popular and largely used; automatic differentiation techniques made the computation of derivatives for functions defined by evaluation programs both easier and faster, especially for complicated functions with thousands of parameters like neural networks.


We compute the full Jacobian matrices of $K_{mat} (\theta)$ using the \texttt{jacfwd} function. It uses the forward-mode automatic differentiation, the most efficient choice when working with ``tall" matrices like in our case. The improvement in flexibility, however, comes with a considerable increase in computational time: the computation by pre-calculated formulae of the derivative is roughly $3$ times faster than the JAX computation; it also comes with a slight decrease in accuracy: on average, the difference between the derivative matrices computed by hand and with JAX is of the order of $10^{-14}$. This is because the code generating $K_{mat}$ contains many element-wise computations, resulting in a large computation graph for automatic differentiation. Thus, it may be preferable when working with synthetic data to use the other two approaches. 


\begin{remark}
The method of automatic differentiation is extremely valuable when working with real-world data despite the increase in computational time and the decrease in terms of accuracy. In many realistic situations, such as weather forecast, we do not have access to the underlying dynamical system, and thus we cannot compute $\nabla_{\theta} K_{mat} (\theta)$ directly. In future work, we plan on using neural networks to approximate the dynamical system from data. Given the large number of parameters and the complex functional form of a deep neural network, it would be impossible to derive $\nabla_{\theta} K_{mat} (\theta)$ explicitly, making the automatic differentiation approach necessary.
\end{remark}

\section{More Numerical Results}\label{sec:more num results}
\subsection{Single Parameter Inversion for the R\"ossler and Chen systems}\label{sec:RosslerChen}
\begin{figure}
    \centering
    \includegraphics[width = 0.3\textwidth]{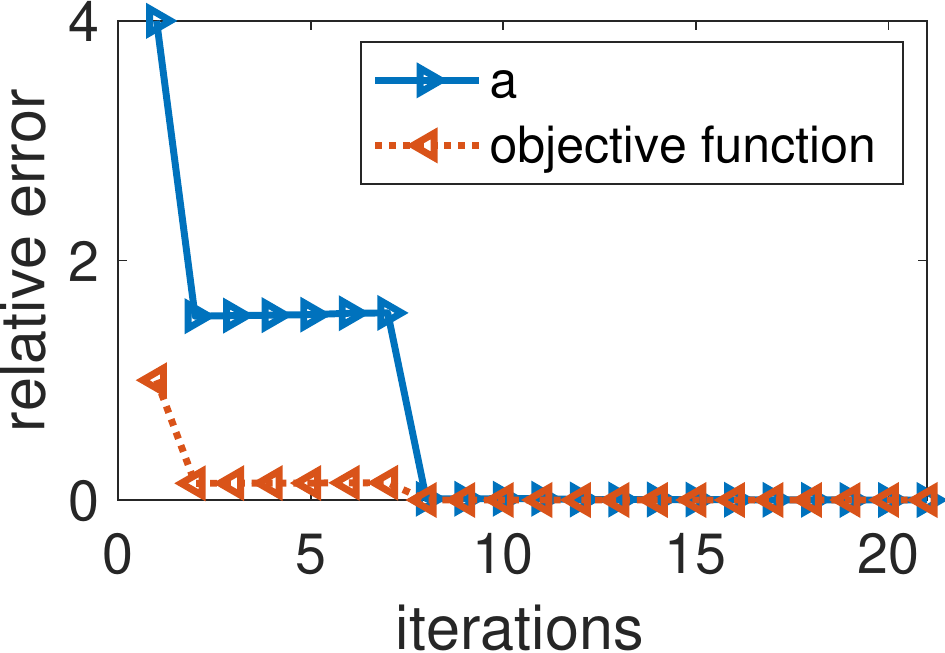}
    \includegraphics[width = 0.3\textwidth]{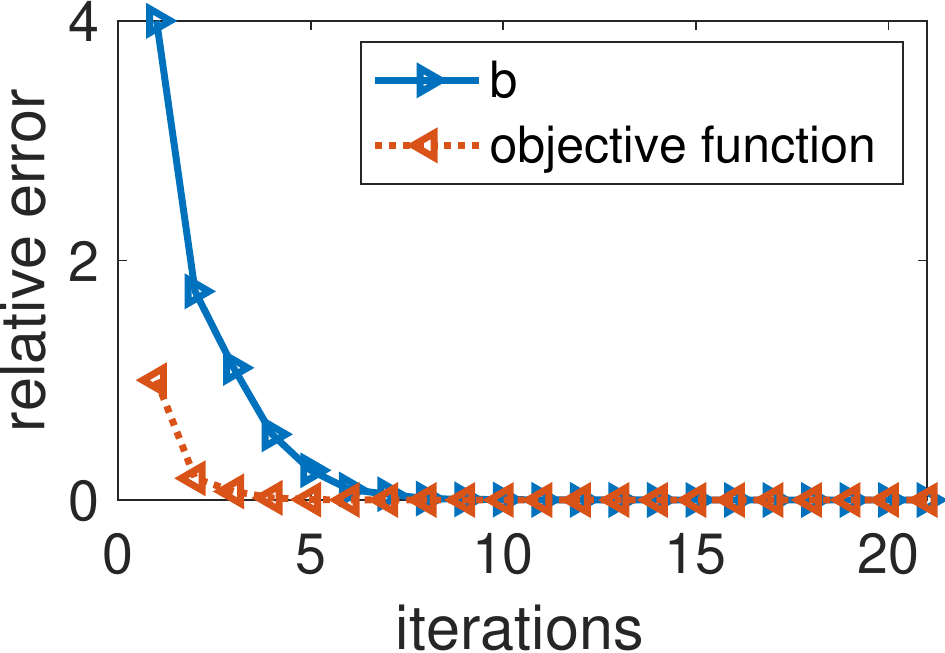}
    \includegraphics[width = 0.3\textwidth]{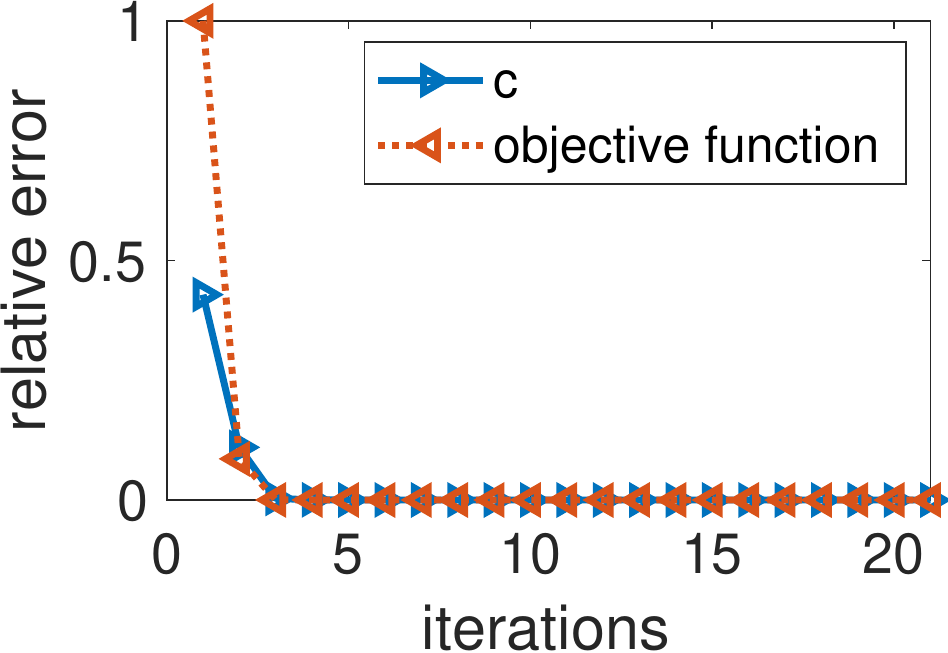}
    \caption{R\"ossler system single-parameter inference starting with $a = 0.5$ (left), $b=0.5$ (middle), $c = 10$ (right), respectively. The reference PDF is generated by the truth $(a,b,c) = (0.1,0.1,14)$ through the same numerical solver for the synthetic data.}
    \label{fig:Rossler_inverse_crime}
\end{figure}

\begin{figure}
    \centering
    \includegraphics[width = 0.3\textwidth]{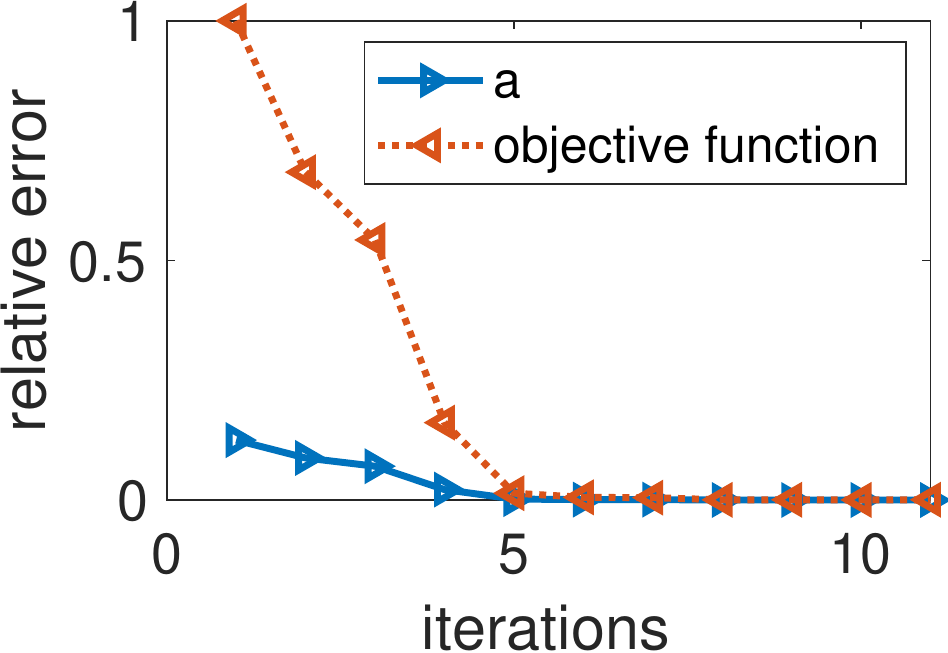}
    \includegraphics[width = 0.3\textwidth]{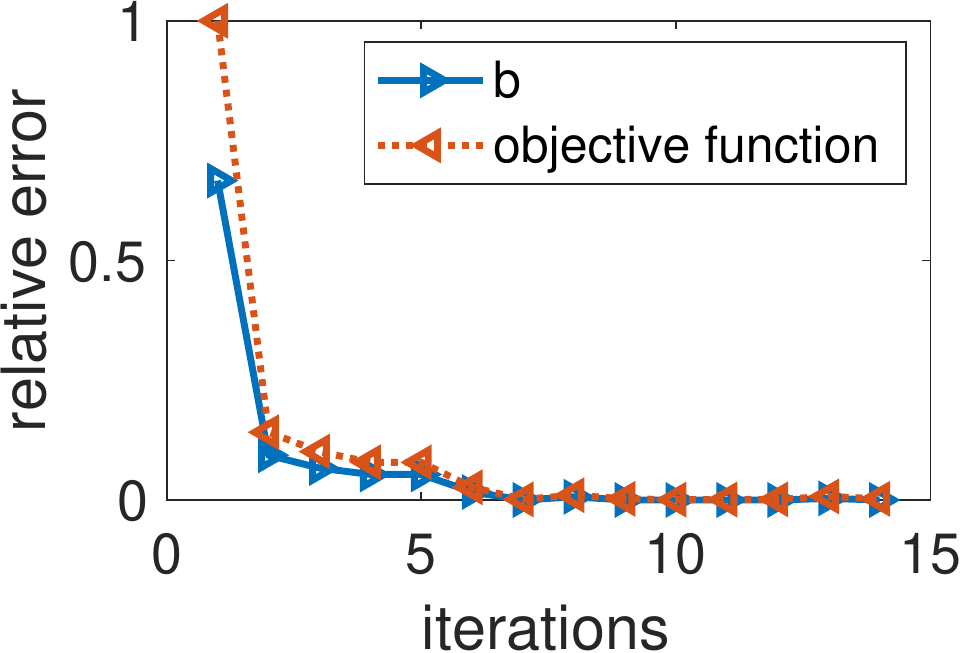}
    \includegraphics[width = 0.3\textwidth]{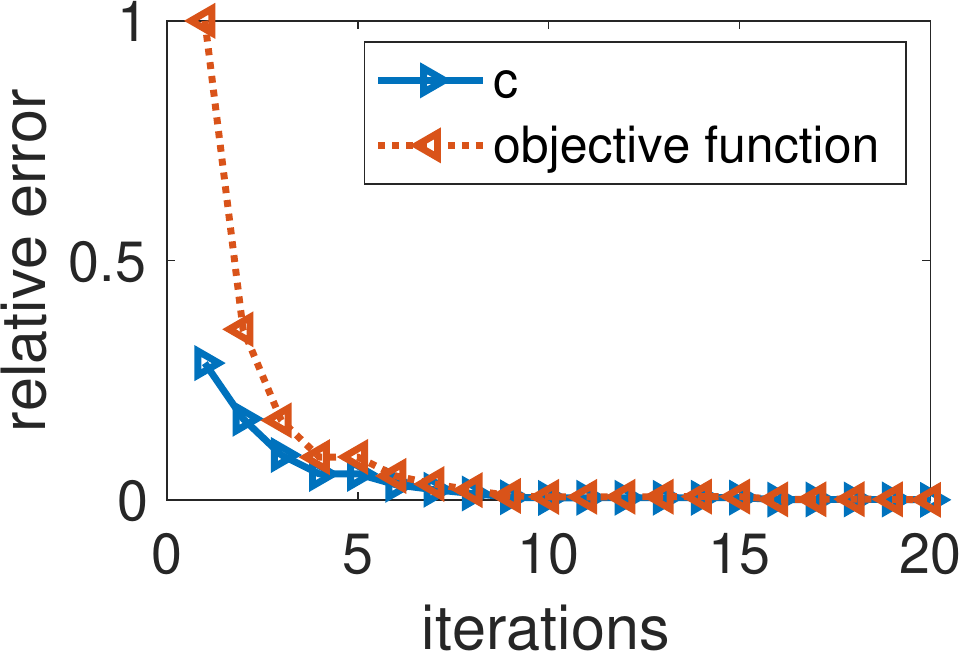}
    \caption{Chen system single-parameter inference starting with $a = 45$ (left), $b=5$ (middle), $c = 20$ (right), respectively. The reference PDF is generated by the truth $(a,b,c) = (40,3,28)$ through the same numerical solver for the synthetic data.}
    \label{fig:Chen_inverse_crime}
\end{figure}

\Cref{fig:Rossler_inverse_crime} and \Cref{fig:Chen_inverse_crime} show the single-parameter inversion for the R\"ossler and Chen systems, respectively. The reference data is produced by the same PDE solver as the synthetic data but evaluated at the true set of parameters.

\subsection{Convergence History of Parameter Inference With Noisy Time Trajectories}\label{sec:inference_with_noise}
\begin{figure}
    \centering
    \subfloat[Single-parameter inversion]{
    \includegraphics[width = 0.3\textwidth]{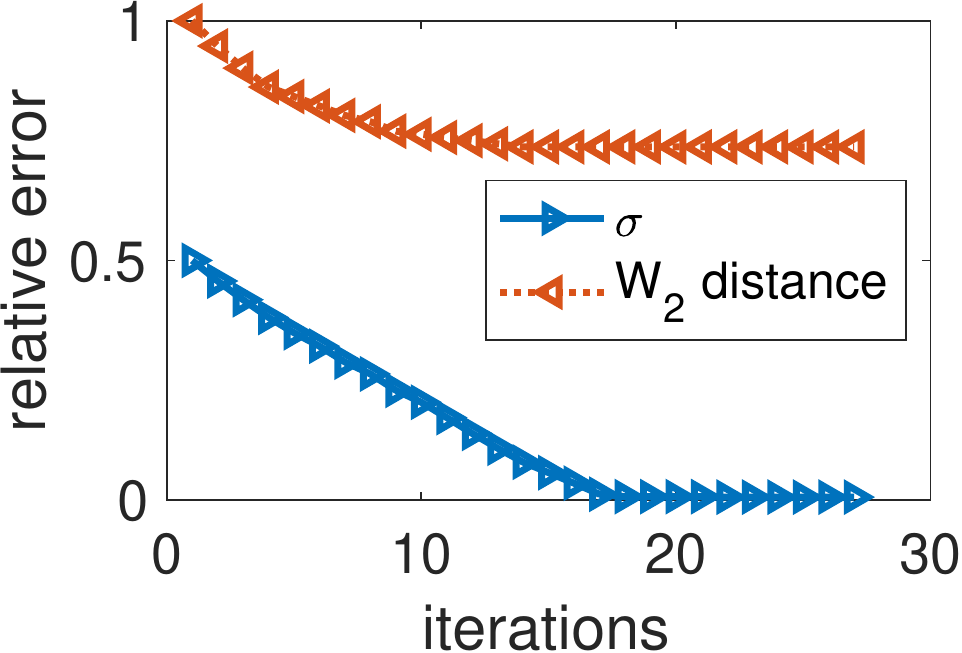}
    \includegraphics[width = 0.3\textwidth]{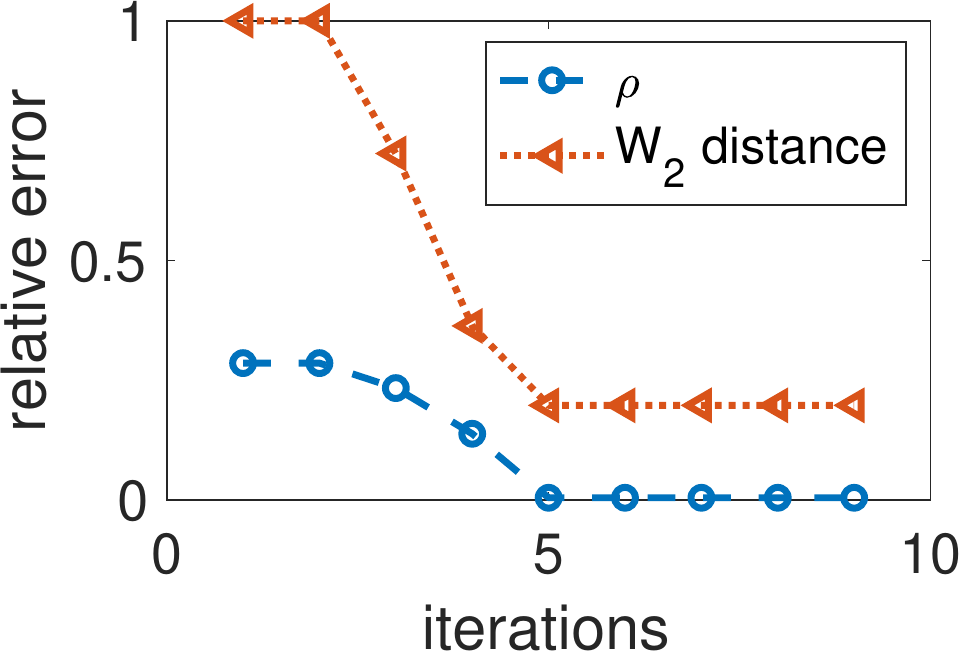}
    \includegraphics[width = 0.3\textwidth]{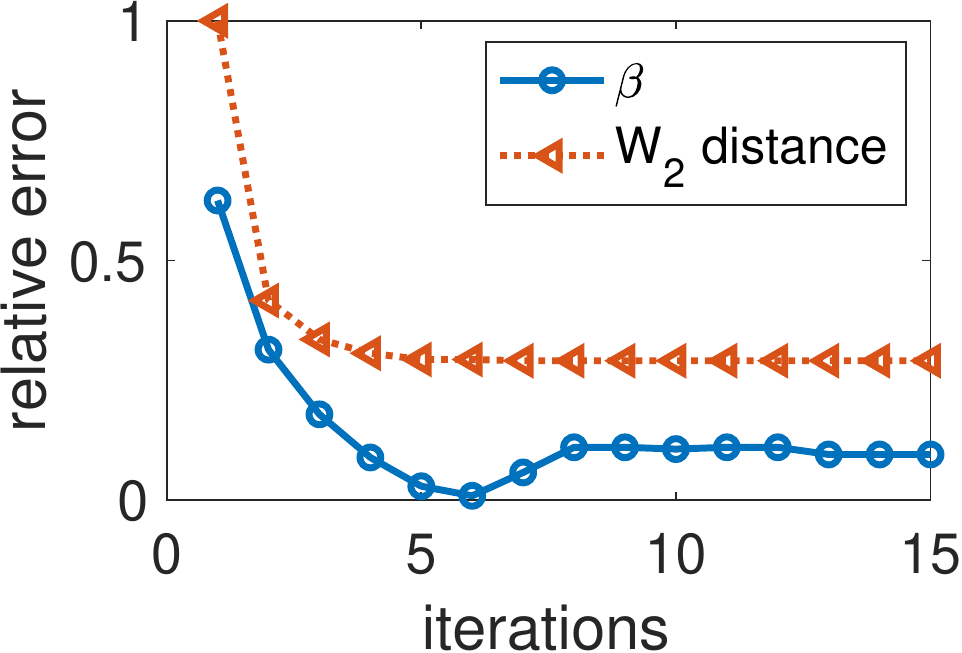} \label{fig:lorenz_intrinsic_inv_single}}\\
\subfloat[Multi-parameter inversion]{\includegraphics[width = 0.9\textwidth]{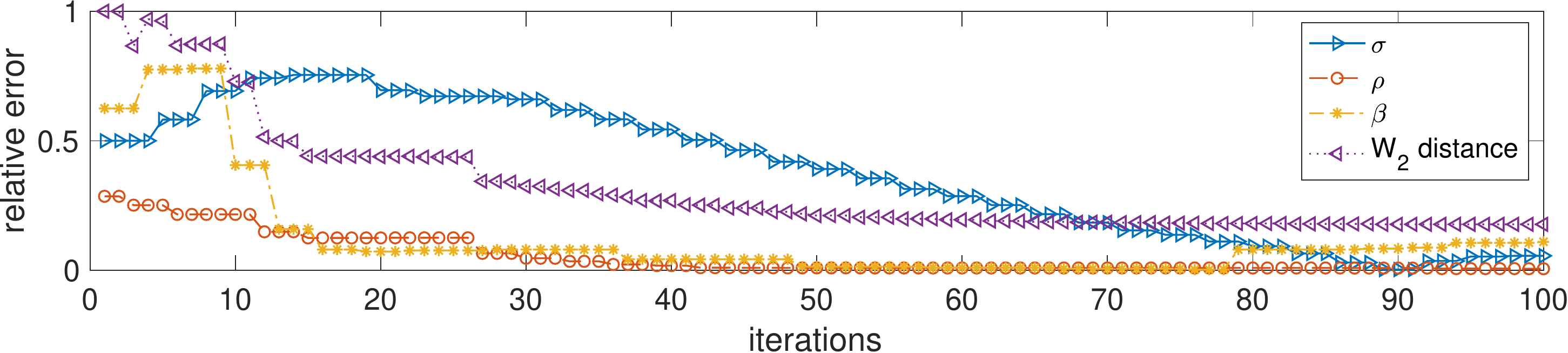}\label{fig:lorenz_intrinsic_inv_mul}}
    \caption{Top row: Lorenz system single-parameter inference starting with $\sigma = 5$ (left), $\rho=20$ (middle), $\beta=1$ (right), respectively. Bottom row: Multi-parameter inference using coordinate gradient descent with initial guess $(\sigma, \rho, \beta) = (5, 20, 1)$. The reference PDF is the histogram from the time trajectory with \textit{intrinsic} noise.}
    \label{fig:Lorenz_intrinsic_noise_W2}
\end{figure}

Fig.~\ref{fig:Lorenz_extrinsic_noise_W2} and Fig.~\ref{fig:Lorenz_intrinsic_noise_W2} are the inversion results where the Lorenz time trajectory is polluted by extrinsic and intrinsic noises, respectively. The properties of the time trajectories that are affected by the intrinsic and extrinsic noises are the same as the ones in~\Cref{fig:Lorenz_pdf}. As one can see from all the single-parameter and multi-parameter inversions, it gets more challenging to achieve reconstruction with high accuracy than the previous noise-free cases. In particular, the over-fitting phenomenon occurs, which can be directly observed for $\beta$ in the single-parameter inversion (the top right plot in both figures) and the three-parameter joint inversion (the bottom plots). As the number of iterations increases, the reconstructed $\beta$ first reaches the actual value but immediately deviates away as the objective function keeps being minimized to fit the noise. 

\end{appendix}
\end{document}